\documentclass{amsart}
\usepackage{amssymb,amscd}
\usepackage{amsmath}
\usepackage{amsthm}

\newtheorem{theorem}{Theorem}
\newtheorem{lemma}[theorem]{Lemma}

\newtheorem{proposition}[theorem]{Proposition}
\newtheorem{corollary}[theorem]{Corollary}
\numberwithin{theorem}{section}

\theoremstyle{definition}

\newtheorem{example}[theorem]{Example}

\theoremstyle{remark}
\newtheorem{remark}[theorem]{Remark}

\numberwithin{equation}{section}
\numberwithin{figure}{section}

\newfont{\germ}{eufm10}

\newcommand{\ot}{\otimes}

\newcommand{\Z}{{\mathbb Z}}
\newcommand{\R}{{\mathbb R}}
\newcommand{\J}{{\mathcal J}}
\newcommand{\Jb}{{\overline{\mathcal J}}}
\newcommand{\Rig}{{\rm Rig}}
\newcommand{\Pth}{{\mathcal P}}

\newcommand{\lc}{{\mathcal L}}
\newcommand{\Aff}{{\rm Aff}}
\newcommand{\wt}{{\rm wt}}
\newcommand{\M}{{\mathcal M}}
\newcommand{\U}{{\mathcal U}}

\def\et#1{\tilde{e}_{#1}}
\def\ft#1{\tilde{f}_{#1}}

\def\veps{\varepsilon}
\def\vphi{\varphi}
\def\ot{\otimes}

\begin{document}

\title{Bethe ansatz and inverse scattering transform \\
in a periodic box-ball system}

\author{A. Kuniba}
\address{Institute of Physics, University of Tokyo, Tokyo 153-8902, Japan}
\email{atsuo@gokutan.c.u-tokyo.ac.jp}

\author{T. Takagi}
\address{Department of Applied Physics, National Defense Academy,
Kanagawa 239-8686, Japan}
\email{takagi@nda.ac.jp}

\author{A. Takenouchi}
\address{Institute of Physics, University of Tokyo, Tokyo 153-8902, Japan}
\email{takenouchi@gokutan.c.u-tokyo.ac.jp}

\begin{abstract}
We formulate the inverse scattering method 
for a periodic box-ball system
and solve the initial value problem.
It is done by a synthesis of 
the combinatorial Bethe ans{\" a}tze 
at $q=1$ and $q=0$, which provides 
the ultradiscrete analogue of
quasi-periodic solutions in soliton equations, e.g., 
action-angle variables, Jacobi varieties, period matrices and so forth.
As an application we establish explicit formulas 
counting the states characterized by conserved 
quantities and the generic and fundamental period 
under the commuting family of time evolutions.
\end{abstract}


\maketitle

\section{Introduction}\label{sec:1}

In \cite{KT1,KT2}, a class of periodic soliton cellular automata
on one dimensional lattice is introduced.
They are associated with crystal basis of 
non-exceptional quantum affine algebras 
$U_q({\mathfrak g}_n)$ \cite{KKM}.
In this paper we focus on the simplest case 
${\mathfrak g}_n = A^{(1)}_1$ known as the 
periodic box-ball system \cite{YT},
and solve the initial value problem by the 
inverse scattering method.
As a result, all the properties and formulas conjectured in \cite{KT1,KT2} 
are established in this case.
Our approach is a synthesis of 
the two versions of the combinatorial Bethe ansatz 
at $q=1$ \cite{KKR,KR} and $q=0$ \cite{KN}.
We provide a self-contained proof 
for all the important statements together with  
several examples.
Here is a typical time evolution 
pattern on the length $13$ lattice: 
\vskip0.2cm
\noindent
\begin{center}
$t=0: \;\; 2 \;  2 \;  2 \;  1 \;  1 \;  1 \;  1 \;  2 \;  2 \;  1 \;  1 \;  2 \;  1 $\\$
 t=1: \;\;  1 \;  1 \;  1 \;  2 \;  2 \;  2 \;  1 \;  1 \;  1 \;  2 \;  2 \;  1 \;  2 $\\$
 t=2:  \;\; 2 \;  2 \;  1 \;  1 \;  1 \;  1 \;  2 \;  2 \;  2 \;  1 \;  1 \;  2 \;  1 $\\$
 t=3:  \;\; 1 \;  1 \;  2 \;  2 \;  2 \;  1 \;  1 \;  1 \;  1 \;  2 \;  2 \;  1 \;  2 $\\$
 t=4:  \;\; 2 \;  2 \;  1 \;  1 \;  1 \;  2 \;  2 \;  2 \;  1 \;  1 \;  1 \;  2 \;  1 $\\$
 t=5:  \;\; 1 \;  1 \;  2 \;  2 \;  1 \;  1 \;  1 \;  1 \;  2 \;  2 \;  2 \;  1 \;  2 $\\$
 t=6:  \;\; 2 \;  2 \;  1 \;  1 \;  2 \;  2 \;  2 \;  1 \;  1 \;  1 \;  1 \;  2 \;  1 $\\$
 t=7:  \;\; 1 \;  1 \;  2 \;  2 \;  1 \;  1 \;  1 \;  2 \;  2 \;  2 \;  1 \;  1 \;  2 $\\$
 t=8:  \;\; 2 \;  2 \;  1 \;  1 \;  2 \;  2 \;  1 \;  1 \;  1 \;  1 \;  2 \;  2 \;  1 $\\$
 t=9:  \;\; 1 \;  1 \;  2 \;  2 \;  1 \;  1 \;  2 \;  2 \;  2 \;  1 \;  1 \;  1 \;  2 $
\end{center}
\vskip0.2cm
Regarding the letter 1 as the background, one observes 
the three solitons traveling periodically to the right 
with the velocity equal to the amplitudes $3, 2$ and $1$
as long as they stay away from each other.
Under the collisions they do not smash into pieces 
nor glue together.
This is due to the underlying quantum group symmetry.
The time evolution from $t=3$ to $t=4$ states 
for example have been determined from the following diagram:
\begin{equation*}
\unitlength 0.8mm
{\small
\begin{picture}(130,20)(8,-11)
\multiput(0.8,0)(11.5,0){13}{\line(1,0){4}}
\multiput(2.8,-3)(11.5,0){13}{\line(0,1){6}}
\put(1.8,5){1}
\put(13.3,5){1}
\put(24.8,5){2}
\put(36.3,5){2}
\put(47.8,5){2}
\put(59.3,5){1}
\put(70.8,5){1}
\put(82.3,5){1}
\put(93.8,5){1}
\put(105.3,5){2}
\put(116.8,5){2}
\put(128.3,5){1}
\put(139.8,5){2}

\put(1.8,-7.5){2}
\put(13.3,-7.5){2}
\put(24.8,-7.5){1}
\put(36.3,-7.5){1}
\put(47.8,-7.5){1}
\put(59.3,-7.5){2}
\put(70.8,-7.5){2}
\put(82.3,-7.5){2}
\put(93.8,-7.5){1}
\put(105.3,-7.5){1}
\put(116.8,-7.5){1}
\put(128.3,-7.5){2}
\put(139.8,-7.5){1}

\put(-6,-1.1){122}
\put(5.5,-1.1){112}
\put(17,-1.1){111}
\put(28.5,-1.1){112}
\put(40,-1.1){122}
\put(51.5,-1.1){222}
\put(63,-1.1){122}
\put(74.5,-1.1){112}
\put(86,-1.1){111}
\put(97.5,-1.1){111}
\put(109,-1.1){112}
\put(120.5,-1.1){122}
\put(132,-1.1){112}
\put(143.5,-1.1){122}

\end{picture}
}
\end{equation*}
Here the local vertices stand for 
the $U_q(A^{(1)}_1)$ quantum $R$ matrix 
sending the spin  $\frac{3}{2} \otimes \frac{1}{2}$ states 
on the NW edges to the spin $\frac{1}{2} \otimes \frac{3}{2}$ 
states on the SE edges.
This is a standard diagram for the row transfer matrix
in the $3 \times 1$ fusion vertex model.
The peculiar feature here is the specialization to $q=0$,
which makes the $R$ a deterministic map called 
the combinatorial $R$ \cite{KMN,NY}.
The above ones correspond to the $l=3$ case in Figure \ref{fig:cr},
and the associated time evolution is named $T_3$.
On the horizontal edge, the leftmost and the rightmost 
ones are the same element $122$ reflecting the 
periodic boundary condition.
Denote the array on the top and the bottom lines by 
$p$ and $p'$, respectively.
It turns out that $p'$ is uniquely determined from $p$
by the combinatorial $R$ and the boundary condition 
said above. The resulting map defines a time evolution 
$p' = T_3(p)$.
The periodic box-ball system is a dynamical system 
endowed with such time evolutions
$T_1, T_2, \ldots$, which are essentially the 
commuting family of fusion transfer matrices at $q=0$.
In other words, it is a solvable vertex model \cite{Ba} 
at $q=0$ on the periodic lattice.
It is the periodic extension of the original 
box-ball system on the infinite lattice \cite{TS}.

In this paper we solve the initial value problem of the 
periodic box-ball system.
It is done by a synthesis of the 
Bethe ans{\" a}tze \cite{Be} at $q=1$ \cite{KKR,KR} and $q=0$ \cite{KN}, 
which yields a periodic ultradiscretization of 
the inverse scattering method in soliton theory \cite{GGKM,AS}.
As an application, 
we are able to answer the interesting question; 
what is the fundamental period of $p$ 
under the time evolutions $T_1, T_2, \ldots$? 
Namely, the minimum positive integer $N_l$
satisfying $T^{N_l}_l(p)=p$.  For instance for the above $p$, 
one has $N_1\!=\!13, N_2\!=\!91, N_3\!=\!273$\footnote{
It is a general feature of the system that 
$T_l(p)=T_\infty(p)$ hence $N_l = N_\infty$
for $l \ge l_0$, where $l_0$ is the 
amplitude of the largest soliton involved in $p$.
See Proposition \ref{pr:ts} and the subsequent remark.}.
An explicit formula for 
the fundamental period under $T_\infty$ 
has been discovered in \cite{YYT}
by a combinatorial argument.
Here we take the inverse scattering approach and 
establish the results 
(\ref{eq:lcm}) and (\ref{eq:nn}) for general $T_l$.
Under our scheme, the period ${\mathcal N}$ acquires the 
intrinsic characterization ${\mathcal N}\vec{h} \in \Gamma$ 
(\ref{eq:ncond}) 
in terms of straight motions on the set
\begin{equation*}
({\mathcal I}_{m_{j_1}}\times 
\cdots \times {\mathcal I}_{m_{j_s}})/\Gamma.
\end{equation*}
This is an `ultradiscrete Jacobi variety'. See (\ref{eq:jv}).
The conjectural formula for the generic 
period in the most general $A^{(1)}_n$ case \cite{KT2} also 
admits a similar interpretation. 

Let us explain the content of the paper in some detail.
In Section \ref{sec:2}, we begin with the rudiments of 
the crystal basis theory, 
formulate the periodic box-ball system and 
explain its basic properties.
The fundamental role is played 
by the combinatorial $R$.
It governs the local dynamics of the system 
and ultimately its whole aspect
together with the global constraint 
imposed by the periodic boundary condition.
The commuting family of time evolutions $T_l$ 
and the conserved quantity $E_l$ called energy are constructed.
It is done by a simple extension of the arguments in 
\cite{HHIKTT,FOY,HKOTY} to periodic versions and  
supplements the original description in \cite{YT}.
The both $T_l$ and the $E_l$ enjoy the symmetry 
under the extended affine Weyl group 
${\widetilde W}(A^{(1)}_1)$ (Proposition \ref{pr:w}).
The time evolution $T_l$ with sufficiently large $l$ 
admits a simple description as a translation in 
${\widetilde W}(A^{(1)}_1)$ (Proposition \ref{pr:ts}).
Although quite straightforward, the proof of these properties 
announced in \cite{KT1,KT2} are 
presented here for the first time.

Section \ref{sec:ism} is a main part of the paper
where we present our inverse scattering formalism. 
We invoke the Kerov-Kirillov-Reshetikhin (KKR) bijection 
between the rigged configurations and highest paths.
We call it the KKR theory or 
combinatorial Bethe ansatz at $q=1$ in this paper.
The necessary facts are available in Appendix \ref{app:kkr}
as well as in the original papers \cite{KKR,KR,KSS} 
and the latest review \cite{Schi}.

Roughly speaking, our action-angle variables are 
the rigged configurations.
However, the KKR theory works only for highest paths
whereas the periodic boundary condition brings 
all the paths into the game. 
To reconcile them in a remarkable harmony 
via the special prescription (\ref{eq:j}) is 
the heart of our inverse scattering formalism.
Proposition \ref{pr:migi} is the crux 
to guarantee the well-definedness of the 
resulting direct and inverse scattering map
$\Phi$ (\ref{eq:pj}).
Our solution of the initial value problem is 
presented in Theorem \ref{th:main}.
It is also applicable, with a drastic simplification, 
to the original box-ball system \cite{TS} by 
taking the system size infinite.
  
In Section \ref{sec:ba} we explain the origin 
of our inverse scattering formalism 
in the other framework of the 
combinatorial Bethe ansatz at $q=0$ \cite{KN}.
The central role is played by the linear congruence 
equation called the string center equation (\ref{eq:sce}).
It is the Bethe equation on the string centers at $q=0$
whose off-diagonal solutions yield 
the weight multiplicities of the 
$sl_2$ module $({\mathbb C}^2)^{\otimes L}$.
By the map $\Psi$ (\ref{eq:psi}), we 
link the logarithmic branch 
in the string center equation with the angle variable in 
Section \ref{sec:ism}.
Then our key construction (\ref{eq:j})
is nothing but the equivalence relation 
among the off-diagonal solutions to the 
string center equation (Theorem \ref{th:ju}).
Through our direct and inverse scattering transforms,
the nonlinear dynamics in 
the periodic box-ball system becomes
a straight motion of the Bethe roots or angle variables.
With these features (Corollary \ref{cor:pju}) in hand,
it is straightforward to derive explicit formulas for  
the generic (\ref{eq:lcm}) and 
the fundamental period (\ref{eq:nn}) under 
any time evolution $T_l$.
The number of states characterized by conserved 
quantities (\ref{eq:cpju}) are obtained and 
the Bethe eigenvalue is shown to be a root of unity 
related to the generic period (Proposition \ref{pr:la}).
The number of disjoint orbits under the commuting 
family of time evolutions is also determined in (\ref{eq:cc}).
Our angle variables 
live in the set (\ref{eq:jv}), which serves as an 
ultradiscrete analog of the Jacobi variety
in the classical theory of quasi-periodic solutions
to soliton equations \cite{DT,DMN}.

The results in 
Sections \ref{sec:ism} and \ref{sec:ba} uncover 
a significant interplay between 
the KKR theory and the string center equation.
They substantially achieve a synthesis of the 
combinatorial Bethe ans{\" a}tze at $q=1$ and $q=0$, 
which was already foreseen in the earlier works \cite{KT1,KT2}
on generalized periodic box-ball systems. 
Section \ref{sec:sum} is a summary.

Appendix \ref{app:kkr} summarizes necessary 
facts on the rigged configurations and their bijective 
correspondence with the highest paths. 
A piecewise linear formula for the KKR bijection $\phi^{-1}$
is available in (\ref{eq:pwf}).

Appendix \ref{app:em} contains the proof of 
Proposition \ref{pr:em}. 
For any path not necessarily highest, 
it essentially identifies 
the two descriptions of the conserved quantity 
in terms of the energy and the configuration.
The former is related to the soliton content and the latter 
to the string content. 

Appendix \ref{app:proof} is devoted to a proof of
Proposition \ref{pr:migi}, which assures the well-definedness 
of the direct/inverse scattering map $\Phi$ (\ref{eq:pj}).
The key is to investigate the relation between
the two rigged configurations 
corresponding to the paths 
$q\otimes r$ and $r\otimes q$ with 
$q$ and $r$ both being highest.

Appendix \ref{app:main}
is the proof of Theorem \ref{th:main},
which asserts the linearization of the time 
evolution in terms of the angle variables.
Our main idea is to realize non highest periodic paths
as a segment of a large highest path having the structure 
$p\otimes p \otimes p \otimes \cdots$ 
and apply the KKR theory to the latter. 

\section{Periodic box-ball system}\label{sec:2}

\subsection{Crystals and combinatorial $R$}\label{subsec:ccr}
We recapitulate the basic facts in the crystal basis theory 
\cite{K,KMN,KKM,NY,HKOTY}.
Let $B_l$ the the crystal of the $l$-fold symmetric 
tensor representation of $U_q(A^{(1)}_1)$.
As a set it is given by 
$B_l = \{x=(x_1,x_2) \in (\Z_{\ge 0})^2 \mid x_1+x_2=l\}$.
The element $(x_1,x_2)$ will also be expressed as 
the length $l$ row shape semistandard tableau 
containing the letter $i$ $x_i$ times.
For example,
$B_1=\big\{\fbox{1},\fbox{2}\big\}, 
B_2 = \big\{\fbox{11}, \fbox{12}, \fbox{22} \big\}$.
(We omit the frames of the tableaux hereafter.)
The action of Kashiwara operators 
${\tilde f}_i, {\tilde e}_i: B \rightarrow B \sqcup \{0\}\, (i=0,1)$ reads 
$({\tilde f}_ix)_j = x_j-\delta_{j,i}+\delta_{j,i+1}$ and 
$({\tilde e}_ix)_j = x_j+\delta_{j,i}-\delta_{j,i+1}$,
where all the indices are in $\Z_2$, and if the result does not 
belong to $(\Z_{\ge 0})^2$, it should be understood as $0$.
The classical part of the weight of $x=(x_1,x_2) \in B_l$ is 
$\wt(x)=l\Lambda_1-x_2\alpha_1 = (x_1-x_2)\Lambda_1$,
where $\Lambda_1$ and $\alpha_1=2\Lambda_1$ are 
the fundamental weight and the simple root of $A_1$.

For any $b\in B$, set 
\[
\veps_i(b)=\max\{m\ge0\mid \et{i}^m b\ne0\},\quad
\vphi_i(b)=\max\{m\ge0\mid \ft{i}^m b\ne0\}.
\]
By the definition one has $\veps_i(x)=x_{i+1}$ and $\vphi_i(x)=x_i$
for $x=(x_1,x_2) \in B_l$.

For two crystals $B$ and $B'$, one can define the tensor product
$B\ot B'=\{b\ot b'\mid b\in B,b'\in B'\}$. 
The operators $\et{i},\ft{i}$ act 
on $B\ot B'$ by
\begin{eqnarray}
\et{i}(b\ot b')&=&\left\{
\begin{array}{ll}
\et{i} b\ot b'&\mbox{ if }\vphi_i(b)\ge\veps_i(b')\\
b\ot \et{i} b'&\mbox{ if }\vphi_i(b) < \veps_i(b'),
\end{array}\right. \label{eq:erule}\\
\ft{i}(b\ot b')&=&\left\{
\begin{array}{ll}
\ft{i} b\ot b'&\mbox{ if }\vphi_i(b) > \veps_i(b')\\
b\ot \ft{i} b'&\mbox{ if }\vphi_i(b)\le\veps_i(b').
\end{array}\right. \label{eq:frule}
\end{eqnarray}
Here $0\ot b'$ and $b\ot 0$ should be understood as $0$. 
The tensor product $B_{l_1}\otimes \cdots \otimes B_{l_k}$
is obtained by repeating the above rule.
The classical part of the weight of $b\in B$ for any 
$B=B_{l_1}\otimes \cdots \otimes B_{l_k}$ is given by 
$\wt(b) = (\varphi_1(b)-\varepsilon_1(b))\Lambda_1
=(\varepsilon_0(b)-\varphi_0(b))\Lambda_1$.

The crystal $B_l$ admits the affinization $\Aff(B_l)$.
It is the infinite set 
$\Aff(B_l) = \{ \zeta^d b \mid b \in B_l, d \in \Z\}$
endowed with the crystal structure
${\tilde e}_i(\zeta^db) = \zeta^{d-\delta_{i,0}}({\tilde e}_ib)$,
${\tilde f}_i(\zeta^db) = \zeta^{d+\delta_{i,0}}({\tilde f}_ib)$.
The parameter $\zeta$ is called the spectral parameter.
The isomorphism of the affine crystal 
$\Aff(B_l) \otimes \Aff(B_k) \overset{\sim}{\rightarrow}
\Aff(B_k) \otimes \Aff(B_l)$ is the unique bijection that  
commutes with Kashiwara operators 
(up to a constant shift of $H$ below).
It is the $q=0$ analogue of the quantum $R$ and
called the combinatorial $R$.
Explicitly it is given by
$R: \zeta^d x\otimes \zeta^e y\mapsto
\zeta^{e+H(x\ot y)}\tilde{y}\ot \zeta^{d-H(x\ot y)}\tilde{x}$ with 
\begin{align}
&{\tilde x}_i = x_i+Q_i(x,y)-Q_{i-1}(x,y),\quad 
{\tilde y}_i = y_i+Q_{i-1}(x,y)-Q_i(x,y),\nonumber\\
&Q_i(x,y) = \min(x_{i+1}, y_i), \nonumber\\
&H(x\ot y) = -Q_0(x,y). \label{eq:h}
\end{align}
Here $x \otimes y \simeq {\tilde y}\otimes {\tilde x}$ 
under the isomorphism 
$B_l \otimes B_k \simeq B_k \otimes B_l$. 
The relation is depicted as

\begin{equation*}
\begin{picture}(90,40)(10,-9)
\put(0,10){\line(1,0){20}}\put(-6,8){$x$}\put(22,8){${\tilde x}$}
\put(10,0){\line(0,1){20}}\put(8,25){$y$}\put(8,-9){${\tilde y}$}
\put(47,8){or}
\put(80,10){\line(1,0){20}}\put(74,8){${\tilde y}$}\put(102,8){$y$}
\put(90,0){\line(0,1){20}}\put(88,25){${\tilde x}$}\put(88,-9){$x$}
\end{picture}
\end{equation*}

For example $B_l \otimes B_1 \simeq B_1 \otimes B_l$ 
is listed as follows:

\begin{figure}[h]\label{fig:cr}
\begin{picture}(130,120)(55,0)
\unitlength 1mm
\multiput(0,0)(0,25){2}{
\multiput(2,10)(65,0){2}{\line(1,0){16}}
\multiput(10,6)(65,0){2}{\line(0,1){8}}}

\put(9.3,40){$\scriptstyle 1$}
\put(-15,34){$\scriptstyle \overbrace{1 \;\cdots \cdots \; 1}^l$} 
\put(20,34){$\scriptstyle \overbrace{1 \;\cdots \cdots \; 1}^l$}
\put(9.3,28){$\scriptstyle 1$}
\put(7,24){$\scriptstyle H=0$}

\put(74.3,40){$\scriptstyle 2$}
\put(50,34){$\scriptstyle \overbrace{2 \;\cdots \cdots \; 2}^l$} 
\put(85,34){$\scriptstyle \overbrace{2 \;\cdots \cdots \; 2}^l$}
\put(74.3,28){$\scriptstyle 2$}
\put(72,24){$\scriptstyle H=0$}

\put(9.3,15){$\scriptstyle 1$}
\put(-18,9){$\scriptstyle \overbrace{1 \cdots 1}^{l-a}
\overbrace{2\cdots 2}^a$} 
\put(20,9){$\scriptstyle \overbrace{1 \cdots 1}^{l-a+1}
\overbrace{2\cdots 2}^{a-1}$}
\put(9.3,3){$\scriptstyle 2$}
\put(4,-1){$\scriptstyle H=0\; (0<a\le l)$}

\put(74.3,15){$\scriptstyle 2$}
\put(47.2,9){$\scriptstyle \overbrace{1 \cdots 1}^{l-a}
\overbrace{2\cdots 2}^a$} 
\put(85,9){$\scriptstyle \overbrace{1 \cdots 1}^{l-a-1}
\overbrace{2\cdots 2}^{a+1}$}
\put(74.3,3){$\scriptstyle 1$}
\put(69,-1){$\scriptstyle H=-1\;(0\le a<l)$}

\end{picture}
\caption{Combinatorial $R: B_l \otimes B_1 \simeq B_1 \otimes B_l$}
\end{figure}

Let $\omega: (x_1,x_2) \mapsto (x_2,x_1)$ be 
the involutive Dynkin digram automorphism of $B_l$.
The combinatorial $R$ enjoys the symmetry:
\begin{equation}\label{eq:dynkin}
(\omega \otimes \omega)R = R(\omega \otimes \omega)
\quad \hbox{on } B_l \otimes B_k.
\end{equation}
We write the highest element $(l,0) \in B_l$ as $u_l$.
The energy function $H$ (\ref{eq:h}) is 
normalized so as to attain the maximum at 
$H(u_l\otimes u_k) = 0$ and ranges over 
$-\min(l,k)\le H \le 0$ on $B_l \otimes B_k$.
The combinatorial $R$ satisfies the Yang-Baxter relation:
\begin{equation}\label{eq:ybr}
(1 \otimes R)(R \otimes 1)(1 \otimes R) =
(R \otimes 1)(1 \otimes R)(R \otimes 1) 
\end{equation}
on $\Aff(B_j) \otimes \Aff(B_l) \otimes \Aff(B_k)$.
Let $\varrho(b_1 \otimes \cdots \otimes b_k) 
= b_k \otimes \cdots \otimes b_1$ be the reverse ordering 
of the tensor product for any $k$.
The combinatorial $R$ has the property:
\begin{equation}\label{eq:rev} 
R \,\varrho = \varrho\,R \quad\hbox{on } B_l \otimes B_k.
\end{equation}

\subsection{Time evolution and conserved quantities}\label{subsec:te}
We fix the integer $L \in \Z_{\ge 1}$ 
corresponding to the system size throughout.
Set 
\begin{equation}\label{eq:pdef}
\Pth = B_1^{\otimes L},\quad 
\Pth_+ = \{ p \in \Pth \mid {\tilde e}_1p = 0\}.
\end{equation}
We will also write $\Aff(\Pth) = \Aff(B_1)^{\otimes L}$.
An element of $\Pth$ ($\Pth_+$) is called 
a path (highest path).
The condition ${\tilde e}_1p = 0$ on the path
$p=b_1 \otimes \cdots \otimes b_L$ is equivalent to the 
simple postulate 
\begin{equation}\label{eq:lp}
\sharp\{1 \le i \le k \mid b_i=1\} \ge 
\sharp\{1 \le i \le k \mid b_i=2\}
\quad \hbox{for all }\; 1 \le k \le L.
\end{equation}
The weight of the path $p=b_1\otimes \cdots \otimes b_L$ is 
given by 
$\wt(p) = \wt(b_1) + \cdots + \wt(b_L)$.
We write 
$\wt(p)>0$ ($\wt(p)<0)$ 
when it belongs to 
$\Z_{>0}\Lambda_1$ $(\Z_{<0}\Lambda_1)$.

The periodic box-ball system is 
a dynamical system on $\Pth$ equipped with the commuting family of 
time evolutions $T_1, T_2, \ldots$, which we shall now introduce.

\begin{proposition}\label{pr:tl}
For any path $p=b_1\otimes \cdots \otimes b_L \in \Pth$ 
and $l \in \Z_{\ge 1}$, there exists an element 
$v_l \in B_l$ such that 
($e = -d_1-\cdots - d_L$)
\begin{equation}\label{eq:bpe}
\zeta^0v_l \otimes (\zeta^0b_1 \otimes \cdots \otimes \zeta^0b_L) 
\simeq 
(\zeta^{d_1}b'_1 \otimes \cdots \otimes \zeta^{d_L}b'_L) \otimes 
\zeta^e v_l
\end{equation}
for some $\zeta^{d_1}b'_1 \otimes \cdots \otimes \zeta^{d_L}b'_L
\in \Aff(\Pth)$ under the isomorphism 
$\Aff(B_l)\otimes \Aff(\Pth) \simeq 
\Aff(\Pth) \otimes \Aff(B_l)$. 
Such $v_l$ is unique except $\wt(p)=0$ case, where 
$\zeta^{d_1}b'_1 \otimes \cdots \otimes \zeta^{d_L}b'_L$ 
is independent of the possibly non unique choice of $v_l$.
\end{proposition}

\begin{proof}
Suppose the relation 
$(x_1,x_2) \otimes p \simeq p' \otimes (y_1,y_2)$ holds 
under the isomorphism $B_l \otimes \Pth \simeq \Pth \otimes B_l$.
Setting $p'=b'_1\otimes \cdots \otimes b'_L$, we depict it 
as Figure \ref{fig:xp}.

\begin{figure}[h]\label{fig:xp}
\begin{picture}(50,50)(40,5)
\unitlength 1mm
\put(0,10){\line(1,0){20}}\put(33,10){\line(1,0){10}}
\put(24.5,9){$\cdots$}
\put(-20,9){$\scriptstyle \overbrace{1 \cdots 1}^{x_1}
\overbrace{2\cdots 2}^{x_2}$}
\put(45,9){$\scriptstyle \overbrace{1 \cdots 1}^{y_1}
\overbrace{2\cdots 2}^{y_2}$}
\put(5,6){\line(0,1){8}} 
\put(4,16){$\scriptstyle b_1$}\put(4,2.2){$\scriptstyle b'_1$}
\put(13,6){\line(0,1){8}}
\put(12,16){$\scriptstyle b_2$}\put(12,2.2){$\scriptstyle b'_2$}
\put(38,6){\line(0,1){8}}
\put(37,16){$\scriptstyle b_L$}\put(37,2.2){$\scriptstyle b'_L$}
\end{picture}
\caption{}
\end{figure}

\noindent
Fixing $p$ and $l$, we regard it as the functional relation 
$y_2=y_2(x_2)$.
Then from Figure \ref{fig:cr} we see that 
$y_2(x_2+1)=y_2(x_2)$ or $y_2(x_2+1)=y_2(x_2)+1$.
Moreover, the latter holds only if 
all the intermediate vertices in Figure \ref{fig:xp}
are the bottom two types in Figure \ref{fig:cr}, namely 
$\wt(p)=-\wt(p')$.

Let $M$ be the number of $2 \in B_1$ contained in $p$.
We first consider the case $L-M>M$, i.e., $\wt(p)>0$.
Set $c=y_2(0)$.
Suppose that $y_2(a)=c$ and $y_2(a+1)=c+1$ hold.
{}From the above observation, the number of $2$ contained in 
$p'$ (for $(x_1,x_2)=(l-a,a)$ in Figure \ref{fig:xp}) 
must be $L-M$. Then the weight conservation demands that
$a+M=c+L-M$. By the assumption $M<L-M$,
this can not happen for $0 \le a \le c$.
This implies that $y_2(a)=c$ for $0 \le a \le \min(c+1,l)$,
which confirms the sought assertion.
The case $L-M<M$ can be shown similarly by
interchanging the role of the letters $1$ and $2$.
Next we consider the case $L-M=M$, i.e., $\wt(p)=0$.
By the same argument as above, we find that 
$y_2(a)=c$ for $0 \le a \le c$ and 
$y_2(a)=a$ can happen for some interval $c<a \le c' (\le l)$.
When $c\le a < c'$, all the intermediate vertices 
in Figure \ref{fig:xp} are the bottom two types in Figure \ref{fig:cr},
whose energy $H$ and the vertical $B_1$ part are independent 
of $a$. This verifies the existence of $v_l$ and uniqueness of 
$\zeta^{d_1}b'_1 \otimes \cdots \otimes \zeta^{d_L}b'_L$.
\end{proof}

We define the time evolution $T_l$ and the energy $E_l$ of 
a path $p$ by $T_l(p) 
= b'_1 \otimes \cdots \otimes b'_L\,(\in \Pth)$ and 
$E_l(p) = e \,(\in \Z_{\ge 0})$ using the notation in (\ref{eq:bpe}).
Proposition \ref{pr:tl} assures that they are solely 
determined from $p$ and $l$.
The definition is summarized by the relation
\begin{equation}\label{eq:tle}
\zeta^0 v_l \otimes p \simeq T_l(p)\otimes \zeta^{E_l(p)}v_l 
\end{equation}
omitting the spectral parameters attached to $p$ and $T_l(p)$. 
Here $v_l \in B_l$ can be constructed from $p \in \Pth$ by 
\begin{equation}\label{eq:vt}
\begin{split}
u_l \otimes p \simeq p^\ast &\otimes v_l 
\quad \hbox{ if } \wt(p)\ge 0,\\
\omega(u_l) \otimes p \simeq p^\ast &\otimes v_l 
\quad \hbox{ if } \wt(p) < 0
\end{split}
\end{equation}
for some $p^\ast \in \Pth$ 
under the isomorphism 
$B_l \otimes \Pth \simeq \Pth \otimes B_l$,
where $u_l=(l,0) \in B_l$ as defined after (\ref{eq:dynkin}).
One may either use the latter relation to define $v_l$ 
when $\wt(p)=0$.
Clearly the time evolutions are 
weight preserving, i.e., $\wt(T_l(p)) = \wt(p)$. 
They are all invertible. By using (\ref{eq:rev}),
the inverse can be found by 
\begin{equation}\label{eq:tinv}
T^{-1}_l(p) = \varrho\, T_l(\varrho(p)).
\end{equation}
In the rest of the paper we will use, 
often without explicitly mentioning,  
the fact that 
(\ref{eq:tle}) determines $T_l(p)$ and $E_l(p)$ 
unambiguously even though $v_l$ is not unique in general.
If $p$ is a highest path, 
the affinization of (\ref{eq:vt}) is available in Lemma \ref{lem:evl}.

Note from Figure \ref{fig:cr} that the combinatorial $R$
on $B_1 \otimes B_1$ is the identity map.
It follows that $T_1$ acts as the cyclic shift or 
`$\exp\big(\sqrt{-1}\,(\hbox{momentum})$\big)':
\begin{equation}\label{eq:t1}
T_1(b_1 \otimes b_2 \otimes \cdots \otimes b_L)
= b_L \otimes b_1 \otimes \cdots \otimes b_{L-1}.
\end{equation}
$T_1$ will play a special role in our 
inverse scattering formalism in Section \ref{sec:ism}.

\begin{theorem}\label{th:te}
The commutativity $T_lT_k(p) = T_kT_l(p)$
and the conservation $E_l(T_k(p))=E_l(p)$ hold.
\end{theorem}

The commutativity with $T_1$ (\ref{eq:t1}) 
is the origin of the 
adjective ``periodic".

\begin{proof}
Take $v_k$ for $p$ and $v_l$ for $T_k(p)$ as in (\ref{eq:vt}).
Set $R(\zeta^0v_l \otimes \zeta^0 v_k) = 
\zeta^{\delta}\overline{v}_k \otimes 
\zeta^{-\delta}\overline{v}_l$ and 
regard $p$ as an element of $\Aff(\Pth)$.
By using the combinatorial $R$, one can reorder 
$\zeta^0v_l \otimes \zeta^0v_k \otimes p$ 
in two ways along the isomorphism 
$\Aff(B_l) \otimes \Aff(B_k) \otimes 
\Aff(\Pth) \simeq 
\Aff(\Pth) \otimes \Aff(B_k) \otimes \Aff(B_l)$
as follows:

\setlength{\unitlength}{1mm}
\begin{picture}(80,55)(-20,0)

\put(0,10){\line(1,0){8}}
\put(0,20){\line(1,0){8}}

\put(11,9){$\cdots$} 
\put(12,19){$\cdots$}

\put(18,10){\line(1,0){8}}
\put(18,20){\line(1,0){8}}

\put(4,7){\line(0,1){16}}
\put(22,7){\line(0,1){16}}

\put(-4,19){$\scriptstyle v_{k}$}
\put(-4,9){$\scriptstyle v_{l}$}

\put(12,25){$\scriptstyle p$}
\put(11,14){$\scriptstyle T_{k}(p)$}
\put(10,4){$\scriptstyle T_lT_{k}(p)$}

\put(28,19){$\scriptstyle \zeta^{E_{k}(p)}{v}_{k}$}
\put(27,9){$\scriptstyle \zeta^{E_{l}(T_{k}(p))}{v}_{l}$}

\put(44,10){\line(1,1){10}}
\put(44,20){\line(1,-1){10}}

\put(55,9){$\scriptstyle \zeta^{E_{k}(p)+\delta}\overline{v}_{k}$}
\put(55,19){$\scriptstyle \zeta^{E_{l}(T_{k}(p))-\delta}\overline{v}_{l}$}

\put(-10,14){$=$}

\put(0.2,35){\line(1,1){10}}
\put(0.2,45){\line(1,-1){10}}

\put(-4,35){$\scriptstyle v_{l}$}
\put(-4,45){$\scriptstyle v_{k}$}

\put(12,35){$\scriptstyle \zeta^{\delta}\overline{v}_{k}$}
\put(11,45){$\scriptstyle \zeta^{-\delta}\overline{v}_{l}$}

\put(20,35){\line(1,0){8}}
\put(20,45){\line(1,0){8}}

\put(31,34){$\cdots$} 
\put(32,44){$\cdots$}

\put(38,35){\line(1,0){8}}
\put(38,45){\line(1,0){8}}

\put(24,32){\line(0,1){16}}
\put(42,32){\line(0,1){16}}

\put(-4,19){$\scriptstyle v_{k}$}
\put(-4,9){$\scriptstyle v_{l}$}

\put(32,50){$\scriptstyle p$}
\put(31,39){$\scriptstyle T_{l}(p)$}
\put(30,29){$\scriptstyle T_{k}T_{l}(p)$}

\put(47,45){$\scriptstyle \zeta^{E_{l}(p)-\delta}\overline{v}_{l}$}
\put(47,35){$\scriptstyle \zeta^{E_{k}(T_{l}(p))+\delta}\overline{v}_{k}$}

\put(75,9){,}

\end{picture}

\noindent
where the equality of the result is due to the 
Yang-Baxter equation (\ref{eq:ybr}).
Here we have identified the 
outputs with $T_kT_l(p), \zeta^{E_{k}(T_{l}(p))+\delta}v_{k}$, etc.
In particular Proposition \ref{pr:tl} 
guarantees that $\overline{v}_k\otimes T_l(p) \simeq 
T_kT_l(p) \otimes \overline{v}_k$ and 
$\overline{v}_l\otimes p\simeq 
T_l(p) \otimes \overline{v}_l$ up to the spectral parameter.
The sought relations 
$T_lT_k(p) = T_kT_l(p)$ and $E_l(T_k(p))=E_l(p)$
are obtained by comparing the two sides.
\end{proof}

\subsection{Extended affine Weyl group invariance}

Let $s_i\, (i=0,1)$ be the Weyl group operator \cite{K}
acting on any crystal $B$ as
\begin{equation*}
s_i(b) = \begin{cases}
{\tilde f}_i^{\varphi_i(b)-\varepsilon_i(b)}(b) & 
\varphi_i(b) \ge \varepsilon_i(b),\\
{\tilde e}_i^{\varepsilon_i(b)-\varphi_i(b)}(b) & 
\varphi_i(b) \le \varepsilon_i(b)
\end{cases}
\end{equation*}
for $b \in B$.
We extend $\omega$ introduced around (\ref{eq:dynkin}) 
to any $B = B_{l_1}\otimes \cdots \otimes B_{l_k}$ 
by $\omega(B) 
= \omega(B_{l_1})\otimes \cdots \otimes \omega(B_{l_k})$.
Then ${\widetilde W}(A^{(1)}_1)
= \langle \omega, s_0, s_1\rangle$ acts on
$\Pth$ as the extended affine Weyl group of type $A^{(1)}_1$.

The action of ${\tilde f}_i, {\tilde e}_i$ and $s_i$ 
is determined in principle by (\ref{eq:erule}) and (\ref{eq:frule}).
Here we explain the {\em signature rule} to find the action 
on any $B_{l_1} \otimes \cdots \otimes B_{l_k}$ 
which is of great practical use.
It will be the basic ingredient in proving Proposition \ref{pr:w}.
The $i$-signature of an element $b \in B_l$ is the symbol
$\overbrace{-\cdots -}^{\varepsilon_i(b)}
\overbrace{+ \cdots +}^{\varphi_i(b)}$.
The $i$-signature of 
the tensor product $b_1 \otimes \cdots \otimes b_k \in 
B_{l_1} \otimes \cdots \otimes B_{l_k}$ is the array of 
the $i$-signature of each $b_j$.
Here is an example from 
$B_5\otimes B_2\otimes B_1 \otimes B_4$:
\begin{equation*}
\underset{1-{\rm signature}}{{\underset{0-{\rm signature}}
{\phantom{1-signature}}}}
\underset{-++++}{\underset{----+}{11112}} \otimes 
\underset{-+}{\underset{-+}{12}} \otimes 
\underset{-}{\underset{+}{2}} \otimes 
\underset{--++}{\underset{--++}{1122}}
\end{equation*}
where $1122$ for example represents $\fbox{1122} \in B_4$ and 
not $\fbox{1}\otimes \fbox{1} \otimes 
\fbox{2} \otimes \fbox{2} \in B^{\otimes 4}_1$, etc.
In the $i$-signature, one eliminates 
the neighboring pair $+-$ (not $-+$) successively to 
finally reach the pattern 
$\overbrace{-\cdots -}^\alpha\overbrace{+ \cdots +}^\beta$
called reduced $i$-signature.
The result is independent of the order of 
the eliminations when it can be done 
simultaneously in more than one places.
The reduced $i$-signature tells that 
$\varepsilon_i(b_i \otimes \cdots \otimes b_k)=\alpha$ and 
$\varphi_i(b_i \otimes \cdots \otimes b_k)=\beta$. 
In the above example, we get
\begin{equation*}
\underset{1-{\rm signature}}{{\underset{0-{\rm signature}}
{\phantom{1-signature}}}}
\underset{-+\phantom{+++}}{\underset{----\phantom{+}}{11112}} 
\otimes 
\underset{}{\underset{}{12}} \otimes 
\underset{}{\underset{}{2}} \otimes 
\underset{\phantom{--}++}{\underset{\phantom{--}++}{1122}}
\end{equation*}
Thus $\varepsilon_0 = 4,\, \varphi_0 = 2,\; 
\varepsilon_1=1$ and $\varphi_1 = 3$.
Finally 
${\tilde f}_i$ hits the component 
that is responsible for the leftmost $+$ 
in the reduced $i$-signature making it $-$.  
Similarly, ${\tilde e}_i$ hits the component 
corresponding to the rightmost $-$ 
in the reduced $i$-signature making it $+$.
If there is no such $+$ or $-$ to hit, the result 
of the action is $0$. 
$s_i$ acts so as to change the reduced $i$-signature  
$\overbrace{-\cdots -}^\alpha\overbrace{+ \cdots +}^\beta$
into 
$\overbrace{-\cdots -}^\beta\overbrace{+ \cdots +}^\alpha$.
In the above example, we have
\begin{align*}
p &= 
11112 \otimes 12 \otimes 2 \otimes 1122\\
{\tilde f}_0(p) &= 
11112 \otimes 12 \otimes 2 \otimes 1112\\
{\tilde f}_1(p) &= 
11122 \otimes 12 \otimes 2 \otimes 1122\\
{\tilde e}_0(p) &= 
11122 \otimes 12 \otimes 2 \otimes 1122\\
{\tilde e}_1(p) &= 
11111 \otimes 12 \otimes 2 \otimes 1122\\
s_0(p) &= 
11222 \otimes 12 \otimes 2 \otimes 1122\\
s_1(p) &= 
11122 \otimes 12 \otimes 2 \otimes 1222.
\end{align*}
For both $i=0$ and $1$, 
note that $\wt(s_i(p))=-\wt(p)$ for any $p$, 
and $s_i(p)=p$ if $\wt(p)=0$.
In order that ${\tilde e}_1p=0$ to hold for a path $p \in \Pth$,
it is necessary and sufficient that 
the reduced $1$-signature to become $+ \cdots +$, which is 
equivalent to the condition (\ref{eq:lp}).
We note that elimination of the $+-$ pairs is 
described by successive applications of the rule
\begin{equation*}
\overbrace{+\cdots +}^{\alpha}\overbrace{-\cdots -}^{\beta}
\quad \longrightarrow \quad 
\overbrace{+\cdots +}^{(\alpha-\beta)_+}
\overbrace{-\cdots -}^{(\beta-\alpha)_+},
\end{equation*}
where $(c)_+=\max(c,0)$.
The time evolutions $T_l$ and the energy $E_l$ 
enjoy the extended affine Weyl group symmetry.

\begin{proposition}\label{pr:w}
For any $w \in {\widetilde W}(A^{(1)}_1), \; p \in \Pth$ 
and $l \in \Z_{\ge 1}$, the commutativity 
$wT_l(p) = T_l(w(p))$ and 
the invariance $E_l(w(p)) = E_l(p)$ are valid.
\end{proposition}

This property persists in the most general $A^{(1)}_n$ case as 
announced in \cite{KT2}. 

\begin{proof}
Let the reduced $i$-signature of $p$ be
$\overbrace{-\cdots -}^\alpha\overbrace{+\cdots +}^\beta$
and the one for $T_l(p)$ be 
$\overbrace{-\cdots -}^{\alpha'}\overbrace{+\cdots +}^{\beta'}$.
Since $T_l$ is weight preserving, one has 
$\alpha-\beta = \alpha'-\beta'$.
First we show the assertion for $w=s_i$.
We assume $\alpha<\beta$, hence
$s_i(p) = {\tilde f}_i^{\beta-\alpha}(p)$ and 
$s_i(T_l(p)) = {\tilde f}_i^{\beta-\alpha}(T_l(p))$. 
The proof for the case $\alpha \ge \beta$ is parallel.
Let $v_l \in B_l$ be the element specified by (\ref{eq:tle})
from $p$ and set $a=\varepsilon_i(v_l)$ and 
$b=\varphi_i(v_l)$.
The $i$-signature of (\ref{eq:tle}) reads
\begin{equation}\label{eq:sig}
\overbrace{-\cdots -}^a\overbrace{+\cdots +}^b\vert
\overbrace{-\cdots -}^\alpha\overbrace{+\cdots +}^\beta
\quad = \quad 
\overbrace{-\cdots -}^{\alpha'}\overbrace{+\cdots +}^{\beta'}
\vert\overbrace{-\cdots -}^a\overbrace{+\cdots +}^b,
\end{equation}
where $\vert$ signifies the position of $\otimes$
that separates $v_l$ with $p$ or $T_l(p)$.

(i) Case $\alpha \ge b$.
Comparing the reduced $+$ signature, 
we see $(\alpha <) \beta=b+(\beta'-a)_+
\le \alpha+(\beta'-a)_+$, compelling $\beta'>a$.
Thus we find $\beta'-a \ge \beta-\alpha$.
By taking these facts into account, (\ref{eq:sig}) 
is reduced to
\begin{equation*}
\overbrace{-\cdots -}^a\vert
\overbrace{- \cdot\cdot -}^{\alpha-b}
\overbrace{+\cdots +}^\beta \;=\; 
\overbrace{-\cdots -}^{\alpha'} 
\overbrace{+\cdot\cdot +}^{\beta'-a}\vert
\overbrace{+\cdots +}^b.
\end{equation*}
In view of this and $\beta'-a \ge \beta-\alpha$,
application of ${\tilde f}_i^{\beta-\alpha}$
to  (\ref{eq:tle}) yields 
$v_l\otimes {\tilde f}_i^{\beta-\alpha}(p) \simeq 
{\tilde f}_i^{\beta-\alpha}(T_l(p)) \otimes 
\zeta^{E_l(p)}v_l$. 
Therefore $v_l \otimes s_i(p) \simeq s_i(T_l(p)) \otimes 
\zeta^{E_l(p)}v_l$, saying 
$T_l(s_i(p)) = s_i(T_l(p))$ and 
$E_l(s_i(p)) = E_l(p)$.

(ii) Case $\alpha < b$.
Comparing the reduced $+$ signature in (\ref{eq:sig}),
we see $(b<) b+\beta-\alpha = b+(\beta'-a)_+$, compelling 
$\beta'-a>0$. Thus we find 
$\beta'-a=\beta-\alpha$. 
By taking these facts into account, (\ref{eq:sig}) 
is reduced to
\begin{equation*}
\overbrace{-\cdots -}^a
\overbrace{+ \cdot\cdot +}^{b-\alpha}\vert
\overbrace{+\cdots +}^\beta \;= \;
\overbrace{-\cdots -}^{\alpha'} 
\overbrace{+\cdot\cdot +}^{\beta-\alpha}\vert
\overbrace{+\cdots +}^b.
\end{equation*}
In view of this, the application of ${\tilde f}_i^{b+\beta-2\alpha}$ 
to (\ref{eq:tle}) yields
$\zeta^{\delta}v'_l
\otimes {\tilde f}_i^{\beta-\alpha}(p) \simeq 
{\tilde f}_i^{\beta-\alpha}(T_l(p))\otimes 
\zeta^{\delta+E_l(p)}v'_l$ with 
$v'_l = {\tilde f}_i^{b-\alpha}v_l$ and $\delta=(b-\alpha)\delta_{i,0}$.
This is equivalent to 
$v'_l \otimes s_i(p) \simeq s_i(T_l(p)) \otimes 
\zeta^{E_l(p)}v'_l$ saying again that 
$T_l(s_i(p)) = s_i(T_l(p))$ and 
$E_l(s_i(p)) = E_l(p)$.

Next we show the assertion for $w=\omega$.
Due to the symmetry (\ref{eq:dynkin}), 
the relation (\ref{eq:tle}) implies 
$\omega(v_l)\otimes \omega(p) \simeq 
\omega(T_l(p)) \otimes \omega(v_l)$ 
at least under $B_l \otimes \Pth \simeq \Pth \otimes B_l$
forgetting the spectral parameter.
Thus we obtain $T_l(\omega(p)) = \omega(T_l(p))$.
In (\ref{eq:tle}), imagine the process of sending $v_l$ to the right 
through $p$ along Figure \ref{fig:xp} using the local rule 
listed in Figure \ref{fig:cr}.
By the definition, $-E_l(p)$ is the sum of $H$ attached to the 
intermediate vertices.
Upon application of $\omega$, 
the four types of vertices in Figure \ref{fig:cr} are interchanged 
horizontally.
But this does not alter the sum of $H$ since 
the number of bottom two types are equal because 
of $\wt(p) = \wt(T_l(p))$.
\end{proof}

\begin{example}\label{ex:s}
The following is a commutative diagram among the 
paths in $\Pth = B^{\otimes 6}_1$.
\begin{align*}
&
212221 \overset{s_0}{\longmapsto} 
211121 \overset{s_1}{\longmapsto} 
222121 \overset{\omega}{\longmapsto} 111212\\
&\!{\scriptstyle T_2} \;
\downarrow \;\;\quad \qquad\quad \downarrow 
\quad  \;\;\qquad \quad \downarrow  
\quad \;\;\;\qquad \quad \downarrow \\
&
121222 \overset{s_0}{\longmapsto} 
121112 \overset{s_1}{\longmapsto} 
122212 \overset{\omega}{\longmapsto} 211121
\end{align*}
For all the paths here, we have $E_l(p)=2$ for any $l \ge 1$.
\end{example}

The time evolution $T_l$ defined by (\ref{eq:tle})
has a simple description for $l$ sufficiently large.

\begin{proposition}\label{pr:ts} 
For any path $p \in \Pth$, there exists $k \in \Z_{\ge 1}$ 
such that $T_l(p)$ is independent of $l$ for $l \ge k$.
Denoting it by $T_\infty(p)$, one has
\begin{equation}\label{eq:ts} 
T_\infty(p)=\begin{cases}
\omega(s_0(p)) & \hbox{ if } \wt(p)\ge 0,\\
\omega(s_1(p)) & \hbox{ if } \wt(p) \le 0.
\end{cases}
\end{equation} 
In particular, $T_\infty(p)=\omega(p)$ is valid if $\wt(p)=0$. 
\end{proposition}

The combination $\omega s_i$ represents a translation in 
${\widetilde W}(A^{(1)}_1)$.
On account of (\ref{eq:tlj}) and Theorem \ref{th:main}, 
the minimum of such $k$ is given by $k=j_s$ 
with $j_s$ specified in (\ref{eq:mh}).
This is the amplitude of the largest soliton involved in $p$.
The key to the connection of the combinatorial $R$ 
and the signature rule is 
\begin{lemma}\label{lem:r}
Let $(l-a,a) \otimes b \simeq b' \otimes (l-a',a')$ under 
the combinatorial 
$R:  B_l \otimes B_1 \simeq B_1 \otimes B_l$.
If $0 \le a < l$, then $a'$ and $b'$ are expressed as
\begin{equation}\label{eq:ab}
a' = (a-\varepsilon_0(b))_+ + \varphi_0(b),\quad 
b' = \omega({\tilde e}^{(\varepsilon_0(b)-a)_+}_0(b)).
\end{equation}
\end{lemma}
\begin{proof}
In Figure \ref{fig:cr}, the top right pattern does not occur.
Then (\ref{eq:ab}) is directly checked case by case.
\end{proof}

\begin{proof}[Proof of Proposition \ref{pr:ts}]
In view of $\omega(s_1(p)) = s_0(\omega(p))$ for any $p$,
the two formulas in (\ref{eq:ts}) are equivalent.
Henceforth we assume $\wt(p)\ge 0$, i.e., 
$\varepsilon_0(p) \ge \varphi_0(p)$.
Notice that $a'$ in (\ref{eq:ab}) is the number of $+$ 
in the reduced $0$-signature of 
$(l-a,a) \otimes b \in B_l \otimes B_1$.
Consider Figure \ref{fig:xp} with 
$p = b_1 \otimes \cdots \otimes b_L$ and 
$(x_1,x_2)=(l-x_2,x_2) \in B_l$ with any fixed $p$ and $x_2$.
It determines 
$b'_1,\ldots, b'_L$ and $(y_1,y_2)$ 
successively from the left 
by the local rule specified in Lemma \ref{lem:r}
by taking $l$ sufficiently large.
{}From the above fact and the signature rule,
we find that 
the choice $(x_1,x_2)=(l,0)$ leads to 
$(y_1,y_2) = (l-\varphi_0(p),\varphi_0(p))$ which is 
called $v_l$ in (\ref{eq:vt}).
Moreover the latter formula in (\ref{eq:ab}) implies that
the choice $(x_1,x_2) = (l-a,a)$ yield
$b'_1 \otimes \cdots \otimes b'_L = 
\omega({\tilde e}_0^{\,\varepsilon_0(p)-a}(p))$ 
provided $a \le \varepsilon_0(p)$.
By the assumption we are allowed to take $a=\varphi_0(p)$,
leading to
$b'_1 \otimes \cdots \otimes b'_L = \omega(s_0(p))$.
But this must be $T_l(p)$ due to (\ref{eq:tle}).
\end{proof}

\begin{example}
Take the path $p = 121221221111 \in B^{\otimes 12}_1$ 
satisfying $\wt(p) >0$.
To compute $T_\infty(p) = \omega(s_0(p))$, we 
display the $0$-signature:
\begin{equation*}
p=
\underset{-}{1}\;\,\,
\underset{(+}{2}\;
\underset{-)}{1}\;
\underset{(+}{2}\;
\underset{(+}{2}\;
\underset{-)}{1}\;
\underset{(+}{2}\;
\underset{(+}{2}\;
\underset{-)}{1}\;
\underset{-)}{1}\;
\underset{-)}{1}\;\,
\underset{-}{1},
\end{equation*}
where the $(+-)$ pairs
to be eliminated successively to get the reduced $0$-signature
are indicated by parentheses.
Those paired remain unchanged under $s_0$, and we thus find 
\begin{equation*}
s_0(p) = 2 \;\; \, 2 \;\;\,  1 \;\; \, \, 2 \;\;\, 2\;\;  \, 
1\;\;  \,  2\;\;\,   2\;\;\,   1\;\;\,   1\;\;\,   1\;\;\,   2.\quad\;\;
\end{equation*}
Interchanging $1$ and $2$ here we obtain 
$T_\infty(p) = 112112112221$. 
In this example the composition $\omega s_0$ has the effect of 
doing nothing for the unpaired $1$ and interchanging 
the paired $1$ and $2$.
When there remain unpaired $+$ in the reduced $0$-signature
like $p=11122$, the effect of $\omega s_0$ 
is described in the same manner 
if those $+$ are paired with $-$ cyclically.  
Thus the formula (\ref{eq:ts}) reproduces the description of $T_\infty$
in terms of the ``arc rule" in \cite{YYT}.
\end{example}

\section{Inverse scattering method}\label{sec:ism}

Here we use the rigged configurations 
and their bijective correspondence $\phi$ 
with the highest paths $\Pth_+$ \cite{KKR,KR}
summarized in Appendix \ref{app:kkr}.

\subsection{Action variable}\label{subsec:acv}

The set of states  $\Pth=B^{\otimes L}_1$ 
of the periodic box-ball system
is decomposed into the disjoint union 
according to the value of the conserved quantities 
$\{E_l\mid l \in \Z_{\ge 1}\}$ which we called energy.
Our aim here is to determine their spectrum
by making a connection with rigged configurations.
It will be attained in Proposition \ref{pr:em} 
and (\ref{eq:dud}).

\begin{lemma}\label{lem:dp+}
For any path $p \in \Pth$ with $\wt(p)\ge 0$, there is 
an integer $d \in \Z$ and a highest path $p_+ \in \Pth_+$ such that 
$p = T^d_1(p_+)$.
\end{lemma}
The pair $(d, p_+)$ is not unique in general even if $d$ is restricted to
$0 \le d < L$.
The proof is elementary, and is illustrated along 
\begin{example}\label{ex:dp+}
Take $p=2211221112122111221$ with length $L=19$.
Regard the letters 1 and 2 in the doubled path
$p \otimes p$ 
as the arrows
$\nearrow$ and $\searrow$ respectively, and construct 
a length $2L$ trail by following them.
In the example, one has the following:

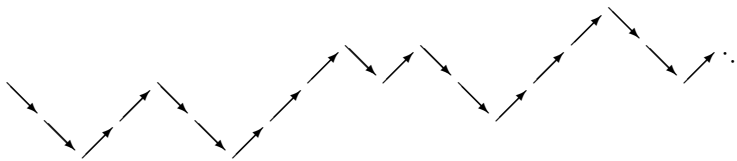
\begin{picture}(200,80)(-30,-40)
\unitlength 0.5mm
\put(0,0){\vector(1,-1){8}}
\put(10,-10){\vector(1,-1){8}}
\put(20,-20){\vector(1,1){8}}
\put(30,-10){\vector(1,1){8}}
\put(40,0){\vector(1,-1){8}}
\put(50,-10){\vector(1,-1){8}}
\put(60,-20){\vector(1,1){8}}
\put(70,-10){\vector(1,1){8}}
\put(80,0){\vector(1,1){8}}
\put(90,10){\vector(1,-1){8}}
\put(100,0){\vector(1,1){8}}
\put(110,10){\vector(1,-1){8}}
\put(120,0){\vector(1,-1){8}}
\put(130,-10){\vector(1,1){8}}
\put(140,0){\vector(1,1){8}}
\put(150,10){\vector(1,1){8}}
\put(160,20){\vector(1,-1){8}}
\put(170,10){\vector(1,-1){8}}
\put(180,0){\vector(1,1){8}}
\put(190,6){$\cdot$}
\put(192,4){$\cdot$}
\put(194,2){$\cdot$}
\end{picture}

\noindent
Find a local minimum 
that are not higher than any 
other points on its right.
In view of $\wt(p)\ge 0$, such minimums can be 
found within the first $L$ steps of the trail. 
To obtain $p_+$, read the trail from 
any one of the minimums to the right for $L$ steps. 
In the example,  
there are three such minimums and accordingly
$p=T^{2}_1(p_1)=T^{6}_1(p_2) = T^{13}_1(p_3)$ 
with $p_1 = 1122111212211122122$,
$p_2 = 1112122111221221122$ and 
$p_3 = 1112212211221112122$.
\end{example}

Now we introduce 
\begin{equation}\label{eq:m}
\M = \{ m=(m_j)_{j\ge 1} \mid m_j \in \Z_{\ge 0},\, 
\sum_{j\ge 1} jm_j \le L/2 \}.
\end{equation}
By the definition $m_j = 0$ for $j \gg 1$, and 
we identify $(m_1,\ldots, m_l, 0, 0, \ldots)$ with 
$(m_1, \ldots, m_l)$.
$\M$ is a finite set.  
We call its elements {\em action variables}.
$m=(m_j)$ is identified with the Young diagram in which
$m_j$ is the number of length $j$ rows.
We call the $m_j \times j$ rectangle consisting of 
the length $j$ rows a {\em block}. 

Based on the KKR bijection we define the map 
\begin{equation}\label{eq:am}
\begin{split}
\mu: \;\Pth &\longrightarrow \M\\
p \;& \longmapsto \; m
\end{split}
\end{equation}
as follows. 
Find $(d, p_+) \in \Z\times \Pth_+$  
such that $p = T_1^d(p_+)$ if $\wt(p)\ge 0$ and 
$\omega(p) = T_1^d(p_+)$ if $\wt(p)<0$
according to Lemma \ref{lem:dp+}.
Then $m$ is obtained as the configuration part of the
rigged configuration $(m,J) = \phi(p_+)$. 
As noted in Lemma \ref{lem:dp+},  the choice of $(d,p_+)$ is not unique. 
Therefore to make sense of the above definition, 
we have to guarantee that  
if $T_1^d(p_+)=T_1^{d'}(p'_+)$ and 
$\phi(p_+) = (m,J), \phi(p'_+) = (m',J')$, then $m=m'$ holds.
To see this, note from the highest condition (\ref{eq:lp}) that 
the situation $T_1^{d-d'}(p_+)=p'_+$ can happen
only when $p_+ = r\otimes q$ and $p'_+ = q \otimes r$ 
for some shorter highest paths $q \in B^{\otimes d-d'}_1$ 
and $r \in B^{\otimes L-d+d'}_1$.
(Without loss of generality, $0 \le d-d'<L$ may be assumed.)
Then the assertion $m=m'$ is included in Lemma \ref{lem:rotate}.

\begin{example}
Consider $p$ in Example \ref{ex:dp+}, which can be expressed 
by the highest paths $p_1, p_2$ and $p_3$.
Computing $\phi(p_i)$, one finds that they all lead to 
$m=(m_1,m_2,m_3)=(2,2,1)$.  The result is depicted as follows: 
\begin{equation*}
\begin{picture}(60,53)(-60,-50)
\put(0,0){\line(0,-1){50}}
\put(10,0){\line(0,-1){50}}
\put(20,0){\line(0,-1){30}}
\put(30,0){\line(0,-1){10}}

\put(-27,-18){$\mu$}
\put(-135,-25){$2211221112122111221 \;\;\longmapsto$}

\put(0,0){\line(1,0){30}}

\put(0,-10){\line(1,0){30}}

\put(0,-20){\line(1,0){20}}

\put(0,-30){\line(1,0){20}}

\put(0,-40){\line(1,0){10}}

\put(0,-50){\line(1,0){10}}
\end{picture}
\end{equation*} 
\end{example}

Let us relate the action variables $m$, namely 
the configurations in the KKR theory, to the energy 
$E_l(p)$ of a path $p$ determined by (\ref{eq:tle}).

\begin{proposition}\label{pr:em}
For any path $p \in \Pth$, its energy is expressed as
\begin{equation}\label{eq:em}
E_l(p) = \sum_{k \ge 1}\min(l,k)m_k
\end{equation}
in terms of the action variable $m=(m_j) = \mu(p) \in \M$. 
\end{proposition}

The proof reduces to the highest case $p \in \Pth_+$, 
and is given in Appendix \ref{app:em}.
According to Remark \ref{rem:sol}, 
$m=(m_j)$ has the meaning of {\em soliton content}, i.e., 
there are $m_j$ solitons with length $j$.
See also Example \ref{ex:omega}.
The right hand side of (\ref{eq:em}) 
is the number of boxes in the first $l$ columns of 
the Young diagram corresponding to $m$.
Proposition \ref{pr:em} determines the range of 
the energy, i.e., the spectrum of the periodic 
box-ball system, in terms of the 
action variable $m \in \M$ (\ref{eq:m}).
In particular it implies the property:
\begin{equation*}
0 \le E_1(p) < E_2(p) < \cdots < E_s(p) = E_{s+1}(p) 
= \cdots = M,
\end{equation*}
where $M=\sum_{k\ge 1}km_k$ and 
$s$ is the greatest integer such that $m_s >0$, 
or equivalently, the length of the longest (top) row of 
the Young diagram for $m$.
The relation (\ref{eq:em}) 
also determines $m$ from $\{E_l\}$.
Combining Proposition \ref{pr:em} with 
Theorem \ref{th:te} and Proposition \ref{pr:w}, we obtain 
\begin{corollary}\label{cor:mu}
The action variable is invariant under the time evolutions
and the extended affine Weyl group. Namely,
$\mu(T_l(p)) = \mu(w(p)) = \mu(p)$ 
for any $l \in \Z_{\ge 1}$ and 
$w \in {\widetilde W}(A^{(1)}_1)$.
\end{corollary}

For each $m = (m_j) \in \M$,
we introduce the corresponding ``level set", namely, 
the set of paths $p$ characterized by $\mu(p)=m$.
Rephrasing the condition in terms of the energy 
by Proposition \ref{pr:em}, we put
\begin{equation}\label{eq:p}
{\widehat {\mathcal P}}(m) = \{ p \in \Pth \mid 
E_l(p) = \sum_{k\ge 1}\min(l,k)m_k \,\hbox{ for any } l \},
\quad 
\Pth_+(m) = {\widehat {\mathcal P}}(m) \cap \Pth_+.
\end{equation}
One has the disjoint union decomposition:
\begin{equation}\label{eq:dud}
\Pth = \sqcup_{m \in \M} {\widehat {\mathcal P}}(m),\quad
\Pth_+ = \sqcup_{m \in \M} \Pth_+(m).
\end{equation}

{}From Corollary \ref{cor:mu}, we see that 
each ${\widehat {\mathcal P}}(m)$ 
is invariant under time evolutions $T_l$ as well as 
the extended affine Weyl group ${\widetilde W}(A^{(1)}_1)$.
(On the other hand highest paths can not remain 
highest under them in general.)
In the KKR algorithm, $M = \sum_{k\ge 1}km_k$ is 
the number of $2 \in B_1$ contained in a highest 
(hence $\wt(p)\ge 0$) path.
Thus we find
\begin{equation*}
\wt({\widehat {\mathcal P}}(m)) = 
\{p_\infty \Lambda_1, -p_\infty \Lambda_1\},\quad
\wt(\Pth_+(m)) = \{p_{\infty}\Lambda_1\},
\end{equation*}
where $p_\infty = L-2M \ge 0$ 
is the limiting (minimum) value of the vacancy number (\ref{eq:va}). 
Thus we have a further decomposition with respect to the weights:
\begin{equation*}
{\widehat {\mathcal P}}(m) = \begin{cases}
\Pth(m) \sqcup \omega(\,\Pth(m)) & \hbox{ if } L>2M,\\
\Pth(m) & \hbox{ if } L=2M,
\end{cases}
\end{equation*}
where $M = \sum_{k\ge 1}km_k$ as above and 
\begin{equation*}
\Pth(m) = \{p \in {\widehat {\mathcal P}}(m) \mid \wt(p) \ge 0\}
= \{p \in {\widehat {\mathcal P}}(m) \mid \wt(p) = p_\infty\Lambda_1\}
\end{equation*}
is a fixed weight subset of ${\widehat {\mathcal P}}(m)$.
One has $\Pth(m) \supset \Pth_+(m)$ and 
\begin{align}
&\{p \in \Pth \mid \wt(p) \ge 0 \} = 
\sqcup_{m \in \M}\Pth(m), \label{eq:pd}\\
&\vert \Pth(m) \vert = 
\frac{\vert {\widehat {\mathcal P}}(m) \vert}
{\vert \wt({\widehat {\mathcal P}}(m)) \vert}\nonumber
\end{align}
for any $m \in \M$.
This cardinality will be evaluated explicitly in (\ref{eq:cpju})
and (\ref{eq:omega2}).
The set $\Pth(m)$ is still invariant under 
any time evolution $T_l$. 
Now Lemma \ref{lem:dp+} is refined into
\begin{lemma}\label{lem:dpm}
For any path $p \in \Pth(m)$, there is 
an integer $d \in \Z$ and a highest path $p_+ \in \Pth_+(m)$ 
such that $p = T^d_1(p_+)$.
Conversely, for any $p_+ \in \Pth_+(m)$ and $d \in \Z$, 
one has $T^d_1(p_+) \in \Pth(m)$.
\end{lemma}

\subsection{Angle variable}\label{subsec:aa}
Here we construct the set of angle variables 
$\J(m)$ for each prescribed value of the 
action variable $m \in \M$.
We introduce the vacancy numbers 
\begin{equation}\label{eq:va}
p_j = L - 2\sum_{k\ge 1}\min(j,k)m_k.
\end{equation}
Given an element $m \in \M$, put
\begin{equation}\label{eq:mh}
H:=\{j \in \Z_{\ge 1} \mid m_j > 0\} =\{ j_1 < j_2 < \cdots < j_s \}.
\end{equation}
The vacancy numbers satisfy
$L=p_0 > p_1 > p_2 > \cdots > p_{j_s} = \cdots =  p_\infty \ge 0$. 
We introduce
\begin{align}
\Jb &= \Jb(m) = \J_{j_1} \times \J_{j_2} 
\times \cdots \times \J_{j_s}, \label{eq:jbar}\\
\J_j &= \{(J_i)_{i \in \Z} \mid 
J_i \in \Z, \; J_i \le J_{i+1}, \; J_{i+m_j} = J_i + p_j\;
\hbox{ for all } i \}.\label{eq:jcomp}
\end{align}
Note that only $m_j$ of $(J_i)$, say 
$J_1, J_2, \ldots, J_{m_j}$, are independent, and 
the Young diagram 
for the partition 
$(J_{i+m_j}-J_{i}, J_{i+m_j-1}-J_{i},\ldots, J_{i+1}-J_{i})$ 
is contained in the $m_j \times p_j$ rectangle for any $i$.

For $k \in \Z_{\ge 1}$ we introduce a map 
\begin{equation}\label{eq:slide}
\begin{split}
&\sigma_k:  \quad \J_j \,\;\,\;\;\;\longrightarrow 
\;\;\;\J_j\\
&\qquad (J_i)_{i \in \Z} \,\;\longmapsto \;\;(J'_i)_{i \in \Z},\\
&J'_i = J_{i+\delta_{j,k}}+2\min(j,k).
\end{split}
\end{equation}
We extend $\sigma_k$ to the map
$\Jb \rightarrow \Jb$ 
via 
$\sigma_k(\Jb) = 
\sigma_k(\J_{j_1}) \times \sigma_k(\J_{j_2})
\times \cdots \times 
\sigma_k(\J_{j_s})$.
Obviously the inverse $\sigma_k^{-1}$ exists and 
$\sigma_k\sigma_l = \sigma_l\sigma_k$ holds 
for any $k, l \in \Z_{\ge 1}$.
Thus the set of elements 
$\{\sigma_{j_1}^{n_1}\sigma_{j_2}^{n_2}\cdots 
\sigma_{j_s}^{n_s} \mid n_1, n_2, \ldots, n_s \in \Z\}$ 
forms an infinite Abelian group 
${\mathcal A}$ isomorphic to $\Z^s$.
We call an element of ${\mathcal A}$ a {\em slide},
having in mind the index shift 
$i \rightarrow i+\delta_{j,k}$ in (\ref{eq:slide}).
The origin of this curious map will be clarified in 
Lemma \ref{lem:tslide} in connection with the Bethe ansatz.
Given two elements $J, K \in \Jb$, 
say that $J$ and $K$ are {\em equivalent} and denote by
$J \simeq K$ if $J = \sigma(K)$ for some 
$\sigma \in {\mathcal A}$.
Now we define our main object $\J = \J(m)$ by
\begin{equation}\label{eq:j}
\J = \Jb/\simeq.
\end{equation}
An element of $\J$ is called an {\em angle variable}.
For $J \in \Jb$, we will mostly use the same symbol $J \in \J$ 
to denote its image in $\J$ and $[J]$ 
when emphasis is favorable. 
They are the data of the form
$J=(J^{(j)}_i)$ with ${i \in \Z}$ and $j \in H$.
For $J \in \J$ (resp. $J \in \Jb$), 
it is easily seen that 
$J+d := (J^{(j)}_i+d)_{i \in \Z}$ with $d \in \Z$ 
remains in $\J$ (resp. $\Jb$).
For $J \in \Jb(m)$ it is a good exercise to show
\begin{equation}\label{eq:sl}
\sigma_{j_1}^{m_{j_1}} \cdots \sigma_{j_s}^{m_{j_s}}(J)
= J+L.
\end{equation}
Therefore $[J] = [J+L]$ represent the same 
angle variable.
$\J$ is a finite set.
Its cardinality, namely the `volume of the iso-level subset in the 
phase space' will be determined in 
(\ref{eq:cpju}) and (\ref{eq:omega}).
We introduce the time evolution of angle variables 
as follows:
\begin{equation}\label{eq:tlj}
\begin{split}
&T_l:  \quad \quad \;\J(m) \,\;\,\;\;\;\longrightarrow 
\;\;\quad \quad \J(m)\\
&\qquad (J^{(j)}_i)_{i \in \Z, j\in H} \,\;\longmapsto \;\;
(J^{(j)}_i + \min(j,l))_{i \in \Z, j \in H}.
\end{split}
\end{equation}
This is a linear flow on the angle variables. 
It enjoys the commutativity $T_lT_k = T_kT_l$.
The $T_l$ here will be identified with the time evolution 
(\ref{eq:tle}) on the paths $\Pth$ in Theorem \ref{th:main}.
Note that $T_1^d(J) = J+d$ hence 
$T_1^L(J) = J \in \J(m)$.
We also define $T_l$ on $\Jb(m)$ by the same formula as 
(\ref{eq:tlj}).

It is convenient to depict 
the element $(J^{(j)}_i)_{i \in \Z, j \in H} \in \Jb(m)$ 
as the Young diagram $m$ whose rows are assigned 
with these numbers as follows:
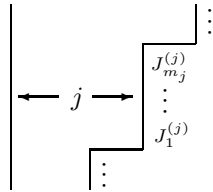
\begin{equation*}
\begin{picture}(50,90)
\put(-30,5){\line(0,1){70}}

\put(-12,40){\vector(-1,0){15}} 
\put(-7,37){$j$}
\put(1,40){\vector(1,0){15}}

\put(40,60){\line(0,1){15}}
\put(44,70){$\cdot$}\put(44,66){$\cdot$}\put(44,62){$\cdot$}
\put(20,60){\line(1,0){20}}
\put(20,20){\line(0,1){40}}
\put(23,50){$\scriptstyle{J^{(j)}_{m_j}}$}
\put(27,40){$\cdot$}\put(27,36){$\cdot$}\put(27,32){$\cdot$}
\put(24,23){$\scriptstyle{J^{(j)}_{1}}$}
\put(0,20){\line(1,0){20}}
\put(4,12){$\cdot$}\put(4,8){$\cdot$}\put(4,4){$\cdot$}
\put(0,5){\line(0,1){15}}
\end{picture}
\end{equation*}
This finite data is enough to recover the whole sequence $(J^{(j)}_i)_{i \in \Z}$ 
on account of the quasi-periodicity 
$J^{(j)}_{i+m_j} = J^{(j)}_i + p_j$ (\ref{eq:jcomp}).
Then the slide $\sigma_{j_k}^l$ (\ref{eq:slide}) translates
the infinite sequence $(J^{(j_k)}_i)_{i \in Z}$ by $l$ 
beside the block-wise shift $2\min(j_k, \alpha)$ on $\J_\alpha$. 
It is helpful to depict the action as
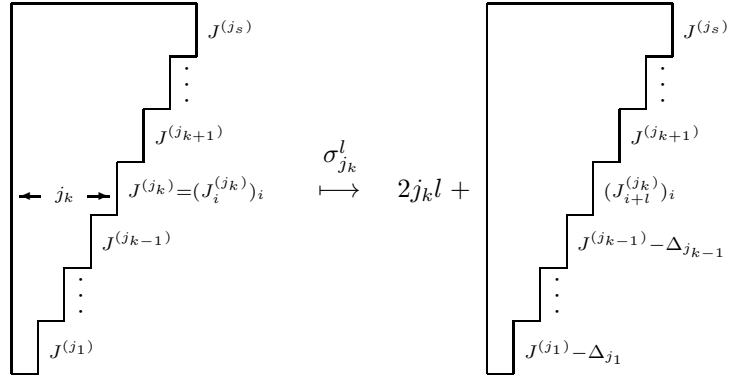
\begin{equation*}
\begin{picture}(100,145)(20,0)
\put(0,140){\line(1,0){70}}\put(0,0){\line(0,1){140}}
\multiput(0,0)(10,20){7}{
\put(0,0){\line(1,0){10}}\put(10,0){\line(0,1){20}}}

\put(74,127){$\scriptstyle J^{(j_{s})}$}

\put(65,113){$\cdot$}
\put(65,107){$\cdot$}
\put(65,101){$\cdot$}

\put(54,87){$\scriptstyle J^{(j_{k+1})}$}
\put(44,67){$\scriptstyle J^{(j_k)}=(J^{(j_k)}_i)_i$}

\put(11,67){\vector(-1,0){8}}
\put(16.5,66){$\scriptstyle j_k$}
\put(29,67){\vector(1,0){8}}
\put(34,47){$\scriptstyle J^{(j_{k-1})}$}

\put(25,33){$\cdot$}
\put(25,27){$\cdot$}
\put(25,21){$\cdot$}

\put(14,6){$\scriptstyle J^{(j_1)}$}

\put(118,80){$\sigma^l_{j_k}$}
\put(116,67){$\longmapsto$}
\end{picture}
\qquad
\begin{picture}(100,145)(-40,0)
\put(0,140){\line(1,0){70}}\put(0,0){\line(0,1){140}}
\multiput(0,0)(10,20){7}{
\put(0,0){\line(1,0){10}}\put(10,0){\line(0,1){20}}}

\put(74,127){$\scriptstyle J^{(j_{s})}$}

\put(65,113){$\cdot$}
\put(65,107){$\cdot$}
\put(65,101){$\cdot$}

\put(54,87){$\scriptstyle J^{(j_{k+1})}$}
\put(44,67){$\scriptstyle (J^{(j_k)}_{i+l})_i$}
\put(34,47){$\scriptstyle J^{(j_{k-1})} - \Delta_{j_{k-1}}$}

\put(25,33){$\cdot$}
\put(25,27){$\cdot$}
\put(25,21){$\cdot$}

\put(14,6){$\scriptstyle J^{(j_1)} - \Delta_{j_1}$}

\put(-34,67){$2j_kl \;+ $}
\end{picture}
\end{equation*}
Here $\Delta_{j} = 2(j_k-j)l$.
Observe that $\sigma_{j_k}^{l}$ only causes the uniform shift 
$2j_kl$
on the upper blocks $\J_{j_{k+1}}, \ldots, \J_{j_s}$,
whereas it relatively induces the back flow $\Delta_{\alpha}$
on the lower blocks $\J_\alpha=\J_{j_1}, \ldots, \J_{j_{k-1}}$. 
A rigged configuration $(m, J)$ with 
$J = (J^{(j)}_i)_{1 \le i \le m_j, j \in H} \in \Rig(m)$ (\ref{eq:rc}) 
also has the data structure that can be depicted as in the above 
picture.

\subsection{Direct and inverse scattering transforms}\label{subsec:dis}

The direct and inverse scattering transforms are the maps
between the path $\Pth$ and the action-angle variables.
In Section \ref{subsec:acv} we have constructed 
the map $\mu$ (\ref{eq:am})   
to find the action variable $m \in \M$ for a given path 
$p \in \Pth$. 
It can be found either from the energy $E_l(p)$ by (\ref{eq:em})
or from the rigged configuration of a highest path described 
under (\ref{eq:am}).
It remains invariant under any time evolution as 
noted in Corollary \ref{cor:mu}.

In the remainder of Section \ref{sec:ism}, we shall only consider 
angle variables supposing that the action 
variable $m$ has been determined.
We assume that $\wt(p) \ge 0$ and formulate the 
inverse scattering method for each $\Pth(m)$ 
appearing in the decomposition (\ref{eq:pd}).
The other case $\wt(p)<0$ is reduced to it by  
\begin{equation*}
T_l(p) = \omega T_l(\omega(p)),
\end{equation*}
owing to  Proposition \ref{pr:w}.

First we 
formulate the direct scattering transform $\Phi$ as follows:
\begin{equation}\label{eq:pj}
\begin{split}
\Phi: \quad &\Pth(m) \longrightarrow \;\Z \times \Pth_+(m)\;\;
\longrightarrow \;\;\quad \Jb(m) \;\;\;\;\longrightarrow 
\quad \;\J(m)\\
&\;\;\;p \;\quad \longmapsto \;\;\;\;(d, p_+) \;\;\quad\longmapsto
\quad \iota(J)+d \;\;\longmapsto \;\;[\iota(J)+d]
\end{split}
\end{equation}
Here the pair $(d,p_+)$ is the one satisfying 
$p = T^d_1(p_+)$ whose existence 
is assured in Lemma \ref{lem:dpm}.
Then the rigging $J \in \Rig(m)$ is specified by 
the KKR bijection $\phi(p_+) = (m,J)$.  The appearance of $m$ here is 
guaranteed by Lemma \ref{lem:dpm}. 
The map $\iota$ is defined by
\begin{equation}\label{eq:iota}
\begin{split}
\iota : \qquad \qquad\qquad\quad\quad
\Rig(m) \qquad \quad\;\;
\qquad&\longrightarrow \qquad\quad
\Jb(m)\\
\left((J^{(j_1)}_i)_{1 \le i \le m_{j_1}}, \ldots, 
(J^{(j_s)}_i)_{1 \le i \le m_{j_s}}\right)
&\mapsto 
\left(({J}^{(j_1)}_i)_{i \in \Z}, \ldots, 
({J}^{(j_s)}_i)_{i \in \Z}\right),
\end{split}
\end{equation}
where the infinite sequence $({J}^{(j)}_i)_{i \in \Z}$ 
is the one that extends 
$(J^{(j)}_i)_{1 \le i \le m_j}$ quasi-periodically as 
${J}^{(j)}_{i+m_j}={J}^{(j)}_i+p_j$ for all $i \in \Z$.
The resulting sequence $({J}^{(j)}_i)_{i \in \Z}$ 
automatically satisfies (\ref{eq:jcomp}).
The map $\iota$ is just the embedding of the 
quasi-periodic extension, hence an injection.

Note that the $(d,p_+)$ satisfying 
$p = T^d_1(p_+)$ is not unique for a given $p$.
Therefore to assure the well-definedness of $\Phi$,
one has to show the $\Rightarrow$ part of 
\begin{proposition}\label{pr:migi}
Let $p_+, p'_+ \in \Pth_+(m)$ be the highest paths and 
$J, J' \in \Rig(m)$ be the corresponding rigging, 
namely, $\phi(p_+) = (m,J)$ and $\phi(p'_+) = (m,J')$.
Then the following relation is valid:
\begin{equation}\label{eq:key}
T^d_1(p_+) = T^{d'}_1(p'_+) \Longleftrightarrow 
\iota(J) + d \simeq \iota(J')+d' \in \Jb(m),
\end{equation}
\end{proposition}

The proof is available in Appendix \ref{app:proof}.

\begin{example}\label{ex:neqn}
Consider the length $L=25$ path
$p=2122112211221111222111122$.
According to Lemma \ref{lem:dp+}, it can be 
expressed as
$p = T_1^{12}(p_+) = T^{19}_1(p'_+)$ in terms of the 
highest paths:
\[
p_+ = 1111222111122212211221122
\quad
p'_+ =1111222122112211221111222.
\]
One applies the KKR bijection to get $\phi(p_+)=(m,J)$ and 
$\phi(p'_+) = (m,J')$. 
The elements $\iota(J)+12$ and $\iota(J')+19$ 
of $\Jb(m)$ are depicted as
\begin{equation*}
\begin{picture}(60,53)(-10,-50)
\put(0,0){\line(0,-1){50}}
\put(10,0){\line(0,-1){50}}
\put(20,0){\line(0,-1){40}}
\put(30,0){\line(0,-1){20}}
\put(40,0){\line(0,-1){10}}

\put(-33,-25){$12\;+$}

\put(0,0){\line(1,0){40}}
\put(43,-9){1}
\put(0,-10){\line(1,0){40}}\put(-7,-8){1}
\put(33,-19){1}
\put(0,-20){\line(1,0){30}}\put(-7,-18){3}
\put(23,-29){7}
\put(0,-30){\line(1,0){20}}\put(-7,-33){7}
\put(23,-39){7}
\put(0,-40){\line(1,0){20}}\put(-12,-48){15}
\put(13,-49){10}
\put(0,-50){\line(1,0){10}}
\end{picture}
\qquad\qquad\qquad
\begin{picture}(60,53)(-10,-50)
\put(0,0){\line(0,-1){50}}
\put(10,0){\line(0,-1){50}}
\put(20,0){\line(0,-1){40}}
\put(30,0){\line(0,-1){20}}
\put(40,0){\line(0,-1){10}}

\put(-33,-25){$19\;+$}

\put(0,0){\line(1,0){40}}
\put(43,-9){0}
\put(0,-10){\line(1,0){40}}\put(-7,-8){1}
\put(33,-19){3}
\put(0,-20){\line(1,0){30}}\put(-7,-18){3}
\put(23,-29){4}
\put(0,-30){\line(1,0){20}}\put(-7,-33){7}
\put(23,-39){4}
\put(0,-40){\line(1,0){20}}\put(-12,-48){15}
\put(13,-49){5}
\put(0,-50){\line(1,0){10}}
\end{picture}
\end{equation*} 
where $1,3,7,15$ on the left are the vacancy numbers 
(\ref{eq:va}) exhibited for convenience.
The equivalence $\iota(J)+12 \simeq \iota(J')+19$ 
is realized by $\sigma_3\big(\iota(J)+12\big)=\iota(J')+19$.
\end{example}

Next we show that the direct scattering map 
$\Phi$ (\ref{eq:pj}) is invertible, which yield the 
inverse scattering map $\Phi^{-1}$.
For this two properties are to be established.
First, any element in $\Jb(m)$ 
must be equivalent to the form $\iota(J) + d$
for some $J \in \Rig(m)$ and $d \in \Z$, which 
ensures the existence of 
an inverse image for the middle arrow in (\ref{eq:pj}).
($\phi(p_+) = (m,J)$.)
Second, any two equivalent forms 
$\iota(J)+d \simeq \iota(J')+d' \in \Jb(m)$ 
must be pulled back to the same path 
in $\Pth(m)$ in (\ref{eq:pj}).
The second property is nothing but the $\Leftarrow$ part of 
Proposition \ref{pr:migi}.
Thus it remains to verify the first property.
This is done in
\begin{lemma}\label{lem:jrc}
For any $\overline{J} \in \Jb(m)$,
there exist $d \in \Z$ and $J \in \Rig(m)$ 
such that $\overline{J} \simeq \iota(J) + d$.
\end{lemma}
\begin{proof}
Let $H = \{j_1 < \cdots < j_s\}$ 
be the list of lengths of rows in $m$ as in (\ref{eq:mh}).
We give a concrete algorithm to repair $\overline{J}$ by a slide so that  
$\sigma_{j_1}^{n_1}\cdots \sigma_{j_s}^{n_s} \overline{J}
=\iota(J)+d$ holds for some $J \in \Rig(m)$ and $d \in \Z$.
Such a slide is not unique in general but this relation 
can always be achieved gradually from the longer rows in $m$
as explained below. 
First concentrate only on the two blocks 
$\J_{j_{s-1}}$ and $\J_{j_s}$.
In the slide $\sigma_{j_{s-1}}$ 
(\ref{eq:slide}), the quantity $2\min(j,j_{s-1})$ 
is a common constant $2j_{s-1}$ on them.
Therefore by applying $\sigma^{n}_{j_{s-1}}$ for some $n$,
one can make the resulting 
$(J^{(j_{s-1})}_i)_i \in \J_{j_{s-1}}$ and 
$(J^{(j_{s})}_i)_i \in \J_{j_s}$ to 
satisfy $d \le J^{(j_{s-1})}_1 \le \cdots \le J^{(j_{s-1})}_{m_{j_{s-1}}}
\le d+p_{j_{s-1}}$ and 
$d \le J^{(j_{s})}_1 \le \cdots \le J^{(j_{s})}_{m_{j_s}} \le d+p_{j_{s}}$ 
for some $d$. The effect of 
the common change $2\min(j,j_{s-1})$ can be absorbed into $d$.
Next one uses $\sigma_{j_{s-2}}$ similarly to adjust 
$\J_{j_{s-2}}$ to the amended $\J_{j_{s-1}}, \J_{j_s}$ 
up to a redefinition of $d$ without violating the 
foregoing adjustment among the latter two.
This process, although the adjustment 
is not unique in general, works through until 
$d \le J^{(j)}_1 \le \cdots \le J^{(j)}_{m_j} \le d+p_j$ is achieved 
for all $j \in H$ for some $d$. 
\end{proof}

\begin{example}
\begin{equation*}
\begin{picture}(60,53)(10,-50)
\put(0,0){\line(0,-1){50}}
\put(10,0){\line(0,-1){50}}
\put(20,0){\line(0,-1){30}}
\put(30,0){\line(0,-1){10}}

\put(0,0){\line(1,0){30}}\put(33,-9){-1}
\put(0,-10){\line(1,0){30}}\put(-7,-9){1}\put(23,-19){5}
\put(0,-20){\line(1,0){20}}\put(-7,-24){3}\put(23,-29){5}
\put(0,-30){\line(1,0){20}}\put(13,-39){-10}
\put(0,-40){\line(1,0){10}}\put(-7,-44){9}\put(13,-49){-14}
\put(0,-50){\line(1,0){10}}
\end{picture}
\quad
\begin{picture}(60,53)(-3,-50)
\put(-38,-25){$= -1\;+$}

\put(0,0){\line(0,-1){50}}
\put(10,0){\line(0,-1){50}}
\put(20,0){\line(0,-1){30}}
\put(30,0){\line(0,-1){10}}

\put(0,0){\line(1,0){30}}\put(33,-9){0}
\put(0,-10){\line(1,0){30}}\put(23,-19){6}
\put(0,-20){\line(1,0){20}}\put(23,-29){6}
\put(0,-30){\line(1,0){20}}\put(13,-39){-9}
\put(0,-40){\line(1,0){10}}\put(13,-49){-13}
\put(0,-50){\line(1,0){10}}
\end{picture}
\quad
\begin{picture}(60,53)(-16,-50)
\put(-27,-25){$-9\;+$}
\put(0,0){\line(0,-1){50}}
\put(10,0){\line(0,-1){50}}
\put(20,0){\line(0,-1){30}}
\put(30,0){\line(0,-1){10}}

\put(-37,-15){$\scriptstyle{\sigma^{-2}_2}$}
\put(-38,-25){$\simeq$}
\put(0,0){\line(1,0){30}}\put(33,-9){0}
\put(0,-10){\line(1,0){30}}\put(23,-19){3}
\put(0,-20){\line(1,0){20}}\put(23,-29){3}
\put(0,-30){\line(1,0){20}}\put(13,-39){-5}
\put(0,-40){\line(1,0){10}}\put(13,-49){-9}
\put(0,-50){\line(1,0){10}}
\end{picture}
\quad
\begin{picture}(60,53)(-29,-50)
\put(-27,-25){$-3\;+$}
\put(0,0){\line(0,-1){50}}
\put(10,0){\line(0,-1){50}}
\put(20,0){\line(0,-1){30}}
\put(30,0){\line(0,-1){10}}

\put(-37,-15){$\scriptstyle{\sigma^{3}_1}$}
\put(-38,-25){$\simeq$}
\put(0,0){\line(1,0){30}}\put(33,-9){0}
\put(0,-10){\line(1,0){30}}\put(23,-19){3}
\put(0,-20){\line(1,0){20}}\put(23,-29){3}
\put(0,-30){\line(1,0){20}}\put(13,-39){9}
\put(0,-40){\line(1,0){10}}\put(13,-49){4}
\put(0,-50){\line(1,0){10}}
\end{picture}
\end{equation*} 
Here the vacancy numbers $1,3,9$ have been 
shown only in the leftmost diagram. 
\end{example}

{}From 
Lemma \ref{lem:jrc} and Proposition \ref{pr:migi},
it follows that 
the $\Phi$ in (\ref{eq:pj}) is a well-defined and invertible map.
Taking the disjoint union over the action variable $m \in \M$,
we obtain

\begin{theorem}\label{th:phi}
The map $\Phi$ in (\ref{eq:pj}) gives the 
bijection among the set of paths 
$\sqcup_{m\in \M}\Pth(m)$ and the 
set of action-angle variables
$\sqcup_{m \in \M}\{(m,J) \mid J \in \J(m)\}$.
\end{theorem}

\subsection{Solution of initial value problem}\label{subsec:ivp}

Now we present our main theorem in this paper.

\begin{theorem}\label{th:main}
The following commutative diagram is valid:
\begin{equation}\label{eq:cd}
\begin{CD}
{\mathcal P}(m) @>{\Phi}>> \J(m) \\
@V{T_l}VV @VV{T_l}V\\
{\mathcal P}(m) @>{\Phi}>> \J(m) 
\end{CD}
\end{equation}
Here $T_l$ on the left and the right are given by
(\ref{eq:tle}) and (\ref{eq:tlj}), respectively. 
\end{theorem}

The proof is included in Appendix \ref{app:main}. 
Since the map $\Phi$ is invertible, Theorem \ref{th:main}
completes the solution of the initial value problem of
the periodic box-ball system by the inverse scattering method.
The time evolution on the paths $\Pth(m)$ has been linearized 
in terms of the angle variables $\J(m)$.
The number of computational steps required for executing  
$\Phi^{-1} \circ T_l^t \circ \Phi$ is independent of $t$. 

It also contains the solution of the initial value problem in 
the box-ball system on the semi-infinite 
lattice $B_1\otimes B_1 \otimes \cdots$ 
as the case $L \rightarrow \infty$.
This limit is well-defined and drastically simplifies our 
construction so far.
For any path $b_1\otimes b_2 \otimes \cdots$ obeying 
the boundary condition $b_i=1$ for $i \gg 1$,
one can make $1^{\otimes n}\otimes b_1 \otimes b_2 \otimes \cdots$
highest for some $n$.
Thus the degree $d \in \Z$ in (\ref{eq:pj}) is always frozen to
$d=0$.
In view of $p_j \rightarrow \infty$ as $L \rightarrow \infty$,
we do not make a non-trivial identification (\ref{eq:j}) under 
any slide (\ref{eq:slide}).
Consequently (\ref{eq:pj}) just becomes 
$\Phi: p \mapsto J$, which is nothing but $\phi$.
Namely, the direct and the inverse scattering transforms are 
the KKR bijection itself (without an upper bound 
on the rigging), reproducing the results in \cite{KOTY,Ta}
essentially.

\begin{example}\label{ex:ivp}
Let us take the length $L=19$ path $p$ in Example \ref{ex:dp+} 
and derive the time evolution 
\begin{equation}\label{eq:kotae}
T^{1000}_2(p) = 1211221112122211221,\qquad
T^{1000}_3(p) = 2112221211221112112
\end{equation}
based on the inverse scattering formalism in Theorem \ref{th:main}.
As noted in Example \ref{ex:dp+}, we have the expression $p=T^2_1(p_1)$ 
in terms of the highest path $p_1$.
By computing the image of the KKR bijection $\phi$ of $p_1$, 
one finds 
\begin{equation*}
\begin{picture}(60,53)(-10,-50)
\put(0,0){\line(0,-1){50}}
\put(10,0){\line(0,-1){50}}
\put(20,0){\line(0,-1){30}}
\put(30,0){\line(0,-1){10}}

\put(-48,-18){$\Phi$}
\put(-65,-25){$p \;\;\longmapsto$}
\put(-29,-25){$2\;+$}

\put(0,0){\line(1,0){30}}
\put(33,-9){1}
\put(0,-10){\line(1,0){30}}\put(-7,-9){1}
\put(23,-19){1}
\put(0,-20){\line(1,0){20}}\put(-7,-24){3}
\put(23,-29){0}
\put(0,-30){\line(1,0){20}}
\put(13,-39){8}
\put(0,-40){\line(1,0){10}}\put(-7,-44){9}
\put(13,-49){4}
\put(0,-50){\line(1,0){10}}
\end{picture}
\end{equation*} 
The vacancy numbers displayed on the left of the 
diagram for convenience will be omitted in the sequel.
By using the linearized time evolution (\ref{eq:tlj}) one has
\begin{equation*}
\begin{picture}(55,53)(0,-50)
\put(0,0){\line(0,-1){50}}
\put(10,0){\line(0,-1){50}}
\put(20,0){\line(0,-1){30}}
\put(30,0){\line(0,-1){10}}

\put(-80,-25){$T^{1000}_2\Phi(p) =\; 2\; + $}

\put(0,0){\line(1,0){30}}
\put(33,-9){2001}
\put(0,-10){\line(1,0){30}}
\put(23,-19){2001}
\put(0,-20){\line(1,0){20}}
\put(23,-29){2000}
\put(0,-30){\line(1,0){20}}
\put(13,-39){1008}
\put(0,-40){\line(1,0){10}}
\put(13,-49){1004}
\put(0,-50){\line(1,0){10}}
\end{picture}
\quad
\begin{picture}(55,53)(-57,-50)
\put(0,0){\line(0,-1){50}}
\put(10,0){\line(0,-1){50}}
\put(20,0){\line(0,-1){30}}
\put(30,0){\line(0,-1){10}}

\put(-59,-13){$\sigma_1^{222}$}
\put(-60,-25){$\simeq \;\; \;\;2446 \; + $}

\put(0,0){\line(1,0){30}}
\put(33,-9){1}
\put(0,-10){\line(1,0){30}}
\put(23,-19){1}
\put(0,-20){\line(1,0){20}}
\put(23,-29){0}
\put(0,-30){\line(1,0){20}}
\put(13,-39){7}
\put(0,-40){\line(1,0){10}}
\put(13,-49){3}
\put(0,-50){\line(1,0){10}}
\end{picture}
\end{equation*} 
The last rigged configuration corresponds to the 
highest path $p'= 1122112112211121222$. 
Therefore the map $\Phi^{-1}$ sends the 
above scattering data to
$T^{2446}_1(p')=T^{14}_1(p')$ yielding the first result in 
(\ref{eq:kotae}). 
Similarly the calculation of $T^{1000}_3(p)$ goes as follows:
\begin{equation*}
\begin{picture}(55,53)(0,-50)
\put(0,0){\line(0,-1){50}}
\put(10,0){\line(0,-1){50}}
\put(20,0){\line(0,-1){30}}
\put(30,0){\line(0,-1){10}}

\put(-80,-25){$T^{1000}_3\Phi(p) =\; 2\; + $}

\put(0,0){\line(1,0){30}}
\put(33,-9){3001}
\put(0,-10){\line(1,0){30}}
\put(23,-19){2001}
\put(0,-20){\line(1,0){20}}
\put(23,-29){2000}
\put(0,-30){\line(1,0){20}}
\put(13,-39){1008}
\put(0,-40){\line(1,0){10}}
\put(13,-49){1004}
\put(0,-50){\line(1,0){10}}
\end{picture}
\quad
\begin{picture}(55,53)(-77,-50)
\put(0,0){\line(0,-1){50}}
\put(10,0){\line(0,-1){50}}
\put(20,0){\line(0,-1){30}}
\put(30,0){\line(0,-1){10}}

\put(-78,-13){$\sigma_1^{740}\sigma_2^{667}$}
\put(-65,-25){$\simeq \;\; \;\;\;7150 \; + $}

\put(0,0){\line(1,0){30}}
\put(33,-9){1}
\put(0,-10){\line(1,0){30}}
\put(23,-19){2}
\put(0,-20){\line(1,0){20}}
\put(23,-29){0}
\put(0,-30){\line(1,0){20}}
\put(13,-39){4}
\put(0,-40){\line(1,0){10}}
\put(13,-49){0}
\put(0,-50){\line(1,0){10}}
\end{picture}
\end{equation*} 
The last rigged configuration corresponds to the 
highest path $p''= 1211221112112211222$. 
Therefore the map $\Phi^{-1}$ sends the 
above scattering data to
$T^{7150}_1(p'')=T^{6}_1(p'')$ yielding the second result in 
(\ref{eq:kotae}). 
\end{example}

\begin{example}\label{ex:tg}
We take the length $L=26$ path $p=12112211122211121112211111$ 
and compute $T_3^{130}(p)$.
$p$ is already highest.
\begin{equation*}
\begin{picture}(60,53)(-10,-50)
\put(0,0){\line(0,-1){50}}
\put(10,0){\line(0,-1){50}}
\put(20,0){\line(0,-1){30}}
\put(30,0){\line(0,-1){10}}

\put(-48,-18){$\Phi$}
\put(-65,-25){$p \;\;\longmapsto$}

\put(0,0){\line(1,0){30}}
\put(33,-9){0}
\put(0,-10){\line(1,0){30}}\put(-7,-9){8}
\put(23,-19){5}
\put(0,-20){\line(1,0){20}}\put(-12,-24){10}
\put(23,-29){0}
\put(0,-30){\line(1,0){20}}
\put(13,-39){8}
\put(0,-40){\line(1,0){10}}\put(-12,-44){16}
\put(13,-49){0}
\put(0,-50){\line(1,0){10}}
\end{picture}
\end{equation*} 

\begin{equation*}
\begin{picture}(55,53)(0,-55)
\put(0,0){\line(0,-1){50}}
\put(10,0){\line(0,-1){50}}
\put(20,0){\line(0,-1){30}}
\put(30,0){\line(0,-1){10}}

\put(-70,-25){$T^{130}_3\Phi(p) =  $}

\put(0,0){\line(1,0){30}}
\put(33,-9){390}
\put(0,-10){\line(1,0){30}}
\put(23,-19){265}
\put(0,-20){\line(1,0){20}}
\put(23,-29){260}
\put(0,-30){\line(1,0){20}}
\put(13,-39){138}
\put(0,-40){\line(1,0){10}}
\put(13,-49){130}
\put(0,-50){\line(1,0){10}}
\end{picture}
\quad
\begin{picture}(55,53)(-77,-55)
\put(0,0){\line(0,-1){50}}
\put(10,0){\line(0,-1){50}}
\put(20,0){\line(0,-1){30}}
\put(30,0){\line(0,-1){10}}

\put(-78,-13){$\sigma_1^{39}\sigma_2^{26}$}
\put(-69,-25){$\simeq \quad \quad 572 \;\;+ $}

\put(0,0){\line(1,0){30}}
\put(33,-9){0}
\put(0,-10){\line(1,0){30}}
\put(23,-19){5}
\put(0,-20){\line(1,0){20}}
\put(23,-29){0}
\put(0,-30){\line(1,0){20}}
\put(13,-39){8}
\put(0,-40){\line(1,0){10}}
\put(13,-49){0}
\put(0,-50){\line(1,0){10}}
\end{picture}
\end{equation*} 
Since $572=22\times L$, we find $T^{130}_3(p) = p$.
\end{example}

\section{Relation with Bethe ansatz at $q=0$}\label{sec:ba}

Here we show that our inverse scattering 
formalism originates in the Bethe ansatz at $q=0$ \cite{KN}. 
The angle variables stem from the logarithmic branch 
of the string center equation and the 
time evolution is the straight motion of its solution.

\subsection{\mathversion{bold} Bethe ansatz at $q=0$}\label{subsec:ba}

Let us quickly recall the relevant results from the  
Bethe ansatz at $q=0$.
For the precise definitions and statements,
we refer to \cite{KN}.  
The Bethe equation for the spin $1/2$ 
one dimensional XXZ chain on length $L$ periodic lattice reads
\cite{Ba}
\begin{equation}\label{eq:be}
\left(\frac{\sin{\pi\!\left(u_i + \sqrt{-1}\hbar\right)}}
{\sin{\pi\!\left(u_i - \sqrt{-1}\hbar\right)}}\right)^{L}
= - \prod_{j=1}^{M}
\frac{\sin\pi\!\left(u_i - u_j + 2\sqrt{-1}\hbar \right)}
{\sin\pi\!\left(u_i - u_j - 2\sqrt{-1}\hbar \right)}
\end{equation}
for $1 \le i \le M$, 
where $0 \le M \le L/2$ is the number of down spins 
preserved by the Hamiltonian.
The system is associated with the quantum affine 
algebra $U_q(A^{(1)}_1)$ with $q=e^{-2\pi\hbar}$.
Fix $m=(m_j) \in \M$.
String solutions are the ones  
in which $\{u_1, \ldots, u_M\}$ are arranged as
\begin{equation}\label{eq:string}
\bigcup_{j \ge 1}
\bigcup_{1\le \alpha\le m_j}
\bigcup_{u^{(j)}_{\alpha} \in {\mathbb R}}
\{ u^{(j)}_{\alpha} + \sqrt{-1}(j+1-2k)\hbar + \epsilon^{(j)}_{\alpha k}
\mid 1 \le k \le j \},
\end{equation}
where $M = \sum_{j}jm_j$ and 
$\epsilon^{(j)}_{\alpha k}$ stands for a small deviation.
$u^{(j)}_{\alpha}$ is the string center of the $\alpha$ th 
string of length $j$.
In this context, the data $m \in \M$ is referred as the 
{\em string content}. 
Let $H$ be as in (\ref{eq:mh}).
For the generic string solution,  
the Bethe equation is linearized at $q=0$ into a logarithmic form 
called the {\em string center equation}:
\begin{equation}\label{eq:sce}
\sum_{k \in H}\sum_{\beta=1}^{m_k} 
A_{j\alpha,k\beta} u^{(k)}_{\beta} \equiv
\frac{1}{2}(p_j + m_j + 1) \quad \mathrm{mod}\ \Z
\end{equation}
for $j \in H$ and $1 \le \alpha \le m_j$.
Here $p_j$ is the vacancy number (\ref{eq:va}) and 
\begin{equation}\label{eq:a}
A_{j\alpha, k\beta}=\delta_{j,k}\delta_{\alpha,\beta}
(p_j+m_j) + 2\min(j,k)-\delta_{j,k}.
\end{equation}
{}From (\ref{eq:deta}), the matrix $A=(A_{j\alpha, k\beta})$ 
is invertible under the condition $m \in \M$.
 
There are a number of conditions which 
the solutions of the string center equation (\ref{eq:sce}) 
are to satisfy or to be identified thereunder.
First, the Bethe vector depends on $u_i$ only via 
$e^{2\pi\sqrt{-1}u_i}$.
Therefore the string center 
should be understood as $u^{(j)}_{\alpha} \in \R/\Z$ rather than $\R$.
Second, the original Bethe equation (\ref{eq:be}) is symmetric
with respect to $u_1, \ldots, u_M$, but 
their permutation does not lead to a new Bethe vector. 
Consequently, we should regard 
\begin{equation*}
(u^{(j)}_1, u^{(j)}_2, \ldots, u^{(j)}_{m_j})
\in \left(\R/\Z\right)^{m_j}/{\frak S}_{m_j}
\end{equation*}
for each $j$. 
Last, we prohibit 
$u^{(j)}_{\alpha}=u^{(j)}_{\beta}$ for $1 \le \alpha \neq \beta \le m_j$
for any $j$.
This is a remnant of the well-known 
constraint on the Bethe roots so that 
the associated Bethe vector does not vanish.
To summarize, we consider {\em off-diagonal solutions} 
$(u^{(j)}_{\alpha})$ to the string center equation (\ref{eq:sce}) 
that live in
\begin{equation}\label{eq:odc}
(u^{(j)}_1, u^{(j)}_2, \ldots, u^{(j)}_{m_j}) \in
\bigl(\left(\R/\Z\right)^{m_j} - \delta_{m_j} \bigr) /{\frak S}_{m_j}
\quad\hbox{ for each } j,
\end{equation}
where $\delta_n  = \{(v_1, \ldots, v_n)
\in (\R/\Z)^n \mid v_\alpha = v_\beta\ \text{for some}\
1 \le \alpha \neq
\beta \le n\}$.
For simplicity we will often say 
{\em Bethe roots} to mean 
the off-diagonal solutions to the string center equation.
Let $\U(m)$ be the set of the Bethe roots
having the string content $m$.

For $m \in \M$, we introduce 
\begin{align}
\Omega(m) &=
(\det F)\prod_{j \in H} \frac{1}{m_j}
\binom{p_j + m_j - 1}{m_j - 1}.\label{eq:omega}\\
F &= (F_{j,k})_{j,k \in H},\quad 
F_{j,k} = \delta_{j,k}p_j + 2\min(j,k)m_k.\label{eq:f}
\end{align}
In case $H=\emptyset$ (i.e., $m=(0,0,\ldots)$), 
we put $\Omega(m) = 1$.
By expanding the determinant, it is easy to see $\Omega(m) \in \Z$.
Moreover, (\ref{eq:detf}) below tells that 
$\Omega(m) \in \Z_{\ge 1}$ and 
(\ref{eq:omega}) is also expressed as
\begin{equation}\label{eq:omega2}
\Omega(m) = \frac{L}{p_{j_s}}
\prod_{j \in H}\binom{p_j+m_j-1}{m_j},
\end{equation}
where $p_{j_s} = L-2M$ with $M = \sum_jjm_j$.
In case $p_{j_s}=0$, the combination  
$\frac{L}{p_{j_s}}\binom{p_{j_s}+m_{j_s}-1}{m_{j_s}}$ is 
to be understood as $\frac{L}{m_{j_s}}$.
Note that $\Omega(m) = \det F$ given in (\ref{eq:detf}) 
if $m_j = 1$ for all $j \in H$.

\begin{theorem}[\cite{KN} Theorems 3.5, 4.9]\label{th:kn}
\begin{align}
\Omega(m) &= \vert \, \U(m) \vert \quad 
m \in \M,\label{eq:uo}\\
\sum_{m \vdash M}\Omega(m) &= \binom{L}{M} \quad 
0 \le M \le L/2,\label{eq:osum}
\end{align}
where the sum extends over $m_1,m_2, \ldots \in \Z_{\ge 0}$
such that $\sum_jjm_j = M$.
\end{theorem}
The derivation of (\ref{eq:uo}) is due to the M\"obius 
inversion trick.
Actually (\ref{eq:osum}) is known to hold for 
any $M \ge 0$ if the symbol $\binom{\alpha}{\beta}$
is interpreted as the generalized binomial coefficient
$\alpha(\alpha-1)\cdots(\alpha-\beta+1)/\beta!$. 
The combined identity 
$\sum_m \vert \, \U(m) \vert = \binom{L}{M}$ is 
called the combinatorial completeness of the 
string hypothesis at $q=0$.

We include the formulas needed here and in 
Sections \ref{subsec:up} and \ref{subsec:dp}.
Fix $l \in \Z_{\ge 1}$. 
\begin{align}
\det A & = (\det F) 
\prod_{j \in H}(p_j+m_j)^{m_j-1},\label{eq:deta}\\
\det A[k\beta] & = (\det F[k])
\prod_{j \in H}(p_j+m_j)^{m_j-1}\;(k \in H,\, 1 \le \beta \le m_k),
\label{eq:detaj}\\
\det F &= L p_{j_1}p_{j_2}\cdots p_{j_{s-1}},\label{eq:detf}\\
\det F[j_{n+1}]- \det F[j_n] &= 
\frac{p_{i_s}(i_{n+1}-i_n)}
{p_{i_{n+1}}p_{i_n}} \det F \quad (0 \le n \le s-1).\label{eq:ff}
\end{align}
Here the matrix $A[k\beta]$ is obtained from 
$A=(A_{j\alpha, k\beta})$ (\ref{eq:a}) by replacing 
the $k \beta$ th column by $\vec{h}$.
In (\ref{eq:ff}), we are using the notation 
$i_0 =0$ and $i_n = \min(l,j_n)$ for $1\le n \le s$.
$p_{i_0}=p_0=L$ as noted after (\ref{eq:mh}).
The matrix $F[k]$ is 
obtained from $F=(F_{j,k})$ (\ref{eq:f}) by replacing 
the column for $k\in H$ by the column vector $\vec{h}'$.
(We understand $F[j_0]=0$.)
Here the vectors $\vec{h}$ and $\vec{h}'$ are given by
\begin{equation}\label{eq:hvec}
\vec{h} = (\min(j,l))_{j\in H, \,1 \le \alpha\le m_j} 
\in \Z^{m_{j_1}+\cdots +m_{j_s}},
\quad
\vec{h}' = (\min(j,l))_{j\in H} \in \Z^s,
\end{equation}
which are dependent on $l$.
Note that (\ref{eq:ff}) is vanishing unless $0 \le n \le t$, 
where $t$ is the maximum integer 
such that $i_{t+1}>i_t$. 
By the definition $0 \le t \le s\!-\!1$.
The relations (\ref{eq:deta}) and (\ref{eq:detf}) are derived by noting 
\begin{equation}\label{eq:suml}
\sum_{k\in H, 1 \le \beta \le m_k}A_{j\alpha, k\beta} 
= \sum_{k \in H}F_{j,k} = L.
\end{equation}
See also eq. (3.8) and (3.10) in \cite{KN}.
Combining (\ref{eq:deta}) and (\ref{eq:detaj}), we find
\begin{equation}\label{eq:ratio}
\frac{\det A[j\alpha]}{\det A} =
\frac{\det F[j]}{\det F} \quad (j \in H,\, 1 \le \alpha \le m_j).
\end{equation}
The property $\det\!A\gneqq 0$ is obvious from 
(\ref{eq:deta}) and (\ref{eq:detf}) 
under the condition $m \in \M$.

\subsection{Solitons as strings}\label{subsec:up}
Let us uncover the origin of our 
inverse scattering formalism in the Bethe ansatz.
Note that the action variable $m \in \M$ has 
emerged in two independent contexts in $\J(m)$ 
and $\U(m)$.
In the former it represents the soliton content of a path 
whereas in the latter it is the string content of the Bethe equation.
The key to their link is the map:
\begin{equation}\label{eq:psi}
\begin{split}
\Psi:\;\; \Jb(m) & \;\; \longrightarrow \;\; \U(m)\\
J \;\; \; & \;\;\longmapsto \;\; \;\; \vec{u},
\end{split}
\end{equation}
where $\vec{u}=(u^{(j)}_{\alpha})_{j \in H, 1 \le \alpha \le m_j}$ 
is determined 
from $J = (J^{(j)}_i)_{ i \in \Z,j \in H} \in \Jb(m)$
as the solution of the linear equation:
\begin{equation}\label{eq:sce2}
\sum_{k \in H}\sum_{\beta=1}^{m_k} 
A_{j\alpha,k\beta} u^{(k)}_{\beta} = 
\frac{1}{2}(p_j + m_j + 1) + J^{(j)}_\alpha+\alpha-1\quad 
(j \!\in \!H, \; 1 \!\le \!\alpha \!\le\! m_j).
\end{equation}
We write it as $A\vec{u} = \vec{c} + \vec{J} + \vec{\rho}$,
where $\vec{c}=(c^{(j)}_{\alpha})$ 
corresponding to $(p_j \!+\! m_j\! +\! 1)/2$ is a constant vector 
having the components that are independent of the index $\alpha$.
$\vec{\rho} = (\rho^{(j)}_{\alpha} = \alpha-1)$ 
is also a constant 
vector\footnote{The vector $\vec{c}$ plays a role only in 
Proposition \ref{pr:la}, and
$\alpha-1$ in (\ref{eq:sce2}) 
can be replaced with $\alpha+ a$ for any 
integer $a$.}.
The equation (\ref{eq:sce2}) 
is the string center equation (\ref{eq:sce}) that corresponds to 
a prescribed logarithmic branch of the Bethe equation.
To make sense of (\ref{eq:psi}), 
we are to check the off-diagonal condition 
(\ref{eq:odc}) on $\vec{u}$.
To do this it is useful to grasp the structure of the matrix $A$.
For example in case $H=\{1,2,3\}$, it looks as
(${\mathcal P}_j = p_j + m_j$ for short)
{\small
\begin{displaymath}
\left(
\begin{array}{cccccccccccc}
{\mathcal P}_1+1 &  \cdots & 1 &
2 & \cdots & 2 &
2 & \cdots & 2 \\
\vdots &\ddots  & \vdots &
\vdots &   & \vdots  &
\vdots &   & \vdots  \\
1 & \cdots  & {\mathcal P}_1+1 &
2 & \cdots & 2 &
2 & \cdots & 2 \\
2 & \cdots & 2 &
{\mathcal P}_2+3 & \cdots & 3 &
4 & \cdots & 4 \\
\vdots &   & \vdots  &
\vdots & \ddots & \vdots &
\vdots &   & \vdots  \\
2 & \cdots & 2 &
3 & \cdots & {\mathcal P}_2+3 &
4 & \cdots  & 4 \\
2 & \cdots & 2 &
4 & \cdots & 4 &
{\mathcal P}_3+5 &  \cdots & 5 \\
\vdots &   & \vdots  &
\vdots & \ddots  & \vdots  &
\vdots &   & \vdots \\
2 & \cdots  & 2 &
4 & \cdots  & 4 &
5 & \cdots & {\mathcal P}_3+5
\end{array}\right),
\end{displaymath}
}
which consists of 9 sub-matrices of size $m_i \times m_j$
($1 \le i,j \le 3$).
Thus (\ref{eq:sce2}) leads to 
$(p_j+m_j)(u^{(j)}_{\alpha}-u^{(j)}_{\beta}) 
= J^{(j)}_\alpha-J^{(j)}_\beta+\alpha-\beta$. 
{}From this and the condition (\ref{eq:jcomp}), we confirm that 
$\Psi$ produces a specific array of 
real numbers $(u^{(j)}_\alpha)$ that obey
$a\le u^{(j)}_1< u^{(j)}_2 < \cdots < u^{(j)}_{m_j} < a\!+\!1$ 
for some $a \in \R$, 
which may indeed be viewed as a representative 
element of $\U(m)$.

Varying $J \in \Jb(m)$ to $J'$ changes the image $\vec{u} = \Psi(J)$
to $\vec{u}' = \Psi(J')$ as an array of real numbers.
But $\vec{u}=(u^{(j)}_\alpha)$ and 
$\vec{u}'=(u'^{(j)}_\alpha)$ can be regarded as 
the same element in $\U(m)$  
under the identification by ${\mathfrak S}_{m_j}$ in (\ref{eq:odc})
\footnote{A proper treatment of this is to introduce 
the set of array of real numbers 
$\overline{\U}(m)$ that projects onto $\U(m)$
under the identification scheme 
in an analogous way from $\Jb(m)$ to $\J(m)$.  
However we skip it here 
supposing no confusion might arise.}.
{}From the above property of $\Psi$, such an event
takes place if and only if 
$\vec{u}' = 
{\tilde \sigma}^{n_1}_{j_1}\cdots {\tilde \sigma}^{n_s}_{j_s}(\vec{u})$
for some $n_1, \ldots, n_s \in \Z$, 
where ${\tilde \sigma}_k$ is defined by
\begin{equation*}
{\tilde \sigma}_k:\; 
(u^{(j)}_1, u^{(j)}_2, \ldots, u^{(j)}_{m_j}) 
\longmapsto 
\begin{cases}
(u^{(j)}_2, \ldots, u^{(j)}_{m_j}, u^{(j)}_1 + 1) & \hbox{ if } j=k,\\
(u^{(j)}_1, u^{(j)}_2, \ldots, u^{(j)}_{m_j}) & \hbox{ if } j\neq k.
\end{cases}
\end{equation*}
Since the matrix $A$ is invertible, the effect of ${\tilde \sigma}_k$
is translated to that on $J \in \Jb(m)$.
\begin{lemma}\label{lem:tslide}
If $\vec{u}=\Psi(J)$, then 
${\tilde \sigma}_k(\vec{u}) = \Psi(\sigma_k(J))$,
where $\sigma_k$ is the slide defined in (\ref{eq:slide}).
\end{lemma}
\begin{proof}
Set $\vec{u}'={\tilde \sigma}_k(\vec{u}) = \Psi(J')$.
By the definition we have 
$A\vec{u} = \vec{c} + \vec{J} + \vec{\rho}$
and $A\vec{u}' = \vec{c} + \vec{J}' + \vec{\rho}$.
In view of the structure of $A$, we find for $j \neq k$ that
$J'^{(j)}_\alpha = J^{(j)}_\alpha + A_{j\alpha, k m_k}
= J^{(j)}_\alpha + 2\min(j,k)$.
On the other hand, the $k$ th block 
of $\vec{c} + \vec{J}' + \vec{\rho}=A\vec{u}'$
is evaluated as ($c=c^{(k)}_1$) 
\begin{equation*}
\begin{split}
&
\begin{pmatrix}
c\!+\!J'^{(k)}_1 \\
c\!+\!J'^{(k)}_2 + 1\\
\vdots \\
c\!+\!J'^{(k)}_{m_k}\!+\!m_k\! - \!1
\end{pmatrix} \\
&=
\begin{pmatrix}
p_k\!+\!m_k\!+\!2k\!-\!1 & \!\!\!\cdots & \!\!\!2k\!-\!1 \\
\vdots  & \ddots & \vdots \\
\\
2k\!-\!1 & \!\!\!\cdots & \!\!\!p_k\!+\!m_k\!+\!2k\!-\!1
\end{pmatrix}\!\!
\begin{pmatrix}
u^{(k)}_2 \\
\vdots \\
u^{(k)}_{m_k}\\
u^{(k)}_1\!+\!1
\end{pmatrix}+
\sum_{j(\neq k), 1 \le \beta \le m_j}\!\!
\begin{pmatrix}
A_{k1, j\beta}\\
A_{k2, j\beta}\\
\vdots \\
A_{km_k,j\beta}
\end{pmatrix}\!u^{(j)}_\beta \\
&=
\begin{pmatrix}
c\!+\!J^{(k)}_2+1+2k-1 \\
\vdots \\
c\!+\!J^{(k)}_{m_k}\!+\!m_k\! - \!1+2k-1\\
c\!+\!J^{(k)}_1+p_k+m_k+2k-1
\end{pmatrix},
\end{split}
\end{equation*} 
where the last equality is due to 
the $k$ th block of 
the relation 
$A\vec{u} = \vec{c} + \vec{J} + \vec{\rho}$.
Noting the quasi-periodicity $J^{(k)}_1+p_k = J^{(k)}_{m_k+1}$,
this result is expressed as $J'^{(k)}_i = J^{(k)}_{i+1}+2k$
for $1 \le i \le m_k$.
Unifying the formulas we find that $J'=(J'^{(j)}_i)$ is 
obtained from $J=(J^{(j)}_i)$ by 
$J'^{(j)}_i = J^{(j)}_{i+\delta_{j,k}}+2\min(j,k)$
in agreement with (\ref{eq:slide}).
\end{proof}
We have established that 
$\Psi(J)$ and $\Psi(J')$ represent the same element 
in $\U(m)$ if and only if $J \simeq J' \in \Jb(m)$.
{}From the definition (\ref{eq:j}), we obtain
\begin{theorem}\label{th:ju}
The map $\Psi$ (\ref{eq:psi}) induces the bijection 
between the set of angle variables $\J(m)$ 
and the set of off-diagonal solutions $\U(m)$ to the 
string center equation (\ref{eq:sce}).
\end{theorem}

The induced bijection will also be denoted by $\Psi$.
Now we are able to introduce the time evolution $T_l$ 
of the Bethe root $\vec{u} = \Psi(J) \in \U(m)$ by
$T_l(\vec{u}) = \Psi(T_l(J))$ using (\ref{eq:tlj}).

\begin{corollary}\label{cor:pju}
Theorem \ref{th:main} is extended to 
the following commutative diagram:
\begin{equation}\label{eq:cd2}
\begin{CD}
{\mathcal P}(m) @>{\Phi}>> \J(m) @>{\Psi}>> \U(m) \\
@V{T_l}VV @V{T_l}VV @V{T_l}VV\\
{\mathcal P}(m) @>{\Phi}>> \J(m) @>{\Psi}>> \U(m)
\end{CD}
\end{equation}
Their cardinality is given by
\begin{equation}\label{eq:cpju}
\vert \Pth(m) \vert = \vert \J(m) \vert = \vert \U(m) \vert = 
\Omega(m).
\end{equation}
\end{corollary}
The result $\vert \Pth(m) \vert = \Omega(m)$ 
endows the character formula in \cite{KN} with the
quasi-particle interpretation.
The quantity $\Omega(m)$ originally associated 
with the string content $m$ turns out to be 
the number of states in the periodic box-ball system 
having the prescribed soliton content $m$.
Then the identity (\ref{eq:osum}) implies that weight spaces 
of the $sl_2$ module $({\mathbb C}^2)^{\otimes L}$ 
are decomposed into subspaces spanned by iso-energy
states in the periodic box-ball system.  
These observations and Proposition \ref{pr:em} 
are quantitative supports to 
identify strings with solitons.
The identification also agrees with the physical picture 
\cite{Be} that 
the strings represent bound states of magnons
over the ferromagnetic vacuums $111 \ldots 111$ 
or $222 \ldots 222$. 
The fact $\vert \Pth(m) \vert = \Omega(m)$ 
has also been shown in \cite{YYT} by a different approach.

\begin{example}\label{ex:omega}
Consider the length $L=8$ paths $p$ with $\wt(p)=0$.
There are $\binom{8}{4}=70$ such paths, which are grouped into
$\Pth(m)$ according to the 
5 elements $m=(m_1,m_2,\ldots) \in \M$ (\ref{eq:m}).
Their cardinality $\vert \Pth(m) \vert = \Omega(m)$ 
are given as follows:

\begin{center}
\begin{tabular}{c|c|c|c|c|c}
$m\!=\!(m_1,m_2,m_3,m_4)$ 
& $(4,0,0,0)$ & $(2,1,0,0)$ & $(0,2,0,0)$ & $(1,0,1,0)$ & $(0,0,0,1)$ \\
\hline 
&&&&&\vspace{-0.2cm}\\
$\Omega(m)$ 
& $2$ & $24$ & $4$ & $32$ & $8$
\end{tabular}
\end{center}

\vspace{0.2cm}
\noindent
The paths in $\Pth(m)$ are listed as follows: 
\vspace{0.2cm}

\begin{center}
\begin{tabular}{c|l}
$(m_1,m_2,m_3,m_4)$ & $\Pth(m)$ \\
\hline & \vspace{-0.2cm}\\
(4,0,0,0)& 12121212, \; 21212121\\ & \vspace{-0.3cm}\\
(2,1,0,0)& $T_1^n(12121122), \;T_1^n(12112122), \;T_1^n(11212122)$ \\
& \vspace{-0.3cm}\\
(0,2,0,0)& 11221122, \; 21122112, \; 22112211, \; 12211221\\
& \vspace{-0.3cm}\\
(1,0,1,0)& $T_1^n(12111222),  \; T_1^n(11211222), \; T_1^n(11121222),
\; T_1^n(11122122)$\\
& \vspace{-0.3cm}\\
(0,0,0,1)& $T_1^n(11112222)$
\end{tabular}
\end{center}

\vspace{0.2cm}\noindent
where $T^n_1(p)$ 
stands for the 8 paths $(n \in \Z_8)$ obtained by cyclic
shifts of $p$.
The path $11122122$ for $m=(1,0,1,0)$ 
is the intermediate stage of the 
collision of the length 3 and 1 solitons.
\end{example}

\subsection{Periodicity}\label{subsec:dp}

The time evolution $T_l$ is invertible as seen in (\ref{eq:tinv})
and the set of states $\Pth$ is finite.
Therefore every path $p \in \Pth$ possesses the 
property $T_l^N(p)=p$ for some integer $N\ge 1$.
We say any such integer a {\em period} of $p$.
The minimum period is called the 
{\em fundamental period} of $p$ and denoted by 
${\mathcal N}^\ast= {\mathcal N}^\ast_l(p)$.
Every period is a multiple of the fundamental period
${\mathcal N}^\ast$.
Here we establish a formula for 
the fundamental period under any $T_l$
taking advantage of the linearization scheme (\ref{eq:cd2}).
We assume $p \in \sqcup_{m\in \M}\Pth(m)$ with no loss
of generality thanks to (\ref{eq:pd}) and Proposition \ref{pr:w}.

For $p \in \Pth(m)$, let $J=\Phi(p) \in \J(m)$.
The fundamental period ${\mathcal N}^\ast$ of $p$ 
under $T_l$ is a function of $l$, the action variable $m$ and 
the angle variable $J$.
To be expository, we approach 
the fundamental period in two steps.
First we show 
in Theorem \ref{th:dp} that there exists a value 
${\mathcal N} = {\mathcal N}_l(m)$ that is independent of 
$J$ nonetheless making  
$T^{\mathcal N}_l(p)=p$ hold for all $p \in \Pth(m)$. 
We call it the {\em generic period} 
inherent to the whole $\Pth(m)$.
In general ${\mathcal N}$ is yet a multiple of 
the fundamental period ${\mathcal N}^\ast$.
Second we analyze the accidental symmetry 
gained by special $J$ that makes  
${\mathcal N}^\ast$ a divisor of ${\mathcal N}$.
It improves the idea of Theorem \ref{th:dp} and leads to the  
final result in Theorem \ref{th:fp}. 
We begin by presenting the formula for the 
generic period ${\mathcal N}$.
\begin{equation}\label{eq:lcm}
{\mathcal N} = 
{\rm LCM}\!\Bigl(1,
\bigcup_{j\alpha}{}^\prime 
\frac{\det \!A}{\det \!A[j\alpha]}\Bigr)
= {\rm LCM}\!\Bigl(1,
\bigcup_{j}{}^\prime 
\frac{\det \!F}{\det \!F[j]}\Bigr),
\end{equation}
where the matrices $A, A[j\alpha]$ and $F, F[j]$ are defined 
in (\ref{eq:a}), (\ref{eq:f}) and after (\ref{eq:ff}).
In particular, $A[j\alpha]$ and $F[j]$ are dependent on $l$. 
The union $\bigcup_{j\alpha}'$ extends over 
those $j \in H, \,1 \le \alpha \le m_j$ such that 
$\det A[j\alpha]\neq 0$.
The union $\bigcup_j'$ does over
those $j \in H$ such that $\det F[j]\neq 0$.
The least common multiple ${\rm LCM}(r_1, \ldots, r_n)$ 
of nonzero rational numbers $r_1,\ldots, r_n$ 
is the minimum positive rational number 
in $\Z r_1 \cap \cdots \cap \Z r_n$.
Thus ${\rm LCM}(1, r_1, \ldots, r_n)$ 
denotes the minimum positive integer  
in $\Z r_1 \cap \cdots \cap \Z r_n$.
The equality of the two expressions is due to (\ref{eq:ratio}).

\begin{theorem}\label{th:dp}
For any path $p \in \Pth(m)$, 
$T^{\mathcal N}_l(p) = p$ is valid.
\end{theorem}

To be instructive we present two proofs utilizing 
$\U(m)$ and $\J(m)$ in the linearization scheme 
(\ref{eq:cd2}) although they are essentially the same.
In both proofs we first consider a sufficient condition for
$T^{\mathcal N}_l(p) = p$ by regarding ${\mathcal N}$ as
an unknown, and confirm afterwards that the choice $(\ref{eq:lcm})$ 
fulfills it.
The identities and definitions (\ref{eq:deta})--(\ref{eq:hvec})
work effectively.

\begin{proof}[Proof based on $\U(m)$]
Suppose 
$p \overset{\Phi}{\mapsto} J \overset{\Psi}{\mapsto} \vec{u}$
and 
$T_l^{\mathcal N}(p) \overset{\Phi}{\mapsto} J' 
\overset{\Psi}{\mapsto} \vec{u}'$.
We have $A\vec{u} = \vec{c}+\vec{J}+\vec{\rho}$ and 
$A\vec{u}' = \vec{c}+\vec{J'}+\vec{\rho}$.
{}From (\ref{eq:tlj}) and (\ref{eq:hvec}), this leads to 
$A(\vec{u}' - \vec{u}) = \vec{J'}-\vec{J}={\mathcal N}\vec{h}$.
A sufficient condition for $T^{\mathcal N}_l(p) = p$ is that 
all the components of 
$\vec{v}:=\vec{u}'-\vec{u}$ are integers.
Since $\vec{v} = {\mathcal N}A^{-1}\vec{h}$,
its $j\alpha$ th component is 
$v^{(j)}_{\alpha} = {\mathcal N}\det A[j\alpha]/\det A$.
Therefore the ${\mathcal N}$ in (\ref{eq:lcm}) satisfies the  
condition.
\end{proof}

\begin{proof}[Proof based on $\J(m)$]
Take $J \in \Jb(m)$ such that $[J] = \Phi(p) \in \J(m)$. 
We are to consider a sufficient condition 
for $[T^{\mathcal N}_l(J)] = [J]$.
This indeed happens if there exist integers 
$n_1, \ldots, n_s$ such that 
$\sigma^{n_1}_{j_1}\cdots \sigma^{n_s}_{j_s}(J) 
=  T^{\mathcal N}_l(J) \in \Jb(m)$.
We look for them by further presuming that 
$r_a = n_a/m_{j_a}$ is an integer.
Then from (\ref{eq:jcomp}) and (\ref{eq:slide}), 
the above relation on the block $\J_{j_k}$ is expressed as  
\begin{equation}\label{eq:r}
p_{j_k} r_k + 
2\sum_{a=1}^s\min(j_k,j_a)m_{j_a}r_a = {\mathcal N}\min(j_k,l).
\end{equation}
In terms of the column vectors 
$\vec{r} = {}^t(r_1,\ldots, r_s)$, 
$\vec{h}'$ (\ref{eq:hvec}) and 
the matrix $F$ (\ref{eq:f}), 
this is written as $F\vec{r} = {\mathcal N} \vec{h}'$.
Its solution is given as $r_a=  {\mathcal N}\det F[j_a]/\det F$
by using the matrix $F[k]$ appearing in (\ref{eq:detaj}).
Thus the ${\mathcal N}$ in (\ref{eq:lcm})
assures the existence of $\forall r_a \in \Z$.
\end{proof}

The expressions (\ref{eq:lcm}) can be simplified.
We employ the notation used in (\ref{eq:ff}), e.g., 
$H=\{j_1,\ldots, j_s\}$ as in (\ref{eq:mh}),
$i_n=\min(j_n,l)\, (1 \le n \le s), i_0=0$ and 
$0 \le t \le s-1$ is the maximum integer 
such that $i_{t+1}>i_t$.
Since $F[j_{t+1}]=F[j_{t+2}]= \cdots = F[j_s]$ from (\ref{eq:ff}), 
we have
\begin{equation}\label{eq:lcm2}
\begin{split}
{\mathcal N} &={\rm LCM}
\left(1,\frac{\det F}{\det F[j_1]}, \frac{\det F}{\det F[j_2]}, 
\ldots, 
\frac{\det F}{\det F[j_{t+1}]}\right)\\
&= {\rm LCM}
\left(1,\frac{\det F}{\det F[j_1]}, \frac{\det F}{\det F[j_2]-\det F[j_1]}, 
\ldots, 
\frac{\det F}{\det F[j_{t+1}]-\det F[j_t]}\right)\\
&= {\rm LCM}\!\left(1, \bigcup_{n=0}^{t}{}^\prime 
\frac{p_{i_{n+1}}p_{i_n}}{(i_{n+1}-i_n)p_{i_s}}\right),
\end{split}
\end{equation}
where the second equality is an elementary property of ${\rm LCM}$.
The last line is obtained by (\ref{eq:ff}). 
The union $\cup_n'$ extends over those $n$ such that 
$\frac{p_{i_{n+1}}p_{i_n}}{(i_{n+1}-i_n)p_{i_s}}$ is finite.
This caution is needed only if $\wt(p)=0$ and $l \ge j_s$,
where one encounters 
$p_{i_s}=p_{j_s}=0$ in the denominator.
In this case $t=s-1$ and 
only $n=s-1$ is allowed in the union so as to cancel the zero
by $p_{i_{n+1}}=0$ in the numerator.
This leads to ${\mathcal N} = 
{\rm LCM}(1,\frac{p_{j_{s-1}}}{j_s-j_{s-1}})$.
Noting that $p_{j_{s-1}} = p_{j_{s-1}}-p_{j_s}
= 2(j_s-j_{s-1})m_{j_s}$, we get ${\mathcal N} = 2m_{j_s}$.
Thus $T^{\mathcal N}_l(p)=p$ is certainly valid, 
for Proposition \ref{pr:ts} and the 
remark following it imply that $T^2_l(p)=p$.
Another simplification of (\ref{eq:lcm}) 
occurs at $l=1$, where $i_n=1(n>0)$ hence $t=0$.
In this case we have ${\mathcal N}=L$, which is again consistent
with $T^L_1(p)=p$ for any path $p$.

Being able to simplify (\ref{eq:lcm}) into (\ref{eq:lcm2}) is gratifying.
However as we will see in (\ref{eq:ncond}) 
in Section \ref{subsec:ana}, 
it is instead the expression (\ref{eq:lcm})  
that makes the generic period ${\mathcal N}$ conceptual and 
elucidates the essence of the game.
We note that the formula (\ref{eq:lcm}) 
is the simplest case of the most general 
one for $A^{(1)}_n$ conjectured in eq.(8) in \cite{KT2}.
The expression (\ref{eq:lcm2}) with $l=\infty$ was
first obtained in \cite{YYT}.

\begin{example}\label{ex:period}
Take the path $p=21121111221122111111222$ of 
length $L = 23$. It belongs to $\Pth(m)$ with 
$m = \mu(p) = (m_1,m_2,m_3,m_4) = (1,2,0,1)$ hence
$H = \{1,2,4\}$ and $s=3$ in (\ref{eq:mh}).
The vacancy numbers (\ref{eq:va}) are
$(p_1,p_2,p_3,p_4) = (15,9,7,5)$.
For $1 \le l \le 4$, the data $i_1, i_2, i_3$ and $t$ are listed as follows:
\begin{center}
\begin{tabular}{c|cccc}
$l$ & $i_1$ & $i_2$ & $i_3$ & $t$ \\
\hline 
1 & 1 & 1 & 1 & 0 \\
2 & 1 & 2 & 2 & 1 \\
3 & 1 & 2 & 3 & 2 \\
4 & 1 & 2 & 4 & 2
\end{tabular}
\end{center}
By convention we also have $i_0 = 0$ and $p_0 = L = 23$.
According to the last expression in (\ref{eq:lcm2}),
${\mathcal N}$ is calculated as
\begin{align*}
l=1 & : {\rm LCM}(1, p_0) = {\rm LCM}(1, 23) = 23,\\
l=2 & : {\rm LCM}\Bigl(1, \frac{p_1p_0}{p_2}, p_1\Bigr) 
= {\rm LCM}\Bigl(1, \frac{115}{3}, 15\Bigr) = 345,\\
l=3 & : {\rm LCM}\Bigl(1, \frac{p_1p_0}{p_3}, \frac{p_2p_1}{p_3}, p_2\Bigr) 
= {\rm LCM}\Bigl(1, \frac{345}{7}, \frac{135}{7}, 9\Bigr)
=3105,\\
l=4 & : {\rm LCM}\Bigl(1, \frac{p_1p_0}{p_4}, \frac{p_2p_1}{p_4},
\frac{p_2}{2}\Bigr) = {\rm LCM}\Bigl(1, 69, 27,\frac{9}{2}\Bigr) = 621.
\end{align*}
Actually these generic period ${\mathcal N}$ 
coincide with the fundamental period ${\mathcal N}^\ast$ of $p$ 
under the respective time evolutions $T_l$.  
The reason will be explained in (\ref{eq:nns}).
\end{example}

\begin{example}\label{ex:deft}
Take the path 
$p=2122112211221111222111122$ of length $L=25$.
This is the same path as that considered in Example \ref{ex:neqn}.
It belongs to $\Pth(m)$ with 
$m= \mu(p) = (m_1,m_2,m_3,m_4)=(1,2,1,1)$.
The vacancy numbers read
$(p_1,p_2,p_3,p_4) = (15,7,3,1)$.
This time we calculate ${\mathcal N}$ along (\ref{eq:lcm}).
The matrix $F$ (\ref{eq:f}) reads
\begin{equation*}
F = \begin{pmatrix}
p_1+2m_1 & 2m_2 & 2m_3 & 2m_4\\
2m_1 & p_2+4m_2 & 4m_3 & 4m_4\\
2m_1 & 4m_2 & p_3+6m_3 & 6m_4\\
2m_1 & 4m_2 & 6m_3 & p_4+8m_4
\end{pmatrix}
= \begin{pmatrix}
17 & 4 & 2 & 2 \\
 2 & 15 & 4 & 4 \\
 2 & 8 & 9 & 6 \\
 2 & 8 & 6 & 9
\end{pmatrix}.
\end{equation*}
If $l=3$ for instance, the matrix $F[j]$ is obtained 
from $F$ by replacing its $j$ th column with 
$\vec{h}' = {}^t(1,2,3,3)$ (\ref{eq:hvec}). 
Thus $F[1], F[2], F[3]$ and $F[4]$ look as
\begin{equation*}
\begin{pmatrix}
1 & 4 & 2 & 2 \\
 2 & 15 & 4 & 4 \\
 3 & 8 & 9 & 6 \\
 3 & 8 & 6 & 9
\end{pmatrix},
\;\;
\begin{pmatrix}
17 & 1 & 2 & 2 \\
 2 & 2 & 4 & 4 \\
 2 & 3 & 9 & 6 \\
 2 & 3 & 6 & 9
\end{pmatrix},
\;\;
\begin{pmatrix}
17 & 4 & 1 & 2 \\
 2 & 15 & 2 & 4 \\
 2 & 8 & 3 & 6 \\
 2 & 8 & 3 & 9
\end{pmatrix},
\;\;
\begin{pmatrix}
17 & 4 & 2 & 1 \\
 2 & 15 & 4 & 2 \\
 2 & 8 & 9 & 3 \\
 2 & 8 & 6 & 3
\end{pmatrix},
\end{equation*}
which lead to 
\begin{equation*}
\Bigl(1, \frac{\det F}{\det F[1]}, 
\frac{\det F}{\det F[2]}, 
\frac{\det F}{\det F[3]}, 
\frac{\det F}{\det F[4]} \Bigr) = 
\Bigl(1, 125, \frac{875}{32}, \frac{875}{157}, \frac{875}{157} \Bigr).
\end{equation*}
Taking the ${\rm LCM}$ of them, we get ${\mathcal N} =875$.
We list the result of such calculations for $1 \le l \le 4$
in the following table.

\vspace{0.2cm}
\begin{center}
\begin{tabular}{c|ccccc|c}
$l$ & \multicolumn{5}{|c|}{LCM} &  ${\mathcal N}$  \\ \hline
\vspace{-0.3cm} &&&&&&\\
1 & 1,& 25,& 25,& 25,& 25 & 25 \\
\vspace{-0.2cm} &&&&&&\\
2 & 1, &$\frac{375}{7},$ & $\frac{375}{32},$ & $\frac{375}{32},$ &
$\frac{375}{32}$ & 375 \\
\vspace{-0.2cm} &&&&&&\\
3 & 1,  &125, & $\frac{875}{32},$ & $\frac{875}{157},$ &
$\frac{875}{157}$ & 875\\
\vspace{-0.2cm} &&&&&&\\
4 & 1, & 375, & $\frac{2625}{32},$ & $\frac{2625}{157},$ &
$\frac{875}{344}$ & 2625
\end{tabular}
\end{center}
\vspace{0.2cm}

For each $l$, $T^{\mathcal N}_l(p) = p$ holds.
Actually the ${\mathcal N}$ listed here coincides with 
the fundamental period of $p$. 
See (\ref{eq:nns}) for the reason.
\end{example}

As in Examples \ref{ex:period} and \ref{ex:deft}, 
the generic period ${\mathcal N}={\mathcal N}_l(m)$ (\ref{eq:lcm})
turns out to be the fundamental period ${\mathcal N}^\ast_l(p)$ 
for a majority of paths $p \in \Pth(m)$. 
The event 
${\mathcal N}_l(m)/{\mathcal N}_l^\ast(p) \in \Z_{\ge 2}$
indicates some extra symmetry in the path $p$.
The structure of 
our angle variable provides a lucid picture on 
the nature of the extra symmetry, and 
improves (\ref{eq:lcm}) into a formula for  
the fundamental period.

Consider the block $J=(J_i)_{i \in \Z} \in \Jb_j(m)$ (\ref{eq:jcomp})
for $j \in H=\{j_1, \ldots, j_s\}$.
Seek the maximum positive integer $g=g_j$ such that 
\begin{equation}\label{eq:afold}
\frac{p_j}{g}, \; \frac{m_j}{g} \in \Z_{\ge 1}, 
\quad J_{i+\frac{m_j}{g}}-J_i = \frac{p_j}{g} \quad\hbox{for any } i \in \Z,
\end{equation}
which is consistent with the condition 
$J_{i+m_j} = J_i + p_j$ in (\ref{eq:jcomp}).
We say that such $J$ has the order $g$ symmetry,
which means the finer quasi-periodicity 
than the original one by factor $g$.
The case $g=1$ corresponds to the previous treatment.

\begin{theorem}\label{th:fp}
For any path $p \in \Pth(m)$, let 
$(J^{(j)})_{j \in H} = \Phi(p)$ be the angle variable.
Suppose $J^{(j)}$ has the order $g_j$ symmetry. 
Then the fundamental period of $p$ is given by
\begin{equation}\label{eq:nn} 
{\mathcal N}^\ast = 
{\rm LCM}\!\Bigl(1,
\bigcup_{j \in H}\!{}^\prime \;\frac{\det \!F}{g_{j}\!\det \!F[j]}
\Bigr),
\end{equation}
where $\cup'_j$ is the union over those $j$ such that 
$\det \!F[j] \neq 0$.
\end{theorem}

Obviously ${\mathcal N} \ge {\mathcal N}^\ast$ and 
the equality holds if (but not only if) $g_j=1$ for all $j \in H$.

\begin{proof}
A slight modification of the argument around (\ref{eq:r}) suffices.
{}From the definition (\ref{eq:afold}) of the symmetry of $J$, 
to validate $T^N_l(p)=p$ or equivalently 
$\sigma^{n_1}_{j_1}\cdots \sigma^{n_s}_{j_s}(J) 
=  T^N_l(J)$, it is {\em necessary} 
for $r_a = n_a/(m_{j_a}/g_{j_a})$ to be an integer
{\em and} to satisfy 
\begin{equation*}
\frac{p_{j_k}}{g_{j_k}} r_k + 
2\sum_{a=1}^s\min(j_k,j_a)\frac{m_{j_a}}{g_{j_a}}r_a 
=N \min(j_k,l).
\end{equation*}
This is precisely (\ref{eq:r}) with $r_a$ replaced by $r_a/g_{j_a}$
and ${\mathcal N}$ by $N$.
Thus $N$ must be chosen so that 
$r_a = N g_{j_a}\det F[j_a]/\det F$ 
be an integer for all $1 \le a \le s$.  
The ${\mathcal N}^\ast$ in 
(\ref{eq:nn}) is the minimum of such $N$.
\end{proof}

Beside the generic case $\forall g_j =1$, the next simplest 
situation is $\forall g_j = g (\ge 2)$.
It only happens when $L/g \in \Z$, and 
corresponds to the path of the form 
$p=q^{\otimes g}$ for some 
$q \in B^{\otimes L/g}_1$ without any such symmetry. 
(This statement is justified from Lemma \ref{lem:waru}.) 
In this case, ${\mathcal N}^\ast$ (\ref{eq:nn}) coincides with 
the generic period ${\mathcal N}$ (\ref{eq:lcm}) for $q$,
which is consistent with $T_l(p)=T_l(q)^{\otimes g}$.

\begin{example}\label{ex:fp}
We consider the path $p$ treated in Example \ref{ex:tg}.
We know that $T^{130}_3(p)=p$, and actually $130$
is the fundamental period under $T_3$.
Looking at the angle variable 
$J = (J^{(1)}, J^{(2)}, J^{(3)}) = \Phi(p)$ given there,
we find that $J^{(1)}$ and $J^{(2)}$ possess 
the order 2 symmetry. Thus 
$(g_1,g_2,g_3)=(2,2,1)$ in the notation in (\ref{eq:nn}). 
We have $H=\{1,2,3\}$, 
$(m_1,m_2,m_3)=(2,2,1)$ and $(p_1,p_2,p_3) = (16,10,8)$.
Thus the matrix $F$ (\ref{eq:f})  reads
\begin{equation*}
F = \begin{pmatrix}
p_1+2m_1 & 2m_2 & 2m_3\\
2m_1 & p_2+4m_2 & 4m_3\\
2m_1 & 4m_2 & p_3+6m_3
\end{pmatrix}
= \begin{pmatrix}
20 & 4 & 2\\
4 & 18 & 4\\
4 & 8 & 14
\end{pmatrix}.
\end{equation*}
We consider the time evolution $T_3$, for which 
the vector $\vec{h}'$ (\ref{eq:hvec}) reads $\vec{h}' = {}^t(1,2,3)$.
Replacing the columns of $F$ with this, we get
\begin{equation*}
\begin{split}
F[1]& = \begin{pmatrix}
1 & 4 & 2\\
2 & 18 & 4\\
3 & 8 & 14
\end{pmatrix},\quad 
F[2] = \begin{pmatrix}
20 & 1 & 2\\
4 & 2 & 4\\
4 & 3 & 14
\end{pmatrix},\quad 
F[3] = \begin{pmatrix}
20 & 4 & 1\\
4 & 18 & 2\\
4 & 8 & 3
\end{pmatrix}.
\end{split}
\end{equation*}
Substituting the determinants 
$(\det F, \det F[1], \det F[2], \det F[3]) = (4160, 80, 288, 704)$ 
into the formulas 
for the generic period ${\mathcal N}$ (\ref{eq:lcm}) 
and the fundamental period 
${\mathcal N}^\ast$ (\ref{eq:nn}), we find $(l=3)$
\begin{align*}
{\mathcal N} & = {\rm LCM}
\Bigl(1,\frac{4160}{80}, \frac{4160}{288}, \frac{4160}{704}\Bigr) 
=  {\rm LCM}
\Bigl(1, 52, \frac{130}{9}, \frac{65}{11}\Bigr) = 260,\\
{\mathcal N}^\ast & = {\rm LCM}
\Bigl(1,\frac{4160}{2\times 80}, \frac{4160}{2\times 288}, 
\frac{4160}{704}\Bigr) 
=  {\rm LCM}
\Bigl(1, 26, \frac{65}{9}, \frac{65}{11}\Bigr) = 130.
\end{align*}
\end{example}

{}From (\ref{eq:lcm}), (\ref{eq:afold}) and (\ref{eq:nn}),
we find 
\begin{equation}\label{eq:nns}
{\mathcal N} = {\mathcal N}^\ast \quad
\hbox{ if } 
\;{\rm GCD}(p_j, m_j) = 1 \;\;\hbox{ for all }\; j \in H,
\end{equation}
where we employ the convention ${\rm GCD}(0,m_{j_s})=1$
for the greatest common divisor when $p_{j_s}=0$.
Note that (\ref{eq:nns}) is a sufficient
but not a necessary condition for ${\mathcal N} = {\mathcal N}^\ast$.
It explains the reason for ${\mathcal N} = {\mathcal N}^\ast$ 
in Examples \ref{ex:period} and \ref{ex:deft}.
For $l=\infty$, 
the fundamental period has also been studied in \cite{YYT}.

\subsection{Bethe eigenvalue}\label{subsec:be}

The time evolution $T_l$ in the periodic box-ball system 
is the $q=0$ limit of the row transfer matrix $T_l(\zeta)$.
Its eigenvalues are given by the analytic Bethe ansatz \cite{R,KS}.
Let
$Q(\theta)= \prod_{k=1}^M\sinh\pi(\theta-\sqrt{-1}u_k)$ be 
Baxter's $Q$-function, where 
$\{u_k\}$ satisfy the Bethe equation (\ref{eq:be}).
We set $q = e^{-2\pi\hbar}$ as before, 
$\zeta = e^{2\pi \theta}$ and assume $M \le L/2$.
For the string solution (\ref{eq:string}), the $q\rightarrow 0$ 
limit of the Bethe eigenvalue (\cite{KS} eq.(1.8)) coincides with 
that of the top term (\cite{KS} eq.(2.12)):
\begin{equation}\label{eq:la0}
\lim_{q \rightarrow 0}
\frac{Q(\theta-l\hbar)}{Q(\theta+l\hbar)}= 
\zeta^{-E_l}\Lambda_l,
\; \;\;
\Lambda_l=\exp\Bigl(
2\pi\sqrt{-1}\sum_{j\alpha}\min(j,l)(u^{(j)}_\alpha+\frac{1}{2})
\Bigr)
\end{equation}
under an appropriate normalization.
Here $E_l=\sum_k\min(l,k)m_k$ is the right hand side of (\ref{eq:em}), and
the sum $\sum_{j\alpha}$ 
extends over $j \in H$ and $1 \le \alpha \le m_j$.

\begin{proposition}\label{pr:la}
For any solution $\vec{u}=(u^{(j)}_\alpha)$ to the string center 
equation (\ref{eq:sce}), the equality 
$\Lambda^{\mathcal N}_l = 1$ is valid, where 
${\mathcal N}$ is the generic period (\ref{eq:lcm}).
\end{proposition}
\begin{proof}
Put $\vec{v} = (v^{(j)}_\alpha) = {\mathcal N} A^{-1}\vec{h}$,
where $\vec{h}$ is specified in (\ref{eq:hvec}).
Its component is expressed as 
$v^{(j)}_\alpha = {\mathcal N}\det A[j\alpha]/\det A
= {\mathcal N} \det F[j]/\det F$, where we have used 
(\ref{eq:ratio}). This is independent of the index $\alpha$.
Moreover from (\ref{eq:lcm}), we may set $v^{(j)}_\alpha = x_j$ 
for some $x_j \in \Z$.
Further define the vectors $\vec{c} = (c^{(j)})_{j\alpha}$ with
$c^{(j)} = (p_j+m_j+1)/2$, and $\vec{e}=(\frac{1}{2})_{j\alpha}$. 
Then the string center equation (\ref{eq:sce}) is written as
$A \vec{u} = \vec{I}+\vec{c}$ for some integer vector $\vec{I}$.
The $\Lambda_l$ (\ref{eq:la0}) is expressed as
$\Lambda_l = \exp(2\pi\sqrt{-1}\,{}^t\vec{h}(\vec{u}+\vec{e}))
= \exp(2\pi\sqrt{-1}\,{}^t\vec{h}A^{-1}
(\vec{I}+\vec{c}+A\vec{e}))$.
Since $A$ is a symmetric matrix, we have 
${\mathcal N}({}^t\vec{h} A^{-1})_{j\alpha} 
= v^{(j)}_\alpha=x_j$.
Thus 
$\Lambda_l^{\mathcal N} = 
\exp(2\pi\sqrt{-1}\sum_{j\alpha}x_j
(\vec{I}+\vec{c}+A\vec{e})_{j\alpha})
= \exp(2\pi\sqrt{-1}\sum_{j\alpha}x_j
(\vec{c}+A\vec{e})_{j\alpha})$.
Now $(\log\Lambda_l^{\mathcal N})/(2\pi\sqrt{-1})$ is 
calculated using (\ref{eq:a}) as
\begin{equation*}
\begin{split}
\sum_{j\alpha}x_j(\vec{c}+A\vec{e})_{j\alpha}
&=\frac{1}{2}\sum_{j\alpha}x_j\bigl(p_j+m_j+1+
\sum_{k\beta}A_{j\alpha,k\beta}\bigr)\\
&=\frac{1}{2}\sum_{j\alpha}x_j\bigl(p_j+m_j+1+
p_j+m_j+2\sum_k\min(j,k)m_k-m_j\bigr)\\
&=\sum_{j}x_jm_j\bigl(p_j+\frac{m_j+1}{2}+\sum_k\min(j,k)m_k\bigr)
\equiv 0 \mod \Z.
\end{split}
\end{equation*}
\end{proof}

Proposition \ref{pr:la} does not serve as a proof of 
Theorem \ref{th:dp}.
Nevertheless it has opened a route to create a conjectural
formula for the generic period in a large class of generalized 
periodic box-ball systems \cite{KT1,KT2}.
For $l=\infty$, Proposition \ref{pr:la} 
has been shown also in \cite{MIT} independently.

\subsection{Discussion}\label{subsec:ana}

We write the string center equation (\ref{eq:sce2}) 
in the matrix form:
\begin{equation}\label{eq:scemat}
A\vec{u} = \vec{c} + \vec{I} \in \Z^\gamma,
\end{equation}
where $\gamma = m_{j_1}+\cdots+m_{j_s}$ for 
$H = \{j_1, \ldots, j_s\}$ as in (\ref{eq:mh}).
The shifted array  
$\vec{I} = \vec{J}+\vec{\rho}$
will also be called the angle variable here. 
See around (\ref{eq:sce2}) for $\vec{\rho}$.
{}From (\ref{eq:odc}) and (\ref{eq:sce2}), 
the $\vec{I}$ belongs to the set
\begin{equation}\label{eq:jv}
({\mathcal I}_{m_{j_1}}\times 
\cdots \times {\mathcal I}_{m_{j_s}})/\Gamma,
\end{equation}
where ${\mathcal I}_n = (\Z^n-\Delta_n)/{\frak S}_n$ is 
the $n$ dimensional lattice without the diagonal 
points $\Delta_n = \{(z_1,\ldots, z_n) \in \Z^n
\mid z_\alpha = z_\beta \hbox{ for some }
1 \le \alpha \neq \beta \le  n \}$ 
identified under the permutations ${\frak S}_n$. 
$\Gamma= \bigoplus_{k \in H, 1 \le \beta \le m_k}
\Z \vec{A}_{k\beta}$ is the $\gamma$ dimensional lattice 
generated by the column vectors 
$\vec{A}_{k\beta}$ of the matrix $A=(A_{j\alpha, k\beta})$ (\ref{eq:a})
that characterizes the string center equation 
(\ref{eq:scemat})\footnote{
The column vectors $\vec{A}_{k\beta}$ are independent
because of $\det\!A >0$ under the condition $m \in \M$
as noted in the end of Section \ref{subsec:ba}.}.
The division by $\Gamma$ 
originates in the identification of $\vec{u}$ under 
$u^{(k)}_\beta \rightarrow u^{(k)}_\beta+1$ 
in (\ref{eq:scemat}).
The time evolution $T_l$  (\ref{eq:tlj}) on the angle 
variable $\vec{I}$ is expressed as the linear flow 
$T_l(\vec{I}) = \vec{I}+\vec{h}$ in terms of 
the $l$-dependent vector $\vec{h} \in \Z^\gamma$ 
defined in (\ref{eq:hvec}). 
According to Theorem \ref{th:ju}, it induces the 
time evolution of the Bethe roots 
$T_l(\vec{u}) = \vec{u} + A^{-1}\vec{h}$.
They are summarized in the following table:

\begin{equation*}
\begin{tabular}{c|c|c}
\hfil  &  Bethe roots $\vec{u}$ &  Angle variables $\vec{I}$  \\
\hline 
&&\vspace{-0.2cm}\\
$\mod$ & $\Z^\gamma$ & $\Gamma = A\Z^\gamma$ \\
&& \vspace{-0.3cm} \\
\hline
&& \vspace{-0.2cm}\\
$T_l$  & $A^{-1}\vec{h}$
& $\vec{h}$
\end{tabular}
\end{equation*}

\noindent
Here ${\rm mod}$ is the lattice under which the respective 
variables are to be identified save the 
permutations 
${\mathfrak S}_{m_{j_1}}\times 
\cdots \times {\mathfrak S}_{m_{j_s}}$.
Now it is transparent under what condition 
the time evolution $T^{\mathcal N}_l$ becomes trivial.
According to the above table, it is 
presented in the three equivalent forms: 
\begin{equation}\label{eq:ncond}
{\mathcal N}A^{-1}\vec{h} \in \Z^\gamma,\quad \;
{\mathcal N}\vec{h} \in \Gamma,\quad \;
{\mathcal N}F^{-1}\vec{h}' \in \Z^s,
\end{equation}
where the last one is a contracted version of the first 
due to (\ref{eq:ratio}). See (\ref{eq:hvec}) for $\vec{h}'$.
Similarly, the formula (\ref{eq:nn}) 
for the fundamental period is rephrased as
\begin{equation}\label{eq:fcond}
{\mathcal N}^\ast F^{-1}\vec{h}' \in G^{-1}\Z^s, \quad 
G = {\rm diag}(g_{j_1}, \ldots, g_{j_s}).
\end{equation}
In this way we arrive at the intrinsic meaning of the 
generic period (\ref{eq:lcm}) 
(rather than in the form (\ref{eq:lcm2})) and
the fundamental period (\ref{eq:nn}).
Namely, they are the smallest positive integers 
${\mathcal N}$ and ${\mathcal N}^\ast$ 
that make (\ref{eq:ncond}) and (\ref{eq:fcond}) valid.
Depending on the choice of $l$ in (\ref{eq:hvec}), 
the vectors $\vec{h}$ and $\vec{h}'$ (like Hamiltonian) 
encode various direction and 
speed of the straight motions in the set (\ref{eq:jv}) 
corresponding to the time evolution $T_l$.
The simplest among them is 
$\vec{h} = {}^t(1,1,\ldots,1) \in \Z^\gamma$ for $l=1$, 
for which $L\vec{h} \in \Gamma$ holds because of (\ref{eq:suml}).
This fact corresponds to the simple property 
$T^L_1(p)=p$ for any $p$. See (\ref{eq:t1}).
More generally, the generic and fundamental period 
under the combined time evolution $T = \prod_lT_l^{\beta_l}$
$(\beta_l \in \Z)$ 
can be obtained by replacing 
$\vec{h}=\vec{h}_l$ and $\vec{h}'=\vec{h'_l}$ (\ref{eq:hvec}) with 
$\sum_l \beta_l\vec{h}_l$ and $\sum_l \beta_l\vec{h'_l}$
in (\ref{eq:lcm}) and (\ref{eq:nn}), or equivalently,  
in (\ref{eq:ncond}) and (\ref{eq:fcond}).
The set (\ref{eq:jv}) is an analogue of 
the Jacobi variety on which 
the nonlinear dynamics on $\Pth$ 
looks as a straight motion as in the classical theory of 
quasi-periodic solutions to soliton equations \cite{DT,DMN}.
The coefficient $A$ in the string center equation 
plays the role of the period matrix.
Its size $\gamma$ is the total number of solitons, which is 
equal to the first energy $E_1$ (\ref{eq:em}).

Under any set of 
selected time evolutions $T_{l_1}, \ldots, T_{l_k}$,
one can describe the decomposition of $\Pth(m)$ 
into the disjoint union of orbits:   
\begin{equation*}
\Pth(m) = \bigsqcup \,
\{\hbox{orbit under } T_{l_1}, \ldots, T_{l_k}\}.
\end{equation*}
{}From (\ref{eq:jv}), each orbit here is in one to one 
correspondence with an element of
\begin{equation}\label{eq:jvh}
({\mathcal I}_{m_{j_1}}\times 
\cdots \times {\mathcal I}_{m_{j_s}})/
(\Gamma + \Z\vec{h}_{l_1} + \cdots + \Z\vec{h}_{l_k}).
\end{equation}
Equivalently, in terms of the original angle variable,
the orbits are labeled by  
$\J(m)/(T_{l_1}^{\Z}\cdots T_{l_k}^{\Z})$, where the 
division means the identification of the elements 
connected by the time evolution (\ref{eq:tlj}).

Under the single time evolution $T_l$, 
the number of orbits contained in 
$\Pth(m)$ is $\Omega(m)/{\mathcal N}$ 
if the action variable $m \in {\mathcal M}$ 
satisfies the generic condition (\ref{eq:nns}).
See (\ref{eq:cpju}) and (\ref{eq:omega}) for $\Omega(m)$.
Another simple situation is to include the 
whole family $T_1, T_2, \ldots$, which maximizes the 
orbits and minimizes their number.
{}From (\ref{eq:jbar})--(\ref{eq:j}) and (\ref{eq:tlj}) 
it is easy to see that 
$\J(m)/(T_{1}^{\Z}T_{2}^{\Z}\cdots)$ is factorized as
\begin{equation}
\widehat{\J}_{j_1} \times \widehat{\J}_{j_2} \times \cdots \times
\widehat{\J}_{j_s}.
\end{equation}
Here $\widehat{\J}_j= \J_j/\!\!\!\sim$ is obtained from 
$\J_j = \{(J_i)_{i \in \Z} \!\mid \!\!
J_i \in \Z, \, J_i \le J_{i+1}, \, J_{i+m_j}=J_i + p_j\,
\hbox{ for all } i \}$ in (\ref{eq:jcomp}) by 
the identification $(J_i) \sim (J'_i)$ defined by 
$J'_i=J_{i+1}$ for all $i$ or $J'_i = J_i + 1$ for all $i$.
The difference $d_i=J_{i+1}-J_i+1$ 
satisfies $d_i \in \Z_{\ge 1}$ and
$d_{i+1}+d_{i+2}+\cdots+d_{i+m_j}=p_j+m_j$, especially $d_i=d_{i+m_j}$
for all $i \in \Z$.
Thus the set $\widehat{\J}_j$ is in one to one correspondence 
with the arrangements of  
$p_j+m_j$ letters $1$ and $m_j$ letters $2$ as
\begin{equation*}\label{eq:conf}
\overbrace{1\ldots 1}^{d_1} \,2\, 
\overbrace{1\ldots 1}^{d_2}
\,2\,\overbrace{1\ldots 1}^{d_3} \cdots 
\overbrace{1\ldots 1}^{d_{m_j}}\,2
\end{equation*}
that are inequivalent under the periodic boundary condition.
No two letters $2$ are allowed to be adjacent because of $d_i \ge 1$.  
Such arrangements are the states of the 
periodic box-ball system containing amplitude 1 solitons only.
In this way we obtain
\begin{equation}\label{eq:cc}
\vert\J(m)/(T_{1}^{\Z}T_{2}^{\Z}\cdots)\vert
=C(p_{j_1}, m_{j_1}) C(p_{j_2}, m_{j_2}) \cdots C(p_{j_s}, m_{j_s}),
\end{equation}
where $C(p,m) \, (p \ge 0, m \ge 1)$ is the number of 
orbits in the size $p+2m$ periodic box-ball system 
containing $m$ solitons of amplitude 1 only.  
Namely, the count of orbits 
under the entire family $T_1, T_2, \ldots$ splits into 
individual blocks wherein all the solitons 
behave effectively as amplitude 1.
The quantity $C(p, m)$ is characterized 
as the number of 
orbits of the monomials $x^{\mu_1}_1\cdots x^{\mu_m}_m$
under the cyclic shift 
$x_i\!\rightarrow\!x_{i+1} \, (x_{i+m}\!=\!x_i)$  
in the complete symmetric function 
$h_p(x_1,\ldots, x_m)$. 

Finally we comment on the motion of the Bethe roots 
$\vec{u} \mapsto T_l(\vec{u})=\vec{u} + A^{-1}\vec{h}$.
The fusion transfer matrix corresponding to $T_l$ does not change 
the Bethe vectors up to an overall scalar nor 
the associated Bethe roots.
This does not contradict the motion of $\vec{u}$ considered here
since each path that we associate to 
$\vec{u}$ is a monomial in $({\mathbb C}^2)^{\otimes L}$,  
which is {\em not} a Bethe vector at $q=0$ in general.

\section{Summary}\label{sec:sum}
In this paper we solved the 
initial value problem in the periodic box-ball system 
by a unification of 
the combinatorial Bethe ans{\" a}tze at $q=1$ and $q=0$.
Section \ref{sec:2} gives the formulation of the periodic box-ball system 
in terms of crystal basis theory.
The commutativity, 
energy conservation (Theorem \ref{th:te}) and the 
invariance under the extended affine Weyl group 
(Propositions \ref{pr:w} and \ref{pr:ts}) in the periodic setting 
are firstly shown explicitly in this paper.
In Section \ref{sec:ism}, we introduced the 
action and angle variables in (\ref{eq:am}) and (\ref{eq:j}) 
and the linear time evolution (\ref{eq:tlj}) on the latter.
The direct/inverse scattering map is defined in (\ref{eq:pj}), 
which linearizes the dynamics as in Theorem \ref{th:main}.
In Section \ref{sec:ba}, 
our inverse scattering formalism are linked with the Bethe ansatz
by the key relation (\ref{eq:psi}).
The action-angle variables are in one to one correspondence 
with the off-diagonal solutions (\ref{eq:odc}) 
to the string center equation (\ref{eq:sce}) as 
summarized in Corollary \ref{cor:pju}.
It has led to the explicit formula (\ref{eq:cpju}), 
(\ref{eq:omega}) counting the 
states characterized either by soliton content (energy) or 
string content (configuration),
which we identified in Proposition \ref{pr:em} and 
Section \ref{subsec:up}.
As further applications, the generic (\ref{eq:lcm}) and fundamental 
(\ref{eq:nn}) period and the number of disjoint 
orbits (\ref{eq:cc}) under the commuting 
family of time evolutions are obtained 
(Theorems \ref{th:dp}, \ref{th:fp}). 
The Bethe eigenvalue is shown to be a root of unity 
related to the generic period (Proposition \ref{pr:la}). 
These results are derived and understood most naturally 
from the intrinsic picture on the dynamics as a straight motion in the 
`ultradiscrete Jacobi variety' (\ref{eq:jv}).
We expect that the essential features explored in this paper
persist in the generalized periodic box-ball systems \cite{KT1,KT2}.

\vspace{0.3cm}\noindent
{\bf Acknowledgments} \hspace{0.1cm}
This work is partially supported by Grand-in-Aid for Scientific 
Research JSPS No.15540363.

\appendix

\section{KKR bijection}\label{app:kkr}

Put $B = B_{l_1}\otimes \cdots \otimes B_{l_L}$.
The set ${\mathcal P}$ 
of our periodic box-ball system (\ref{eq:pdef})
corresponds to the $B$ with the choice 
$l_1=\cdots = l_L = 1$. Let 
\begin{equation}
\Pth_+ = \{ p \in B\mid {\tilde e}_1p=0 \}
\end{equation}
be the set of highest elements.
For the array of nonnegative integers 
$m=(m_j) = (m_1, m_2,\ldots)$, we put 
$H = \{j\in \Z_{\ge 1}\mid m_j>0\}$.
The data $m$ is called a {\em configuration} if $p_j \ge 0$
for all $j \in H$, where 
\begin{equation}\label{eq:vaa}
p_j = \sum_{i=1}^L\min(j,l_i) - 2\sum_{k\ge 1}\min(j,k)m_k
\end{equation}
is called the vacancy number. 
The set $H$ is necessarily finite and 
we parameterize it as $H = \{j_1< \cdots < j_s \}$.
The data $m$ is identified with the Young diagram
containing $m_j$ rows of length $j$, i.e., the
$m_j \times j$ rectangular block for each $j \in H$.
We let ${\mathcal M}$ denote the set of all the configurations 
$m$.
These definitions agree with the earlier ones for
$\Pth_+$ (\ref{eq:pdef}), ${\mathcal M}$
(\ref{eq:m}),  $p_j$ (\ref{eq:va}) and $H$ (\ref{eq:mh})
when $l_1=\cdots = l_L = 1$.
Define 
\begin{equation}\label{eq:rc}
\Rig(m) = \{(J^{(j)}_i)_{1 \le i \le m_j, j \in H}
\in \Z^{m_{j_1}} \times \cdots \times \Z^{m_{j_s}}
\mid 0 \le J^{(j)}_1 \le \cdots \le J^{(j)}_{m_j} \le p_j \}.
\end{equation}
When $l_1=\cdots = l_L = 1$ and the dependence 
on $L$ is important, we write $\Rig_L(m)$ in 
Appendices \ref{app:proof} and \ref{app:main}.
 
The purpose of this appendix is to describe the 
Kerov-Kirillov-Reshetikhin (KKR) 
bijection $\phi=\phi_B$ \cite{KKR,KR}:
\begin{equation}\label{eq:kkr}
\phi : \Pth_+ \longrightarrow
\sqcup_{m \in \M} \{(m, J)\mid J \in \Rig(m)\}.
\end{equation}

An element of $\Rig(m)$ is called {\em rigging}.
Each value $J^{(j)}_i$ is also said rigging.
Similarly $p_j - J^{(j)}_i$ is called co-rigging.
The combined data $(m,J)$ is called the rigged configuration.
It is customary to depict it as the Young diagram $m$
with each row attached with the rigging.
We often exhibit the vacancy number $p_j$ on the left of 
the block of width $j$.
The riggings are arranged so as to decrease weakly
downward within a block.
Note that these definitions 
depend on the choice of $B$ although we do not 
exhibit it explicitly.
The choice other than 
$l_1=\cdots =l_L =1$ is needed only in 
Proposition \ref{pr:hatten}.

\begin{example}\label{ex:rc}
For the highest paths $p_1, p_2$ and $p_3$ in 
Example \ref{ex:dp+},  one has 
$B = B_1^{\otimes 19}$ and 
\begin{equation*}
\begin{picture}(75,53)(-30,-50)
\put(0,0){\line(0,-1){50}}
\put(10,0){\line(0,-1){50}}
\put(20,0){\line(0,-1){30}}
\put(30,0){\line(0,-1){10}}

\put(-50,-25){$\phi(p_1)=$}
\put(0,0){\line(1,0){30}}\put(33,-9){1}
\put(0,-10){\line(1,0){30}}\put(-7,-9){1}\put(23,-19){1}
\put(0,-20){\line(1,0){20}}\put(-7,-24){3}\put(23,-29){0}
\put(0,-30){\line(1,0){20}}\put(13,-39){8}
\put(0,-40){\line(1,0){10}}\put(-7,-44){9}\put(13,-49){4}
\put(0,-50){\line(1,0){10}}
\end{picture}
\qquad 
\begin{picture}(75,53)(-30,-50)
\put(0,0){\line(0,-1){50}}
\put(10,0){\line(0,-1){50}}
\put(20,0){\line(0,-1){30}}
\put(30,0){\line(0,-1){10}}

\put(-50,-25){$\phi(p_2)=$}
\put(0,0){\line(1,0){30}}\put(33,-9){1}
\put(0,-10){\line(1,0){30}}\put(-7,-9){1}\put(23,-19){3}
\put(0,-20){\line(1,0){20}}\put(-7,-24){3}\put(23,-29){1}
\put(0,-30){\line(1,0){20}}\put(13,-39){6}
\put(0,-40){\line(1,0){10}}\put(-7,-44){9}\put(13,-49){2}
\put(0,-50){\line(1,0){10}}
\end{picture}
\qquad
\begin{picture}(75,53)(-30,-50)
\put(0,0){\line(0,-1){50}}
\put(10,0){\line(0,-1){50}}
\put(20,0){\line(0,-1){30}}
\put(30,0){\line(0,-1){10}}

\put(-50,-25){$\phi(p_3)=$}
\put(0,0){\line(1,0){30}}\put(33,-9){0}
\put(0,-10){\line(1,0){30}}\put(-7,-9){1}\put(23,-19){3}
\put(0,-20){\line(1,0){20}}\put(-7,-24){3}\put(23,-29){2}
\put(0,-30){\line(1,0){20}}\put(13,-39){8}
\put(0,-40){\line(1,0){10}}\put(-7,-44){9}\put(13,-49){3}
\put(0,-50){\line(1,0){10}}
\end{picture}
\end{equation*} 
\end{example}

The original KKR bijection \cite{KKR,KR} is the one between the 
rigged configurations and the Littlewood-Richardson tableaux.
The bijection (\ref{eq:kkr}) is obtained through 
a simple transformation of 
the Littlewood-Richardson tableaux and the highest paths \cite{NY}.
Here we illustrate 
$\phi$ and $\phi^{-1}$ casually along two examples rather than 
the systematic description which is already available in 
\cite{KKR,KR,KSS,Schi}.
Our convention here is opposite from 
\cite{KKR} in the role of the rigging and co-rigging,
and opposite from \cite{Schi} in the order of tensor product. 

Regard a rigged configuration 
as a multi set of the pairs 
$(\hbox{row length}, \hbox{attached rigging})$.
For example, the leftmost one in Example \ref{ex:rc}
is regarded as 
$\{(3,1), (2,1), (2,0), (1,8), (1,4)\}$.
Each element $(j, \alpha)$ 
of the rigged configuration is called a {\em string} 
with length $j$ and rigging $\alpha$.
A string is {\em singular} if the co-rigging is zero, 
namely $\alpha = p_j$, which is 
the maximum allowed value in (\ref{eq:rc}).
We first illustrate the map $\phi^{-1}$.

\begin{example}\label{ex:ph1}
We begin with the basic case $B=B^{\otimes L}_1$, where 
the vacancy number (\ref{eq:vaa}) reduces to (\ref{eq:va}).
Consider the configuration $m=(2,1)$ for example.
For $L=8$, there are 6 rigged configurations 
depicted in the leftmost column of the following:

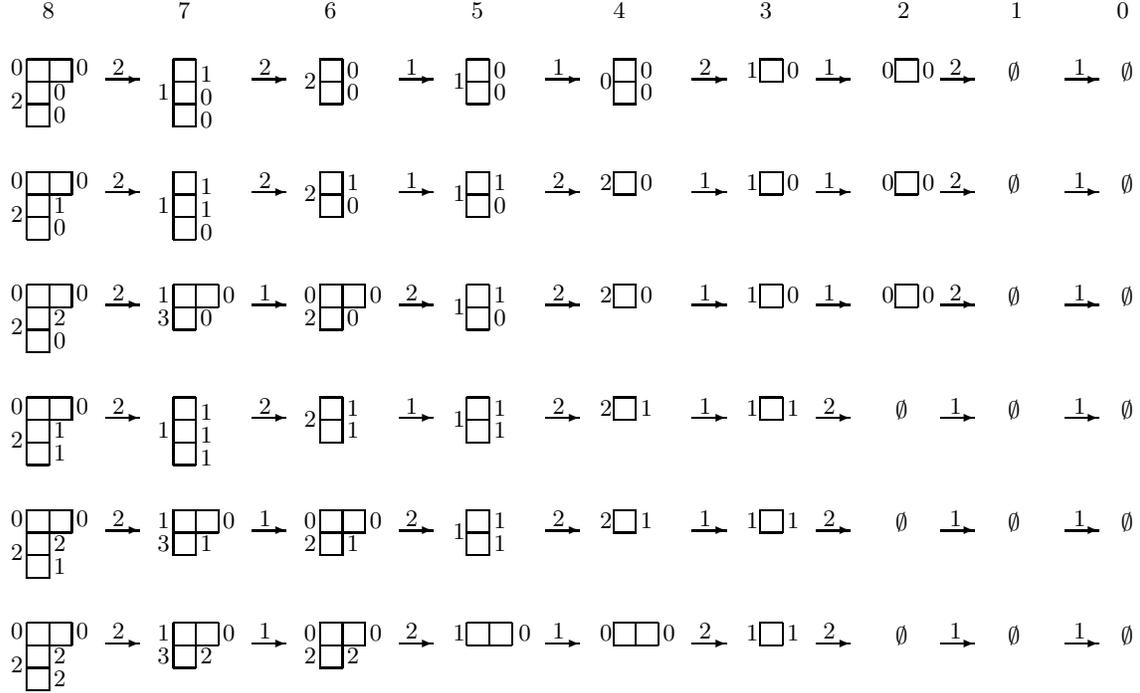
\setlength{\unitlength}{0.3mm}
\begin{picture}(450,343)(40,-11)
\multiput(0,0)(0,50){6}{
\put(0,0){\line(1,0){10}}
\put(0,10){\line(1,0){10}}
\put(0,20){\line(1,0){20}}
\put(0,30){\line(1,0){20}}
\put(0,0){\line(0,1){30}}
\put(10,0){\line(0,1){30}}
\put(20,20){\line(0,1){10}}
\put(-7,8){\small 2} \put(-7,23){\small 0}
\put(22,23){\small 0}
\multiput(35,21)(65,0){5}{\vector(1,0){15}}}
\multiput(0,0)(0,50){6}{
\multiput(350,21)(55,0){3}{\vector(1,0){15}}}

\multiput(0,0)(65,0){2}{
\multiput(0,0)(0,50){2}{
\put(65,10){\line(1,0){10}}
\put(65,20){\line(1,0){20}}
\put(65,30){\line(1,0){20}}
\put(65,10){\line(0,1){20}}
\put(75,10){\line(0,1){20}}
\put(85,20){\line(0,1){10}}}
\multiput(0,0)(0,150){2}{
\put(65,10){\line(1,0){10}}
\put(65,20){\line(1,0){20}}
\put(65,30){\line(1,0){20}}
\put(65,10){\line(0,1){20}}
\put(75,10){\line(0,1){20}}
\put(85,20){\line(0,1){10}}}}

\put(65,100){\line(1,0){10}}
\put(65,110){\line(1,0){10}}
\put(65,120){\line(1,0){10}}
\put(65,130){\line(1,0){10}}
\put(65,100){\line(0,1){30}}
\put(75,100){\line(0,1){30}}
\multiput(0,100)(0,50){2}{
\put(65,100){\line(1,0){10}}
\put(65,110){\line(1,0){10}}
\put(65,120){\line(1,0){10}}
\put(65,130){\line(1,0){10}}
\put(65,100){\line(0,1){30}}
\put(75,100){\line(0,1){30}}}

\multiput(0,0)(0,50){4}{
\put(130,110){\line(1,0){10}}
\put(130,120){\line(1,0){10}}
\put(130,130){\line(1,0){10}}
\put(130,110){\line(0,1){20}}
\put(140,110){\line(0,1){20}}}
\multiput(65,-50)(0,50){5}{
\put(130,110){\line(1,0){10}}
\put(130,120){\line(1,0){10}}
\put(130,130){\line(1,0){10}}
\put(130,110){\line(0,1){20}}
\put(140,110){\line(0,1){20}}}
\multiput(130,150)(0,00){1}{
\put(130,110){\line(1,0){10}}
\put(130,120){\line(1,0){10}}
\put(130,130){\line(1,0){10}}
\put(130,110){\line(0,1){20}}
\put(140,110){\line(0,1){20}}}

\multiput(0,0)(65,0){2}{
\put(195,20){\line(1,0){20}}
\put(195,30){\line(1,0){20}}
\put(195,20){\line(0,1){10}}
\put(205,20){\line(0,1){10}}
\put(215,20){\line(0,1){10}}}

\multiput(0,0)(0,50){4}{
\put(260,70){\line(1,0){10}}
\put(260,80){\line(1,0){10}}
\put(260,70){\line(0,1){10}}
\put(270,70){\line(0,1){10}}}
\multiput(65,-50)(0,50){6}{
\put(260,70){\line(1,0){10}}
\put(260,80){\line(1,0){10}}
\put(260,70){\line(0,1){10}}
\put(270,70){\line(0,1){10}}}
\multiput(125,100)(0,50){3}{
\put(260,70){\line(1,0){10}}
\put(260,80){\line(1,0){10}}
\put(260,70){\line(0,1){10}}
\put(270,70){\line(0,1){10}}}

\multiput(0,0)(0,50){3}{\put(385,21){\small $\emptyset$}}
\multiput(0,0)(50,0){2}{
\multiput(50,0)(0,50){6}{\put(385,21){\small $\emptyset$}}}

\put(7,298){\small 8}
\put(67,298){\small 7}
\put(132,298){\small 6}
\put(197,298){\small 5}
\put(260,298){\small 4}
\put(325,298){\small 3}
\put(386,298){\small 2}
\put(436,298){\small 1}
\put(483,298){\small 0}

\put(12,12){\small 2}\put(12,2){\small 2}

\put(12,62){\small 2}\put(12,52){\small 1}

\put(12,112){\small 1}\put(12,102){\small 1}

\put(12,162){\small 2}\put(12,152){\small 0}

\put(12,212){\small 1}\put(12,202){\small 0}

\put(12,262){\small 0}\put(12,252){\small 0}

\multiput(0,0)(0,50){6}{\put(38,23){\small 2}}

\put(58,22){\small 1}\put(58,12){\small 3}
\put(87,22){\small 0}\put(77,12){\small 2}

\put(58,72){\small 1}\put(58,62){\small 3}
\put(87,72){\small 0}\put(77,62){\small 1}

\put(58,112){\small 1}
\put(77,120){\small 1}\put(77,110){\small 1}\put(77,100){\small 1}

\put(58,172){\small 1}\put(58,162){\small 3}
\put(87,172){\small 0}\put(77,162){\small 0}

\put(58,212){\small 1}
\put(77,220){\small 1}\put(77,210){\small 1}\put(77,200){\small 0}

\put(58,262){\small 1}
\put(77,270){\small 1}\put(77,260){\small 0}\put(77,250){\small 0}

\put(103,23){\small 1}
\put(103,73){\small 1}
\put(103,123){\small 2}
\put(103,173){\small 1}
\put(103,223){\small 2}
\put(103,273){\small 2}

\put(123,22){\small 0}\put(123,12){\small 2}
\put(152,22){\small 0}\put(142,12){\small 2}

\put(123,72){\small 0}\put(123,62){\small 2}
\put(152,72){\small 0}\put(142,62){\small 1}

\put(123,117){\small 2}
\put(142,122){\small 1}\put(142,112){\small 1}

\put(123,172){\small 0}\put(123,162){\small 2}
\put(152,172){\small 0}\put(142,162){\small 0}

\put(123,217){\small 2}
\put(142,222){\small 1}\put(142,212){\small 0}

\put(123,267){\small 2}
\put(142,272){\small 0}\put(142,262){\small 0}

\put(168,23){\small 2}
\put(168,73){\small 2}
\put(168,123){\small 1}
\put(168,173){\small 2}
\put(168,223){\small 1}
\put(168,273){\small 1}

\put(189,22){\small 1}
\put(218,22){\small 0}

\put(189,67){\small 1}
\put(207,72){\small 1}\put(207,62){\small 1}

\put(189,117){\small 1}
\put(207,122){\small 1}\put(207,112){\small 1}

\put(189,167){\small 1}
\put(207,172){\small 1}\put(207,162){\small 0}

\put(189,217){\small 1}
\put(207,222){\small 1}\put(207,212){\small 0}

\put(189,267){\small 1}
\put(207,272){\small 0}\put(207,262){\small 0}

\put(233,23){\small 1}
\put(233,73){\small 2}
\put(233,123){\small 2}
\put(233,173){\small 2}
\put(233,223){\small 2}
\put(233,273){\small 1}

\put(254,22){\small 0}
\put(282,22){\small 0}

\put(254,72){\small 2}
\put(272,72){\small 1}

\put(254,122){\small 2}
\put(272,122){\small 1}

\put(254,172){\small 2}
\put(272,172){\small 0}

\put(254,222){\small 2}
\put(272,222){\small 0}

\put(254,267){\small 0}
\put(272,272){\small 0}\put(272,262){\small 0}

\put(298,23){\small 2}
\put(298,73){\small 1}
\put(298,123){\small 1}
\put(298,173){\small 1}
\put(298,223){\small 1}
\put(298,273){\small 2}

\put(319,22){\small 1}
\put(337,22){\small 1}

\put(319,72){\small 1}
\put(337,72){\small 1}

\put(319,122){\small 1}
\put(337,122){\small 1}

\put(319,172){\small 1}
\put(337,172){\small 0}

\put(319,222){\small 1}
\put(337,222){\small 0}

\put(319,272){\small 1}
\put(337,272){\small 0}

\put(353,23){\small 2}
\put(353,73){\small 2}
\put(353,123){\small 2}
\put(353,173){\small 1}
\put(353,223){\small 1}
\put(353,273){\small 1}

\put(379,172){\small 0}
\put(397,172){\small 0}

\put(379,222){\small 0}
\put(397,222){\small 0}

\put(379,272){\small 0}
\put(397,272){\small 0}

\put(409,23){\small 1}
\put(409,73){\small 1}
\put(409,123){\small 1}
\put(409,173){\small 2}
\put(409,223){\small 2}
\put(409,273){\small 2}

\multiput(464,23)(0,50){6}{\put(0,0){\small 1}}
\end{picture}

\noindent
The procedure to obtain the highest paths by 
applying $\phi^{-1}$ has been shown.
Reading the numbers on the arrows backward,
we find the image of those rigged configurations 
under $\phi^{-1}$ as follows: 
\begin{alignat*}{2}
&\quad 1\otimes 2\otimes 1\otimes 2\otimes 
1\otimes 1\otimes 2\otimes 2, & \\
&\quad 1\otimes 2\otimes 1\otimes 1\otimes 
2\otimes 1\otimes 2\otimes 2, & \\
&\quad 1\otimes 2\otimes 1\otimes 1\otimes 
2\otimes 2\otimes 1\otimes 2, &\\
&\quad 1\otimes 1\otimes 2\otimes 1\otimes 
2\otimes 1\otimes 2\otimes 2, & \\
&\quad 1\otimes 1\otimes 2\otimes 1\otimes 
2\otimes 2\otimes 1\otimes 2, & \\
&\quad 1\otimes 1\otimes 2\otimes 2\otimes 
1\otimes 2\otimes 1\otimes 2. &
\end{alignat*}
All these paths satisfy the highest condition (\ref{eq:lp}).
The KKR algorithm for obtaining $\phi^{-1}$ proceeds recursively 
as $\phi^{-1}_{B^{\otimes L}_1}({\rm rc})
= \phi^{-1}_{B^{\otimes L-1}_1}({\rm rc}')\otimes a$ with 
$a \in B_1 = \{1,2\}$.
This relation is depicted as 
${\rm rc} \overset{a}{\rightarrow} {\rm rc}'$ in the above.
We have $a=1 \in B_1$ and ${\rm rc}'={\rm rc}$ 
if the rigged configuration ${\rm rc}$ is free from singular strings. 
If there exist singular strings in ${\rm rc}$, 
we set $a=2 \in B_1$.
In that case the new rigged configuration ${\rm rc}'$ is 
obtained by replacing any one of the shortest singular string 
$(j,\alpha=p_j)$ by $(j-1, p'_{j-1})$, where
$p'_{j-1}$ is the vacancy number 
(\ref{eq:va}) with $L$ replaced by $L-1$
and $m_i$ replaced $m_i-\delta_{i,j}+\delta_{i,j-1}$.
Namely one removes 
a box from the shortest 
singular string and assings a new rigging to the shortened one
so that it again becomes singular in the new environment.
(If $j=1$, just eliminate it.)
Even when ${\rm rc}={\rm rc}'$, the vacancy number (\ref{eq:va})
is lowered by one by 
the change $B^{\otimes L}_1 \rightarrow B^{\otimes L-1}_1$.
So one must revise $p_j$ in each step 
and keep track of $L$, which we did on the top line. 
\end{example}
To apply our inverse scattering method, the description 
in Example \ref{ex:ph1} suffices for $\phi^{-1}$.
In order to cover the content of  
Proposition \ref{pr:hatten} which is used in 
Proposition \ref{pr:em} and 
Lemma \ref{lem:repeat}, one needs to go beyond 
$B=B^{\otimes L}_1$. We explain it along

\begin{example}\label{ex:ph2}
Take $B = B_2\otimes B_1 \otimes B_2 \otimes B_3 \otimes B_1$.
We show the procedure for obtaining a highest path in $B$.

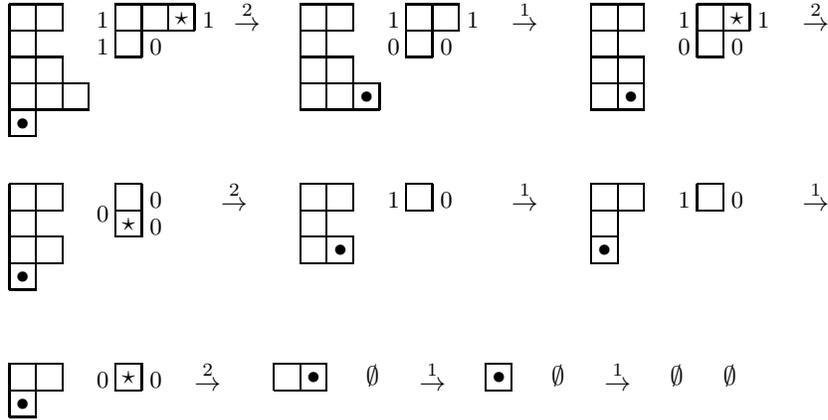
\begin{equation*}
\begin{split}
\setlength{\unitlength}{0.3mm}
&\begin{picture}(80,60)(0,-60)
\put(0,0){\line(1,0){20}}
\put(0,-10){\line(1,0){20}}
\put(0,-20){\line(1,0){20}}
\put(0,-30){\line(1,0){30}}
\put(0,-40){\line(1,0){30}}
\put(0,-50){\line(1,0){10}}
\put(0,0){\line(0,-1){50}}
\put(10,0){\line(0,-1){50}}
\put(20,0){\line(0,-1){10}}
\put(20,-20){\line(0,-1){20}}
\put(30,-30){\line(0,-1){10}}
\put(2.6,-47.5){$\bullet$}\put(62.6,-7.5){$\star$}
\put(33,-9){\small 1}\put(33,-19){\small 1}
\put(40,0){\line(1,0){30}}\put(73,-9){\small 1}
\put(40,-10){\line(1,0){30}}
\put(40,-20){\line(1,0){10}}\put(53,-19){\small 0}
\put(40,0){\line(0,-1){20}}
\put(50,0){\line(0,-1){20}}
\put(60,0){\line(0,-1){10}}
\put(70,0){\line(0,-1){10}}
\put(85,-10){$\overset{2}{\rightarrow}$}
\end{picture}
\qquad\quad
\begin{picture}(80,60)(0,-60)
\put(0,0){\line(1,0){20}}
\put(0,-10){\line(1,0){20}}
\put(0,-20){\line(1,0){20}}
\put(0,-30){\line(1,0){30}}
\put(0,-40){\line(1,0){30}}
\put(0,0){\line(0,-1){40}}
\put(10,0){\line(0,-1){40}}
\put(20,0){\line(0,-1){10}}
\put(20,-20){\line(0,-1){20}}
\put(30,-30){\line(0,-1){10}}
\put(22.6,-37.5){$\bullet$}
\put(33,-9){\small 1}\put(33,-19){\small 0}
\put(40,0){\line(1,0){20}}\put(63,-9){\small 1}
\put(40,-10){\line(1,0){20}}
\put(40,-20){\line(1,0){10}}\put(53,-19){\small 0}
\put(40,0){\line(0,-1){20}}
\put(50,0){\line(0,-1){20}}
\put(60,0){\line(0,-1){10}}
\put(80,-10){$\overset{1}{\rightarrow}$}
\end{picture}
\qquad\quad
\begin{picture}(80,60)(0,-60)
\put(0,0){\line(1,0){20}}
\put(0,-10){\line(1,0){20}}
\put(0,-20){\line(1,0){20}}
\put(0,-30){\line(1,0){20}}
\put(0,-40){\line(1,0){20}}
\put(0,0){\line(0,-1){40}}
\put(10,0){\line(0,-1){40}}
\put(20,0){\line(0,-1){10}}
\put(20,-20){\line(0,-1){20}}
\put(12.6,-37.5){$\bullet$}\put(52.6,-7.5){$\star$}
\put(33,-9){\small 1}\put(33,-19){\small 0}
\put(40,0){\line(1,0){20}}\put(63,-9){\small 1}
\put(40,-10){\line(1,0){20}}
\put(40,-20){\line(1,0){10}}\put(53,-19){\small 0}
\put(40,0){\line(0,-1){20}}
\put(50,0){\line(0,-1){20}}
\put(60,0){\line(0,-1){10}}
\put(80,-10){$\overset{2}{\rightarrow}$}
\end{picture}\\
&\begin{picture}(80,60)(0,-60)
\put(0,0){\line(1,0){20}}
\put(0,-10){\line(1,0){20}}
\put(0,-20){\line(1,0){20}}
\put(0,-30){\line(1,0){20}}
\put(0,-40){\line(1,0){10}}
\put(0,0){\line(0,-1){40}}
\put(10,0){\line(0,-1){40}}
\put(20,0){\line(0,-1){10}}
\put(20,-20){\line(0,-1){10}}
\put(2.6,-37.5){$\bullet$}\put(42.6,-17.5){$\star$}
\put(33,-14){\small 0}
\put(40,0){\line(1,0){10}}\put(53,-9){\small 0}
\put(40,-10){\line(1,0){10}}
\put(40,-20){\line(1,0){10}}\put(53,-19){\small 0}
\put(40,0){\line(0,-1){20}}
\put(50,0){\line(0,-1){20}}
\put(80,-10){$\overset{2}{\rightarrow}$}
\end{picture}
\qquad\quad
\begin{picture}(80,60)(0,-60)
\put(0,0){\line(1,0){20}}
\put(0,-10){\line(1,0){20}}
\put(0,-20){\line(1,0){20}}
\put(0,-30){\line(1,0){20}}
\put(0,0){\line(0,-1){30}}
\put(10,0){\line(0,-1){30}}
\put(20,0){\line(0,-1){10}}
\put(20,-20){\line(0,-1){10}}
\put(12.6,-27.5){$\bullet$}
\put(33,-9){\small 1}
\put(40,0){\line(1,0){10}}\put(53,-9){\small 0}
\put(40,-10){\line(1,0){10}}
\put(40,0){\line(0,-1){10}}
\put(50,0){\line(0,-1){10}}
\put(80,-10){$\overset{1}{\rightarrow}$}
\end{picture}
\qquad\quad
\begin{picture}(80,60)(0,-60)
\put(0,0){\line(1,0){20}}
\put(0,-10){\line(1,0){20}}
\put(0,-20){\line(1,0){10}}
\put(0,-30){\line(1,0){10}}
\put(0,0){\line(0,-1){30}}
\put(10,0){\line(0,-1){30}}
\put(20,0){\line(0,-1){10}}
\put(40,0){\line(1,0){10}}\put(53,-9){\small 0}
\put(40,-10){\line(1,0){10}}
\put(2.6,-27.5){$\bullet$}
\put(33,-9){\small 1}
\put(40,0){\line(0,-1){10}}
\put(50,0){\line(0,-1){10}}
\put(80,-10){$\overset{1}{\rightarrow}$}
\end{picture}\\
&\begin{picture}(80,30)(0,-30)
\put(0,0){\line(1,0){20}}
\put(0,-10){\line(1,0){20}}
\put(0,-20){\line(1,0){10}}
\put(0,0){\line(0,-1){20}}
\put(10,0){\line(0,-1){20}}
\put(20,0){\line(0,-1){10}}
\put(40,0){\line(1,0){10}}\put(53,-9){\small 0}
\put(40,-10){\line(1,0){10}}
\put(2.6,-17.5){$\bullet$}\put(42.6,-7.5){$\star$}
\put(33,-9){\small 0}
\put(40,0){\line(0,-1){10}}
\put(50,0){\line(0,-1){10}}
\put(70,-10){$\overset{2}{\rightarrow}$}
\end{picture}
\quad\quad
\begin{picture}(60,30)(0,-30)
\put(0,0){\line(1,0){20}}
\put(0,-10){\line(1,0){20}}
\put(0,0){\line(0,-1){10}}
\put(10,0){\line(0,-1){10}}
\put(20,0){\line(0,-1){10}}
\put(12.6,-7.5){$\bullet$}
\put(35,-8){$\emptyset$}
\put(55,-10){$\overset{1}{\rightarrow}$}
\end{picture}
\quad\quad
\begin{picture}(60,30)(0,-30)
\put(0,0){\line(1,0){10}}
\put(0,-10){\line(1,0){10}}
\put(0,0){\line(0,-1){10}}
\put(10,0){\line(0,-1){10}}
\put(2.6,-7.5){$\bullet$}
\put(25,-8){$\emptyset$}
\put(45,-10){$\overset{1}{\rightarrow}$}
\put(70,-8){$\emptyset$}
\put(90,-8){$\emptyset$}
\end{picture}
\end{split}
\end{equation*}
This time we have drawn a pair of diagrams in each step.
The right ones are the rigged configurations.
The left ones keep track of the `shape' of the paths.
For example in the top left diagram, the list of its row 
lengths $2, 1, 2, 3, 1$ encodes the indices in 
$B = B_2\otimes B_1 \otimes B_2 \otimes B_3 \otimes B_1$
where we start from.
They are removed one by one from the bottom right 
as marked with $\bullet$. The process 
corresponds to the reduction of the path shape as 
\begin{align*}
&B_2\otimes B_1 \otimes B_2 \otimes B_3 \otimes B_1
\rightarrow 
B_2\otimes B_1 \otimes B_2 \otimes B_3 \rightarrow 
B_2\otimes B_1 \otimes B_2 \otimes B_2 \rightarrow \\
&
B_2\otimes B_1 \otimes B_2 \otimes B_1 \rightarrow 
B_2\otimes B_1 \otimes B_2 \rightarrow 
B_2\otimes B_1 \otimes B_1 \rightarrow \\
&
B_2 \otimes B_1 \rightarrow B_2 \rightarrow 
B_1 \rightarrow \emptyset.
\end{align*}
In each step 
the vacancy numbers are revised by 
adopting the intermediate path shape as $\{l_i\}$ in (\ref{eq:vaa}).
When a box (marked with $\bullet$)  
is removed from the $k$ th column, in accordance with 
$(\cdots) \otimes B_k \rightarrow (\cdots)\otimes B_{k-1}$,
the rigged configuration is left unchanged 
if there is no singular string with length $\ge k$.
In this case we proceed to the next 
step by $\overset{1}{\rightarrow}$.
If there exist singular strings with length $\ge k$,
we have $\overset{2}{\rightarrow}$ and 
remove a box (as marked with $\star$)  
from any one of the shortest such string.
The shortened new string shall be attached with the rigging 
so that it becomes singular in the new environment.
The algorithm ends up with the pair of $\emptyset$.
Reading the letters on the arrows backward,
we get the sequence 
$11\vert 2 \vert 11 \vert 221 \vert 2$,
where the symbol $\vert$ separates those letters 
coming from different rows in the shape diagrams.
Reversing the letters within $\vert \cdot \vert$,
we find the image under $\phi^{-1}_B$ as
$11 \otimes 2 \otimes 11 \otimes 122 \otimes 2 \in B$,
which is highest.
\end{example}
The algorithm illustrated in Example \ref{ex:ph2}
reduces to the simpler one in Example \ref{ex:ph1}
when $l_1=\cdots = l_L=1$.

We have seen that the KKR algorithm 
for $\phi^{-1}_B$ is a removal process of 
a rigged configuration 
creating a highest path from its rightmost component.
Naturally, the map $\phi_B$ is an addition process
building a rigged configuration by using  
the information of a highest path 
from its leftmost component. 
The rule of addition is easily inferred 
by looking at Example \ref{ex:ph2}
backward.
One regards the given path 
$11 \otimes 2 \otimes 11 \otimes 122 \otimes 2 \in B$ 
as the word $11\vert 2 \vert 11 \vert 221 \vert 2$
and adds a box ($\bullet$) carrying these letters 
one by one from the left to form the prescribed shape diagram
having the row lengths $2,1,2,3,2$ from the top.
If the letter is $1$, one does nothing on the rigged configuration. 
If the letter is $2$ and the change of the path shape is
$(\cdots)\otimes B_{k} \rightarrow (\cdots)\otimes B_{k+1}$,
one adds a box ($\star$) to any one of the {\em longest} singular 
string $(j,\alpha=p_j)$ among those $j \ge k$.
The new string of length $j\!+\!1$ should 
be assigned with the maximal rigging 
so as to become singular in the new environment.
If there is no such strings, one creates a singular string of length 1.

For the highest path $p \in \Pth=B^{\otimes L}_1$ 
such that $\phi_\Pth(p) = (m,J)$,
it is easy to see 
\begin{equation*}
\phi_{\Pth\otimes B^{\otimes n}_1}
(p\otimes \overbrace{1\otimes \cdots \otimes 1}^n) = (m,J).
\end{equation*}

\begin{proposition}\label{pr:hatten}
Let $p = b_1\otimes  \cdots \otimes b_L 
\in \Pth = B^{\otimes L}_1$ 
be a highest path such that $\phi(p) = (m,J)$.
For $n$ sufficiently large, define 
$\xi \in \Pth \otimes B^{\otimes n}_1$ by
$u_l \otimes (p\otimes 1^{\otimes n}) \simeq 
\xi \otimes u_l$.
Then $\phi_{\Pth\otimes B^{\otimes n}_1}(\xi) = 
(m,I)$, where the rigging $I=(I^{(j)}_i)$ is given by
$I^{(j)}_i = J^{(j)}_i +\min(l,j)$.
\end{proposition}
In $\xi$, sufficiently many components on the right 
are also $1\in B_1$.

\begin{proof}
We first show that 
$\phi_{B_l\otimes \Pth\otimes B^{\otimes n}_1}
(u_l \otimes p\otimes 1^{\otimes n})=(m,I)$.
In fact, the change from $\phi_{\Pth \otimes B^{\otimes n}_1}$ 
to $\phi_{B_l \otimes \Pth \otimes B^{\otimes n}_1}$ 
increases the vacancy number $p_j$ (\ref{eq:vaa}) to $p_j+\min(l,j)$.
Exactly the same increment has occurred from the rigging $J$ to $I$.
Therefore the co-rigging of $(m,J)$ for 
$\Pth \otimes B^{\otimes n}_1$ and the co-rigging of 
$(m,I)$ for $B_l \otimes \Pth \otimes B^{\otimes n}_1$ coincide,
leading to $\phi_{B_l\otimes \Pth\otimes B^{\otimes n}_1}
(u_l \otimes p\otimes 1^{\otimes n})=(m,I)$. 
On the other hand we have 
$\phi_{B_l\otimes \Pth\otimes B^{\otimes n}_1}
(u_l \otimes p\otimes 1^{\otimes n}) = 
\phi_{\Pth\otimes B^{\otimes n}_1\otimes B_l}
(\xi \otimes u_l)$ by Lemma 8.5 in \cite{KSS}.
In view of the KKR algorithm, the last
one is equal to $\phi_{\Pth\otimes B^{\otimes n}_1}(\xi)$. 
\end{proof}

The image 
\begin{equation*}
\phi^{-1}_B((m,J))=(l_1-x_1,x_1)\otimes \cdots 
\otimes (l_L-x_L,x_L) \in 
B=B_{l_1}\otimes \cdots \otimes B_{l_L}
\end{equation*}
can also be described by the piecewise linear formula:
\begin{equation}\label{eq:pwf}
\begin{split}
x_n &= \tau_1(n)-\tau_1(n-1)-\tau_0(n)+\tau_0(n-1),\\
\tau_i(n)&= \max_{\nu}
\{\sum_k\bigl(\,\sum_{g=1}^n\min(l_g,\nu_k)+(i-2)\nu_k-I_k\bigr)
-2\sum_{j<k}\min(\nu_j,\nu_k)\},
\end{split}
\end{equation}
where the maximum is taken over the subset 
$\nu = \{(\nu_1, I_1), (\nu_2, I_2), \ldots \}$ 
of the rigged configuration $(m,J)$
regarded as the multi set of 
$(\hbox{row length}, \;\hbox{attached rigging})$.

\section{Proof of Proposition \ref{pr:em}}\label{app:em}

Since $E_l(\omega(p))=E_l(p)$ by Proposition \ref{pr:w}, we 
may assume that $\wt(p)\ge 0$.
Such a path can be 
expressed as $p = T^d_1(p_+)$ for some $d \in \Z$ and $p_+ \in \Pth$
due to Lemma \ref{lem:dp+}.
Then $E_l(p) = E_l(T^d_1(p_+))=E_l(p_+)$, where the last 
equality is due to Theorem \ref{th:te}.
Therefore it suffices to show (\ref{eq:em}) 
for $p \in \Pth_+$, which we shall assume in the sequel.

For any element of the form 
$b \in B=B_{l_1} \otimes \cdots \otimes B_{l_k}$  
we define the quantity $D_l(b) \in \Z_{\ge 0}$ by
\begin{equation*}
\zeta^0 u_l \otimes b \simeq b^\ast \otimes 
\zeta^{D_l(b)}w\quad 
(b^\ast \in B,\; w \in B_l)
\end{equation*}
under the isomorphism 
$\Aff(B_l)\otimes \Aff(B) \simeq \Aff(B)\otimes \Aff(B_l)$.
(We have omitted the spectral parameters for $b,b^\ast$.)
When $p$ is highest, its energy $E_l(p)$ emerges not only in (\ref{eq:tle})
but already in the relation (\ref{eq:vt}) that produces $v_l$.
Namely, 
\begin{lemma}\label{lem:evl}
For any highest path $p \in \Pth_+$, one has 
$D_l(p)=E_l(p)$.
\end{lemma}

The following proof is direct but not intrinsic.
\begin{proof}
Fixing a highest path $p = b_1 \otimes \cdots \otimes b_L$,
we regard $D_l(p)$ and $y_2$ in Figure \ref{fig:xp} 
as functions of $x_2$ there. 
The quantity $D_l(p)$ is the number of  
vertices of the bottom right type in Figure \ref{fig:cr}.
Call them scoring vertices.
We keep track of the scoring vertices 
appearing in Figure \ref{fig:xp} as
$x_2$ is increased from $0$ to $c=y_2(l)\le l$.
$D_l(p)=E_l(p)$ is shown if the scoring vertices for $x_2=0$ 
remain scoring and no new scoring vertices are 
created during the increment of $x_2$.
To see this, recall the property 
$y_2(x_2+1)=y_2(x_2)$ or $y_2(x_2+1)=y_2(x_2)+1$ 
as noted in the proof of the Proposition \ref{pr:tl}.
Since $b_1, \ldots, b_L$  are fixed, a little inspection of 
Figure \ref{fig:cr} tells that
the increment of $x_2$ never transforms the other type of 
vertices into scoring ones.
On the other hand a scoring vertex may change into the 
top right type in Figure \ref{fig:cr} if $a=l-1$.
Suppose this firstly happened as $x_2=r$ is increased to 
$x_2 = r+1$ for some $r<c \,(\le l)$. 
Consider the leftmost such vertex 
that has ceased to be scoring at $x_2=r+1$.
Such a situation is realized only if 
all the vertices (at $x_2=r+1$) 
on its left are bottom two types in Figure \ref{fig:cr}.
Let $\alpha$ and $\beta$ be the number of 
the bottom left types and the bottom right (scoring) types in them.
Then one has $r+1-\alpha+\beta = l$.
On the other hand from the highest path condition 
(\ref{eq:lp}), one also has $\alpha \ge \beta+1$, where the last $+1$
is the contribution of the very vertex that has ceased to be 
scoring.
Thus we obtain $r=l+\alpha-\beta-1\ge l$, which is a contradiction.
\end{proof}

\begin{proof}[Proof of Proposition \ref{pr:em}]
Let $p \in \Pth_+$ be a highest path 
such that $\phi(p)=(m,J)$ hence $\mu(p)=m$.
Suppose that 
$u^{\otimes A}_a \otimes p \otimes 
\overbrace{1 \otimes \cdots \otimes 1}^n
\simeq p^\ast \otimes 1 \otimes \cdots \otimes 1 \otimes 
u^{\otimes A}_a$ is valid. 
Here we take $1 \ll A \ll n$. 
In general $p^\ast$ is a highest path longer 
than $L$ (but much shorter than $n$). 
By Lemma \ref{lem:evl} and 
$H(u_l \otimes u_k)=0$, we have 
$E_l(p) = D_l(p) 
= D_l(u^{\otimes A}_a \otimes p \otimes 1 \otimes \cdots \otimes 1)
= D_l(p^\ast \otimes 1 \otimes \cdots \otimes 1 \otimes 
u^{\otimes A}_a)
= D_l(p^\ast \otimes 1 \otimes \cdots \otimes 1)$.
By Proposition \ref{pr:hatten}, we know that
$\phi(p^\ast \otimes 1 \otimes \cdots \otimes 1) = 
(m,I)$, where $I=(I^{(j)}_i)$ reads 
$I^{(j)}_i = J^{(j)}_i + A\min(a,j)$.
Therefore by taking $A$ and $a$ sufficiently large, 
one can achieve the situation
$1 \ll  \cdots \ll I^{(j)}_1 \le \cdots \le I^{(j)}_{m_j}  
\ll I^{(j+1)}_1 \le \cdots \le I^{(j+1)}_{m_{j+1}} \ll \cdots$.
For such a rigged configuration $(m,I)$, the KKR algorithm 
produces the path $p^\ast \otimes 1 \otimes \cdots \otimes 1
= \phi^{-1}((m,I))$ of the form:
\begin{equation}\label{eq:ss}
\ldots\; 
2^{\otimes j} \ldots \;2^{\otimes j} \, \ldots \;2^{\otimes j} \,
\overbrace{\ldots\ldots\ldots}^{L_j}
\;2^{\otimes j\!+\!1} \,\ldots \;2^{\otimes j\!+\!1} \,\ldots
\;2^{\otimes j\!+\!1} \,\ldots
\;2^{\otimes j\!+\!1}\ldots
\end{equation}
Here $\ldots$ means an array of $1 \in B_1$ and 
$2^{\otimes j} \in B^{\otimes j}_1$.
The number of appearance of $2^{\otimes j}$ is $m_j$.
One can satisfy $L_j\gg l$ for all $L_j$ by taking $A$ large.
Moreover from the KKR algorithm, 
there are at least $1^{\otimes j}$ between any two $2^{\otimes j}$.
For such a path $p^\ast\otimes 1 \otimes  \cdots \otimes 1$,
it is straightforward to check 
$D_l(p^\ast \otimes 1 \otimes \cdots \otimes 1)
=\sum_{j\ge 1}\min(l,j)m_j$. 
\end{proof}

\begin{remark}\label{rem:sol}
In the pattern like (\ref{eq:ss}), 
$2^{\otimes j}$ is a soliton with length (amplitude) $j$.
In this context the data $m=(m_j)$ tells 
that there are $m_j$ solitons
with length $j$.
\end{remark}

\section{Proof of Proposition \ref{pr:migi}}\label{app:proof}

\begin{lemma}\label{lem:waru}
Let $q \in B^{\otimes d}_1$ be a highest path of length $d$ and 
$r \in B^{\otimes L-d}_1$ be a highest path of length $L-d$.
Suppose that their rigged configurations are 
$\phi(q) = (l,I)$ and $\phi(r)= (n,K)$.
Then the rigged configuration of the 
highest path $q\otimes r \in B^{\otimes L}_1$ is given by
$\phi(q \otimes r) = (l\cup n, I \cup K')$, where 
$K'=(K'^{(j)}_i)$ is given by
\begin{equation}\label{eq:kp}
K'^{(j)}_i = K^{(j)}_i -2\sum_k\min(j,k)l_k+d,
\end{equation}
where $(l\cup n, I \cup K')$ means the union regarding 
$(l,I)$ and $(n,K')$ 
as multi-sets of rows assigned with rigging.  
\end{lemma}
\begin{proof}
The new vacancy number $p'_j$ for the width $j$ block 
in $(l\cup n, I \cup K')$ reads 
$p'_j = L-2\sum_k\min(j,k)l_k - 2 \sum_k\min(j,k)n_k$.
When applying the KKR map $\phi^{-1}$
to $(l\cup n, I \cup K')$, the co-rigging of the $i$-th row of 
the width $j$ block in $(n,K')$ is 
\begin{equation}
p'_j - K'^{(j)}_i = (L-d-2\sum_k\min(j,k)n_k)-K^{(j)}_i.
\end{equation}
The right hand side is equal to the co-rigging of the same 
row in the rigged configuration $(n,K)$.
Meanwhile the co-rigging of the $(l,I)$ part is not less than 
that in the original $q$ plus the vacancy number of $r$.
Therefore the algorithm of the map $\phi^{-1}$ proceeds as
$\phi^{-1}((l\cup n, I \cup K')) = 
\phi^{-1}((l,I))\otimes r = q \otimes r$.
\end{proof}

\begin{example}\label{ex:waru}
Take $q = 1112122, r=111221221122$, hence $d=7, L=19$.
The highest path $p_2$ in Examples \ref{ex:dp+} and \ref{ex:rc} 
is expressed as $p_2 = q \otimes r$. 
\begin{equation}
\begin{picture}(60,45)(0,-30)
\put(0,0){\line(1,0){20}}\put(23,-9){1}
\put(0,-10){\line(1,0){20}}\put(13,-19){2}
\put(0,-20){\line(1,0){10}}

\put(-35,-14){$\phi(q)=$}
\put(0,0){\line(0,-1){20}}
\put(10,0){\line(0,-1){20}}
\put(20,0){\line(0,-1){10}}
\end{picture}
\qquad 
\begin{picture}(60,45)(0,-30)
\put(0,0){\line(1,0){30}}\put(33,-9){0}
\put(0,-10){\line(1,0){30}}\put(23,-19){2}
\put(0,-20){\line(1,0){20}}\put(13,-29){3}
\put(0,-30){\line(1,0){10}}

\put(-35,-14){$\phi(r)=$}
\put(0,0){\line(0,-1){30}}
\put(10,0){\line(0,-1){30}}
\put(20,0){\line(0,-1){20}}
\put(30,0){\line(0,-1){10}}
\end{picture}
\end{equation}
The rigged configuration for $\phi(p_2)$ in Example \ref{ex:rc}
contains this $\phi(q)$ indeed.  The rest of it corresponding 
to $(n,K')$ is related to $\phi(r)=(n,K)$ as
\begin{equation*}
\begin{picture}(60,40)(0,-30)
\put(0,0){\line(1,0){30}}\put(33,-9){1}
\put(0,-10){\line(1,0){30}}\put(23,-19){3}
\put(0,-20){\line(1,0){20}}\put(13,-29){6}
\put(0,-30){\line(1,0){10}}

\put(0,0){\line(0,-1){30}}
\put(10,0){\line(0,-1){30}}
\put(20,0){\line(0,-1){20}}
\put(30,0){\line(0,-1){10}}
\end{picture}
\begin{picture}(60,40)(0,-30)
\put(0,0){\line(1,0){30}}\put(33,-9){$0-2\times 3+7$}
\put(0,-10){\line(1,0){30}}\put(23,-19){$2-2\times 3+7$}
\put(0,-20){\line(1,0){20}}\put(13,-29){$3-2\times 2 + 7$}
\put(0,-30){\line(1,0){10}}
\put(-15,-14){$=$}
\put(0,0){\line(0,-1){30}}
\put(10,0){\line(0,-1){30}}
\put(20,0){\line(0,-1){20}}
\put(30,0){\line(0,-1){10}}
\end{picture}
\end{equation*}
in agreement with (\ref{eq:kp}).
\end{example}

\begin{lemma}\label{lem:rotate}
Let $q \in B^{\otimes d}_1$ be a highest path of length $d$ and 
$r \in B^{\otimes L-d}_1$ be a highest path of length $L-d$.
Set $\phi(r\otimes q) = (m,J)$ and 
$\phi(q\otimes r)= (m',J')$.
Then $m=m'$ and 
$\iota(J')\simeq \iota(J)+d \in \Jb(m)$ are valid.
\end{lemma}
See (\ref{eq:j}) and (\ref{eq:iota}) for the definitions of 
$\simeq$ and $\iota$.
\begin{proof}
We use the same notation as in Lemma \ref{lem:waru}.
Thus $m'=l \cup n$ and $J'=I\cup K'$. 
Applying Lemma \ref{lem:waru}
with $(l,I)$ and $(n,K)$ interchanged, we find 
$m = l\cup n = m'$ and $J=I' \cup K$ with 
$I'^{(j)}_i = I^{(j)}_i-2\sum_k\min(j,k)n_k+L-d$.
In terms of the vacancy number $p'_j$ for $m$ in the proof of
Lemma \ref{lem:waru}, the relation between the rigging 
$J$ and $J'$ is summarized as
\begin{align}
I'^{(j)}_i&=I^{(j)}_i+2\sum_k\min(j,k)l_k+p'_j-d,\label{eq:ipi}\\
K^{(j)}_i&=K'^{(j)}_i+2\sum_k\min(j,k)l_k-d.\label{eq:kkp}
\end{align}
See the picture after (\ref{eq:ijk}).
Now switch to the extended sequences $\iota(J)$ and $\iota(J')$ 
by (\ref{eq:iota}).
Apart from $-d$, this is exactly the effect of the slide 
$\prod_k\sigma_k^{l_k}$ (\ref{eq:slide}) on the width $j$ block
of $\iota(J')$.
Therefore we conclude $\iota(J) = (\prod_k\sigma_k^{l_k})\iota(J')-d$.
\end{proof}

\begin{proof}[Proof of $\Rightarrow$ in Proposition \ref{pr:migi}]
Without loss of generality we may assume 
$d'=0$ and $0 \le d < L$.
Since $p_+$ and $p'_+$ are highest paths, they must have the form 
$p_+ = r \otimes q$ and $p'_+ = q \otimes r$, where the both paths 
$q \in B^{\otimes d}_1$ and $r \in B^{\otimes L-d}_1$ 
are highest ones.
Then the assertion follows from Lemma \ref{lem:rotate}.
\end{proof}

\vspace{0.3cm}

To show $\Leftarrow$ in Proposition \ref{pr:migi},
we explore the full implication of the right hand side of 
(\ref{eq:key}).
Let $H=\{j_1<\cdots < j_s\}$ be the set of lengths of rows of $m$
as in (\ref{eq:mh}).

\begin{lemma}\label{lem:kagiru}
Suppose 
$\iota(J) + e \simeq \iota(J')+e' \in \Jb(m)$ for some 
$m=(m_j) \in \M$, $J, J' \in \Rig(m)$ and $e, e' \in \Z$.
Then there exist $0 \le l_j \le m_j$ $(j \in H)$ and $0 \le d < L$ 
such that the relation 
$\sigma_{j_1}^{l_{j_1}}\cdots \sigma_{j_s}^{l_{j_s}}\iota(J')
=\iota(J)+d$
holds.  
\end{lemma}

\begin{proof}
{}From the assumption, there exist $l_j$'s and $d$ such that 
$\sigma_{j_1}^{l_{j_1}}\cdots \sigma_{j_s}^{l_{j_s}}\iota(J')
=\iota(J)+d$.
Our task is to show that they can always be chosen within 
the range $0 \le l_j \le m_j$ $(j \in H)$ and $0 \le d < L$.
{}From the last comment in Section \ref{subsec:aa}, 
to achieve the relation
$\sigma_{j_1}^{l_{j_1}} \cdots \sigma_{j_s}^{l_{j_s}}\iota(J')
=\iota(J)+d$ by the successive transformations 
$\iota(J') \rightarrow \sigma_{j_s}^{l_{j_s}}\iota(J') 
\rightarrow \sigma_{j_{s-1}}^{l_{j_{s-1}}}\sigma_{j_s}^{l_{j_s}}\iota(J')
\rightarrow \cdots 
\rightarrow \sigma_{j_1}^{l_{j_1}} \cdots \sigma_{j_s}^{l_{j_s}}\iota(J')$,
one must apply $\sigma^{l_{j_k}}_{j_k}$ so that  
the blocks $\J_{j_k}, \J_{j_{k+1}}, \ldots, \J_{j_s}$ 
in $\sigma^{l_{j_k}}_{j_k}\cdots \sigma^{l_{j_s}}_{j_s}\iota(J')$
already coincides with that in $\iota(J)$ up to an overall additive 
constant.
Having this in mind we proceed to showing 
$0 \le l_j \le m_j$.

Certainly, $l_{j_s}$ can be taken as 
$0 \le l_{j_s} < m_{j_s}$ by using (\ref{eq:sl}).
Since the slides are invertible, one may also assume that 
the first non-zero number in the sequence 
$l_{j_s}, l_{j_{s-1}},\ldots$ is positive, namely,
$l_{j_s}=l_{j_{s-1}}=\cdots =l_{j_{t+1}}=0$ and $l_{j_t}>0$.
If there is no such $t$, it follows that $\iota(J') = \iota(J)+d$ and
we are done because $d \ge 0$ may be assumed without 
loss of generality and then $d < L$ is obvious 
from $L>p_j$ (\ref{eq:va}) and (\ref{eq:rc}).
Henceforth we assume that $l_{j_t}>0$ exists 
for some $1 \le t \le s$.
We first claim that $l_{j_t} \le m_{j_t}$.
In fact, the case  $t=s$ is within our assumption. 
If $l_{j_t}> m_{j_t}$ for $t<s$, 
the upper blocks $\J_\beta \, (\beta > j_t)$ 
acquire the uniform shift $2j_tl_{j_t}$ under $\sigma^{l_{j_t}}_{j_t}$. 
On the other hand the rigging 
$J'^{(j_t)}_{m_{j_t}}$ gets shifted at least by $2p_{j_t}+2j_tl_{j_t}$.
Therefore to adjust $\sigma^{l_{j_t}}_{j_t}\iota(J')$ to some $\iota(J)$, 
one must extract at least $p_{j_t}+2j_tl_{j_t}$ 
from the new rigging on the upper blocks $\J_\beta$ 
as an overall constant.
But this fails since the original rigging there is not greater than
$p_\beta (< p_{j_t})$ hence the result of the above 
extraction is not greater than 
$(p_\beta+2j_tl_{j_t}) - (p_{j_t}+2j_tl_{j_t}) < 0$.
Thus we have verified $l_{j_t} \le m_{j_t}$.

Next we prove 
$0 \le l_{j_a} \le m_{j_a}$ $(1 \le a \le t)$ by induction on $a$ 
assuming that 
$0 \le l_{j_b} \le m_{j_b}$ for all $a+1 \le b \le t$.
Setting $K = (K^{(j)}_i) = 
\sigma^{l_{j_{a}}}_{j_{a}}\cdots \sigma^{l_{j_t}}_{j_t} \iota(J')$, 
we have
$K^{(j_b)}_i = J'^{(j_b)}_{i+l_{j_b}}+2\sum_{c=a}^t\min(j_b,j_c)l_{j_c}$. 
Note that 
\begin{align}
K^{(j_t)}_{m_{j_t}} &\ge  p_{j_t} + 2\sum_{b=a}^t\min(j_t,j_b)l_{j_b},
\label{eq:pjt}\\
K^{(j_s)}_1 & \le p_{j_s}+2\sum_{b=a}^t\min(j_s,j_b)l_{j_b},
\label{eq:pjs}
\end{align}
where the former follows from $l_{j_t} >0$ and 
the latter does from $l_{j_s}<m_{j_s}$ and $J'^{(j_s)}_{m_{j_s}} \le p_{j_s}$.
As explained above, the blocks $\J_{\beta}$ $(\beta \ge j_a)$ in 
$K-{\tilde d}$
should already coincide with those in $\iota(J)$ for some 
${\tilde d}$.
To show $l_{j_a} \ge 0$, suppose $l_{j_a}<0$ on the contrary.
Then we have $K^{(j_a)}_1\le J'^{(j_a)}_{m_{j_a}}-p_{j_a}
+ 2\sum_{b=a}^t\min(j_a,j_b)l_{j_b} \le 
2\sum_{b=a}^t\min(j_a,j_b)l_{j_b}$.
{}From this and (\ref{eq:pjt}),  
$K-{\tilde d}$ can coincide with some $\iota(J)$ only if 
\begin{equation*}
2\sum_{b=a}^t\min(j_t,j_b)l_{j_b}  \le {\tilde d} 
\le 2\sum_{b=a}^t\min(j_a,j_b)l_{j_b}.
\end{equation*}
But this is impossible because of $l_{j_t}>0$ and the induction 
assumption, verifying $l_{j_a}\ge 0$.
To show $l_{j_a}\le m_{j_a}$, suppose $l_{j_a}>m_{j_a}$
on the contrary. Then we have 
$K^{(j_a)}_{m_{j_a}} \ge 
2p_{j_a} + 2\sum_{b=a}^t\min(j_a,j_b)l_{j_b}$.
{}From this and (\ref{eq:pjs}), 
$K-{\tilde d}$ can coincide with some $\iota(J)$ only if 
\begin{equation*}
p_{j_a} + 2\sum_{b=a}^t\min(j_a,j_b)l_{j_b} \le {\tilde d} \le
p_{j_s}+2\sum_{b=a}^t\min(j_s,j_b)l_{j_b}.
\end{equation*}
Again it is easy to check that this is impossible by the 
reason similar to the previous inequality, proving $l_{j_a} \le m_{j_a}$.
 
Finally we consider $d$ in the relation 
$\iota(J) +d = 
\sigma_{j_1}^{l_{j_1}} \cdots \sigma_{j_t}^{l_{j_t}}\iota(J')
=\sigma_{j_1}^{l_{j_1}} \cdots \sigma_{j_s}^{l_{j_s}}\iota(J')$.
{}From (\ref{eq:pjt}) and (\ref{eq:pjs}) with $a=1$, we get 
$J^{(j_t)}_{m_{j_t}} +d \ge 
p_{j_t} + 2\sum_{b=1}^t\min(j_t,j_b)l_{j_b}$ and 
$J^{(j_s)}_1 +d\le p_{j_s}+2\sum_{b=1}^t\min(j_s,j_b)l_{j_b}$.
Therefore in order that $(J^{(j)}_i) \in \Rig(m)$ to hold,
$d$ must satisfy  
$0 \le 2\sum_{b=1}^t\min(j_t,j_b)l_{j_b} \le d \le
p_{j_s}+2\sum_{b=1}^t\min(j_s,j_b)l_{j_b}
<p_{j_s}+2\sum_{b=1}^s\min(j_s,j_b)m_{j_b}=L$.
In particular, $d=0$ can happen only if 
$\forall l_{j_b}=0$.
\end{proof}

\begin{proof}[Proof of $\Leftarrow$ in Proposition \ref{pr:migi}]
By Lemma \ref{lem:kagiru}, $J$ and $J'$ are 
connected by the relation 
$(\prod_k\sigma^{l_k}_k)\iota(J')
=\iota(J)+d$ for some $0 \le l_j \le m_j$ and 
$0 \le d < L$.
If $\forall l_k=0$, the assertion follows easily 
from the definition of the KKR bijection $\phi$.
Henceforth we assume that $l_j>0$ for some $j \in H$.
Let $n=(n_j)$ and $l=(l_j)$ 
be the Young diagrams giving the 
decomposition $m=l \cup n$, where
$0 \le n_j \le m_j$ is defined by $n_j=m_j-l_j$.
By the last remark in the proof of Lemma \ref{lem:kagiru}
we have $0 < d < L$.

Define $I^{(j)}=(I^{(j)}_i)_{1 \le i \le l_j}$ and 
$K^{(j)}=(K^{(j)}_i)_{1 \le i \le n_j}$ by
\begin{equation}\label{eq:ijk}
(I^{(j)}_1,\ldots,I^{(j)}_{l_j}) = (J'^{(j)}_1,\ldots, J'^{(j)}_{l_j}),\quad 
(K^{(j)}_1,\ldots,K^{(j)}_{n_j}) = (J^{(j)}_1,\ldots, J^{(j)}_{n_j}).
\end{equation}
Then define further $(I'^{(j)}_i)_{1 \le i \le l_j}$ and 
$(K'^{(j)}_i)_{1 \le i \le n_j}$ 
by (\ref{eq:ipi}) and (\ref{eq:kkp}), respectively.
By construction we have $J = I' \cup K$ and 
$J'=I \cup K'$, where
$(J^{(j)}_{n_j+1},\ldots, J^{(j)}_{m_j}) = 
(I'^{(j)}_1,\ldots, I'^{(j)}_{l_j})$ and 
$(J'^{(j)}_{l_j+1},\ldots, J'^{(j)}_{m_j}) = 
(K'^{(j)}_1,\ldots,K'^{(j)}_{n_j})$. 
The following figure is helpful to grasp these relations.
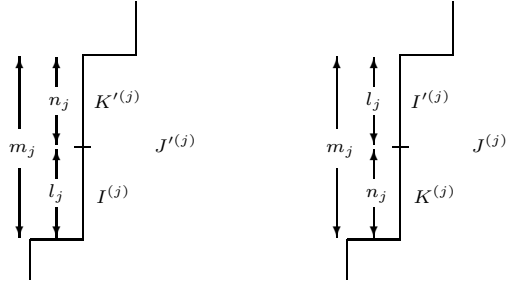
\begin{equation*}
\begin{picture}(60,110)(20,0)
\put(40,85){\line(0,1){20}}
\put(20,85){\line(1,0){20}}
\put(20,15){\line(0,1){70}}
\put(24,65){$\scriptstyle K'^{(j)}$}
\put(10,73){\vector(0,1){11}}
\put(7,66){$\scriptstyle n_j$}
\put(10,62){\vector(0,-1){11}}
\put(47,48){$\scriptstyle J'^{(j)}$}
\put(17,50){\line(1,0){6}}
\put(25,29){$\scriptstyle I^{(j)}$}
\put(10,38){\vector(0,1){11}}
\put(7,31){$\scriptstyle l_j$}
\put(10,27){\vector(0,-1){11}}
\put(0,15){\line(1,0){20}}
\put(0,0){\line(0,1){15}}

\put(-4,57){\vector(0,1){27}}
\put(-8,48){$\scriptstyle m_j$}
\put(-4,42){\vector(0,-1){26}}
\end{picture}
\qquad
\begin{picture}(60,110)(-20,0)
\put(40,85){\line(0,1){20}}
\put(20,85){\line(1,0){20}}
\put(20,15){\line(0,1){70}}
\put(24,65){$\scriptstyle I'^{(j)}$}
\put(10,73){\vector(0,1){11}}
\put(7,66){$\scriptstyle l_j$}
\put(10,62){\vector(0,-1){11}}
\put(47,48){$\scriptstyle J^{(j)}$}
\put(17,50){\line(1,0){6}}
\put(25,29){$\scriptstyle K^{(j)}$}
\put(10,38){\vector(0,1){11}}
\put(7,31){$\scriptstyle n_j$}
\put(10,27){\vector(0,-1){11}}
\put(0,15){\line(1,0){20}}
\put(0,0){\line(0,1){15}}

\put(-4,57){\vector(0,1){27}}
\put(-8,48){$\scriptstyle m_j$}
\put(-4,42){\vector(0,-1){26}}
\end{picture}
\end{equation*}
We claim that $I\in \Rig_{d}(l)$ and $K \in \Rig_{L-d}(n)$.
In fact, 
$0 \le I^{(j)}_1\le \cdots \le I^{(j)}_{l_j}$ is obvious and 
$I^{(j)}_{l_j} = I'^{(j)}_{l_j}-2\sum_k\min(j,k)l_k-p'_j+d
\le d-2\sum_k\min(j,k)l_k$ which is the vacancy number 
for $B^{\otimes d}_1$.
(Here $p'_j = L-2\sum_k\min(j,k)m_k$ is the same as that 
in the proof of Lemmas \ref{lem:waru} and \ref{lem:rotate}.)
Similarly, 
$0 \le K^{(j)}_1\le \cdots \le K^{(j)}_{n_j}$ is obvious and 
$K^{(j)}_{n_j}=K'^{(j)}_{n_j}+
2\sum_k\min(j,k)l_k-d
\le p'_j + 2\sum_k\min(j,k)l_k-d = 
L-d-2\sum_k\min(j,k)n_k$ which is the vacancy number 
for $B^{\otimes L-d}_1$.

Set $q=\phi^{-1}((l,I)) \in B^{\otimes d}_1$ and 
$r= \phi^{-1}((n,K)) \in B^{\otimes L-d}_1$.
Then Lemma \ref{lem:waru} tells that
the highest paths $p_+$ and $p'_+$ in Proposition \ref{pr:migi}
are given by $p_+ = \phi^{-1}((m,J)) = r \otimes q$ and 
$p'_+ = \phi^{-1}((m,J')) = q \otimes r$.
Therefore we obtain $p'_+ = T^{d}_1(p_+)$. 
\end{proof}

The proof of Proposition \ref{pr:migi} is finished.

\section{Proof of Theorem \ref{th:main}}\label{app:main}

For distinction we write the time evolution of the 
angle variable (\ref{eq:tlj}) as $\tau_l$ within this appendix.
Obviously $\tau_l\tau_k = \tau_k\tau_l$ is valid.

It is straightforward to check the commutativity 
of the diagram (\ref{eq:cd}) for $l=1$.
In fact, if $p=T^d_1(p_+)$ for $p_+ \in \Pth_+(m)$ with
$\phi(p_+) = (m,J)$,  we have 
$\Phi(T_1(p))=\Phi(T^{d+1}_1(p_+)) = [\iota(J)+d+1]$.
On the other hand due to $\tau_1(\iota(J))=\iota(J)+1$, 
we obtain $\tau_1(\Phi(p))=\tau_1([\iota(J)+d])=[\iota(J)+d+1]$. 

The commutativity $\tau_1 \Phi = \Phi  T_1$ and 
Lemma \ref{lem:dpm} 
reduce the proof of Theorem \ref{th:main} 
to the highest paths $p_+ \in \Pth_+(m)$.
In fact, the equality $\Phi T_l(p) = \tau_l \Phi(p)$ for general 
path $p=T^d_1(p_+) \in \Pth(m)$ is deduced from 
$\Phi T_l(p_+) = \tau_l \Phi(p_+)$ by multiplying $\tau_1^d$
on the both sides and using the commutativity 
$T_1T_l = T_1T_l$ and $\tau_1\tau_l = \tau_l\tau_1$.

In the remainder of this appendix we fix the 
Young diagram $m = (m_j) \in \M$ and assume that 
$p \in \Pth_+(m) \subset B^{\otimes L}_1$.
($m_j$ is the number 
of the length $j$ rows in $m$ as in the main text.)
Then we determine $J, K, I, d$ and $q$ successively 
as follows:
\begin{align}
\phi(p) &= (m,J), \quad J \in \Rig_L(m),\label{eq:jdef}\\
\tau_l\iota(J) &=K \in \Jb(m),\label{eq:kdef}\\
 [K]& = [\iota(I)+d]  \in \J(m),\quad I \in \Rig_L(m), \; 
d \ge 0,\label{eq:id}\\
\phi^{-1}((m,I)) & = q \in \Pth_+(m),\label{eq:q}
\end{align}
where the definition of $\iota$ 
is available in (\ref{eq:iota}).
The role of these objects will be seen clearly in (\ref{eq:mawaru}).
In (\ref{eq:id}), the choice of $I$ and $d$ is not unique.
However what matters in the following proof is the combination  
$T^d_1(q)$ whose uniqueness has been assured by 
Proposition \ref{pr:migi}.

For $N \in \Z_{\ge 1}$ 
denote by $m^N=(Nm_j)$ the $N$-fold repetition of $m$ and 
define the map 
\begin{equation}\label{eq:pi}
\begin{split}
\pi_N: \qquad\quad\quad\quad
\Jb(m) \qquad \quad
\qquad&\longrightarrow \qquad\quad
\Z^{Nm_{j_1}} \times \cdots \times \Z^{Nm_{j_s}}\\
\left((J^{(j_1)}_i)_{i \in \Z}, \ldots, (J^{(j_s)}_i)_{i \in \Z}\right)
&\mapsto 
\left((J^{(j_1)}_i)_{1 \le i \le Nm_{j_1}}, \ldots, 
(J^{(j_s)}_i)_{1 \le i \le Nm_{j_s}}\right).
\end{split}
\end{equation}
In what follows we set $B = B^{\otimes \lc}_1$ and assume that  
$\lc > (N+1)L$.  
\begin{lemma}\label{lem:repeat}
For $\pi_N(\iota(J)) \in \Rig_{NL}(m^N)$, one has 
\begin{align}
\phi_B^{-1}((m^N,\pi_N(\iota(J)))) &= 
\overbrace{
p \otimes p \otimes p \otimes  \quad \;\;\;\,
\cdots\cdots \;\;\; \quad \otimes p}^{N}
\otimes 1 \otimes 1 \otimes \cdots 
\otimes 1,\label{eq:ppp}\\
\phi_B^{-1}((m^N,\pi_N(K))) &=
p^\ast \otimes \overbrace{
T_l(p) \otimes T_l(p) \otimes \cdots \otimes T_l(p)}^{N-1}
\otimes p^{\ast\!\ast} 
\otimes 1 \otimes \cdots \otimes 1,\label{eq:tptp}
\end{align} 
where $p^\ast \in B^{\otimes L}_1$ 
is defined by (\ref{eq:vt}) and 
$v_l \otimes 1^{\otimes L} \simeq p^{\ast\!\ast}\otimes u_l$.
\end{lemma}
\begin{proof}
The relation (\ref{eq:ppp}) is derived by repeated use of 
Lemma \ref{lem:waru}.
To see (\ref{eq:tptp}), note from (\ref{eq:vt}) that 
\begin{equation*}
u_l \otimes \overbrace{p \otimes \cdots \otimes p}^N
\otimes 1 \otimes \cdots \otimes 1
\simeq 
p^\ast \otimes \overbrace{
T_l(p) \otimes \cdots \otimes T_l(p)}^{N-1}
\otimes p^{\ast\!\ast} 
\otimes 1 \otimes \cdots \otimes 1\otimes u_l.
\end{equation*}
{}From this and (\ref{eq:ppp}), the relation 
(\ref{eq:tptp}) is obtained by 
applying Proposition \ref{pr:hatten} with $p$ and $(m,J)$
replaced with $p^{\otimes N}$ and $(m^N,\pi_N(\iota(J)))$.
\end{proof}

Take the slide 
$\sigma=\sigma_{j_1}^{l_{j_1}} \cdots \sigma_{j_s}^{l_{j_s}} 
\in {\mathcal A}$ 
that achieves (\ref{eq:id}) via $\sigma(K) = \iota(I)+d$.
By the argument similar to the proof of Lemma \ref{lem:kagiru}
(especially after (\ref{eq:pjs})), one can show that 
$l_j$'s can be chosen as $0 \le l_j \le m_j$.
We introduce the Young diagram 
${\tilde m}^N = ({\tilde m}^N_j)$ by
${\tilde m}^N_j = Nm_j-l_j$ and 
let ${\tilde \pi}_N$ be the map (\ref{eq:pi}) with 
$Nm_j$ replaced by ${\tilde m}^N_j$.

\begin{lemma}\label{lem:kimo}
As  $\lc$ and $N$ grow large satisfying $\lc > NL \gg 1$, 
the path $\phi^{-1}_{B}(({\tilde m}^N,{\tilde \pi}_N(\sigma(K))))$ 
takes the form $\xi \otimes \eta$,
where the left part $\xi$ is independent of 
$\lc, N$ and the right part $\eta$ growing with 
$\lc, N$ is the same as the corresponding part in  
$\phi^{-1}_{B}((m^N,\pi_N((K)))$ (\ref{eq:tptp}).
\end{lemma}

\begin{proof}
{}From the definition (\ref{eq:slide}) and 
$\sigma=\sigma_{j_1}^{l_{j_1}} \cdots \sigma_{j_s}^{l_{j_s}}$, 
the rigged configuration $({\tilde m}^N,{\tilde \pi}_N(\sigma(K)))$ 
is obtained from $(m^N, \pi_N(K))$ by 
removing the bottom $l_j$ rows and adding 
$2\sum_{b=1}^s\min(j,j_b)l_{j_b}$ to the rigging for each 
block of length $j$ rows. This operation does not 
change the co-rigging of the remaining part.
Therefore the assertion follows from the KKR algorithm. 
\end{proof}

\begin{proof}[Proof of Theorem \ref{th:main}]
We are to check $\overset{?}{=}$ in the following diagram:
\begin{equation}\label{eq:mawaru}
\begin{picture}(100,70)(0,-15)
\put(42,46){$\scriptstyle \Phi$}
\put(0,40){$p$}
\put(90,40){$[\iota(J)]$}
\put(22,42){\vector(1,0){50}}\put(22,40){\line(0,1){4}}

\put(3,30){\vector(0,-1){25}}\put(1,30){\line(1,0){4}}
\put(98,30){\vector(0,-1){25}}\put(96,30){\line(1,0){4}}

\put(-8,15){$\scriptstyle T_l$} \put(103,15){$\tau_l$}

\put(63,-4){$\scriptstyle \Phi^{-1}$}
\put(80,-8){\vector(-1,0){24}}\put(80,-10){\line(0,1){4}}

\put(90,-10){$[K]=[\iota(I)+d]$}
\put(-10,-10){$T_l(p) \overset{?}{=} T^d_1(q)$}
\end{picture}
\end{equation}
Consider the rigged configuration 
$({\tilde m}^N,{\tilde \pi}_N(\sigma(K)))
=({\tilde m}^N,{\tilde \pi}_N(\iota(I)+d))$.
The latter expression tells that 
this is also obtained from 
$(m^N, \pi_N(\iota(I)+d))$ 
by removing the top $l_j$ rows 
in each block of length $j$ rows.
Thus {}from the KKR algorithm,  
the paths $\phi^{-1}_B(({\tilde m}^N,{\tilde \pi}_N(\iota(I)+d)))$ 
and $\phi^{-1}_{B}((m^N,\pi_N(\iota(I)+d)))$ are 
the same except only some components that are 
located on the right and remain unchanged with growing $\lc, N$.
On the other hand from the definition (\ref{eq:q}) 
and Lemma \ref{lem:repeat} (\ref{eq:ppp}) with 
$(p,J)$ replaced by $(q,I)$, we find 
\begin{equation*}
\phi^{-1}_{B}((m^N,\pi_N(\iota(I)+d)))
= \overbrace{1 \otimes \cdots \otimes 1}^d \otimes 
\overbrace{q \otimes q \otimes \cdots \otimes q}^N 
\otimes 1 \otimes \cdots \otimes 1.
\end{equation*}
Combining the foregoing argument and Lemma \ref{lem:kimo},
we can equate the middle part of this with that in (\ref{eq:tptp}),
leading to $T_l(p) = T^d_1(q)$.
\end{proof}

\end{document}